\documentclass[letterpaper,11pt]{article}
\usepackage[margin=1in]{geometry}
\usepackage{amsmath,amssymb,amsfonts,amsthm}
\usepackage{algorithmicx,algorithm,algpseudocode}
\usepackage{graphicx,booktabs,enumitem,subcaption}
\usepackage[usenames,dvipsnames]{color}
\usepackage{dsfont,bm,titling}
\usepackage[labelfont=bf,skip=0.3em]{caption}
\usepackage{natbib}
\usepackage[T1]{fontenc}
\usepackage[pdftex]{hyperref}
\hypersetup{
  bookmarksopen={true},
  unicode={true},
  pdffitwindow={false},
  pdfstartview={FitH},
  colorlinks={true},
  linkcolor={blue},
  citecolor={blue},
  urlcolor={blue},
  linktocpage={true},
}

\theoremstyle{definition}
\newtheorem{definition}{Definition}[section]
\newtheorem{assumption}{Assumption}[section]
\newtheorem{remark}{Remark}[section]
\newtheorem{example}{Example}[section]

\theoremstyle{plain}
\newtheorem{lemma}{Lemma}[section]
\newtheorem{proposition}{Proposition}[section]
\newtheorem{theorem}{Theorem}[section]
\newtheorem{corollary}{Corollary}[section]

\newcommand{\bX}{\ensuremath{\bm{X}}}
\newcommand{\bY}{\ensuremath{\bm{Y}}}
\newcommand{\bA}{\ensuremath{\bm{A}}}
\newcommand{\bB}{\ensuremath{\bm{B}}}
\newcommand{\bC}{\ensuremath{\bm{C}}}
\newcommand{\bD}{\ensuremath{\bm{D}}}
\newcommand{\bH}{\ensuremath{\bm{H}}}
\newcommand{\bI}{\ensuremath{\bm{I}}}
\newcommand{\bO}{\ensuremath{\bm{O}}}
\newcommand{\bP}{\ensuremath{\bm{P}}}
\newcommand{\bR}{\ensuremath{\bm{R}}}
\newcommand{\bT}{\ensuremath{\bm{T}}}

\newcommand{\ba}{\ensuremath{\bm{a}}}

\newcommand{\bu}{\ensuremath{\bm{u}}}
\newcommand{\bv}{\ensuremath{\bm{v}}}
\newcommand{\bw}{\ensuremath{\bm{w}}}

\newcommand{\cM}{\ensuremath{\mathcal{M}}}
\newcommand{\cA}{\ensuremath{\mathcal{A}}}

\newcommand{\sfY}{\ensuremath{\mathsf{Y}}}
\newcommand{\bsfX}{\ensuremath{\bm{\mathsf{X}}}}

\newcommand{\RR}{\ensuremath{\mathbb{R}}}

\newcommand{\dsone}{\ensuremath{\mathds{1}}}

\newcommand{\bbeta}{\ensuremath{\bm{\beta}}}
\newcommand{\beps}{\ensuremath{\bm{\varepsilon}}}
\newcommand{\bSig}{\ensuremath{\bm{\Sigma}}}
\newcommand{\bpi}{\ensuremath{\bm{\pi}}}
\newcommand{\bmu}{\ensuremath{\bm{\mu}}}

\newcommand{\E}{\ensuremath{\mathsf{E}}}
\newcommand{\Var}{\ensuremath{\mathsf{Var}}}
\newcommand{\Cor}{\ensuremath{\mathsf{Cor}}}
\newcommand{\RSS}{\ensuremath{\mathrm{RSS}}}
\newcommand{\SOIL}{\ensuremath{\mathrm{SOIL}}}

\newcommand{\Norm}{\mathcal{N}}
\newcommand{\Unif}{\textit{Uniform}}

\newcommand{\Bin}{\textit{Binomial}}
\newcommand{\HyperGeom}{\textit{HyperGeom}}

\DeclareMathOperator*{\argmax}{arg\,max}
\DeclareMathOperator*{\argmin}{arg\,min}

\DeclareMathOperator{\GIC}{GIC}
\DeclareMathOperator{\ord}{ord}

\usepackage{tikz}

\usetikzlibrary{arrows, shapes, shapes.multipart, shadows, positioning, fit, backgrounds}
\tikzset{
  cloud/.style={draw, ellipse, fill=red!20, node distance=2.5cm, minimum height=4em, minimum width=8em}, 
  frame/.style={rectangle, draw, text centered, rounded corners}, 
  test/.style={base, diamond, aspect=2, text width=5em}, 
  line/.style={draw, -latex', rounded corners=3mm}, 
}

\algdef{SE}[DOWHILE]{Do}{doWhile}{\algorithmicdo}[1]{\algorithmicwhile\ #1}

\makeatletter
\DeclareRobustCommand\widecheck[1]{{\mathpalette\@widecheck{#1}}}
\def\@widecheck#1#2{%
  \setbox\z@\hbox{\m@th$#1#2$}%
  \setbox\tw@\hbox{\m@th$#1%
    \widehat{%
      \vrule\@width\z@\@height\ht\z@
      \vrule\@height\z@\@width\wd\z@}$}%
  \dp\tw@-\ht\z@
  \@tempdima\ht\z@ \advance\@tempdima2\ht\tw@ \divide\@tempdima\thr@@
  \setbox\tw@\hbox{%
    \raise\@tempdima\hbox{\scalebox{1}[-1]{\lower\@tempdima\box
\tw@}}}{\ooalign{\box\tw@ \cr \box\z@}}}
\makeatother

\let\what=\widehat
\let\wtil=\widetilde

\usepackage[doublespacing]{setspace}
\numberwithin{equation}{section}

\title{Enhancing Multi-model Inference with Natural Selection}

\author{
  Ching-Wei Cheng%
  \thanks{\linespread{1}\selectfont PhD student. Department of Statistics, Purdue University, West Lafayette, IN 47906; e-mail:~\href{mailto:cheng138@purdue.edu}{cheng138@purdue.edu}.}
  \and Guang Cheng%
  \thanks{\linespread{1}\selectfont
    Corresponding Author. Professor, Department of Statistics, Purdue University, IN 47906; e-mail:~\href{mailto:chengg@purdue.edu}{chengg@purdue.edu}.
    Guang Cheng gratefully acknowledges NSF DMS-1712907, DMS-1811812, DMS-1821183, and Office
    of Naval Research, (ONR N00014-18-2759).}}

\date{\vspace*{-2em}}

\begin{document}

\maketitle

\begin{abstract}
\addcontentsline{toc}{section}{Abstract}
Multi-model inference covers a wide range of modern statistical applications such as variable selection, model confidence set, model averaging and variable importance. The performance of multi-model inference depends on the availability of candidate models, whose quality has been rarely studied in literature. In this paper, we study genetic algorithm (GA) in order to obtain high-quality candidate models. Inspired by the process of natural selection, GA performs genetic operations such as selection, crossover and mutation iteratively to update a collection of potential solutions (models) until convergence. The convergence properties are studied based on the Markov chain theory and used to design an adaptive termination criterion that vastly reduces the computational cost. In addition, a new schema theory is established to characterize how the current model set is improved through evolutionary process. Extensive numerical experiments are carried out to verify our theory and demonstrate the empirical power of GA, and new findings are obtained for two real data examples.
\end{abstract}

\noindent{\bf Keywords:}
Convergence analysis; evolvability; genetic algorithm; Markov chain; multi-model inference; schema theory.
\section{Introduction}
\label{sec:intro}

A collection of candidate models serves as a first and important step of multi-model inference, whose spectrum covers variable selection, model confidence set, model averaging and variable importance \citep{BA04:book,A08}. The importance of a candidate model set is highlighted in \citet{LD09}: ``all results of the multi-model analyses are conditional on the (candidate) model set.'' However, in literature, candidate models are either given (e.g., \citet{HLN11,H14}) or generated without any justifications (e.g., \citet{AL14,SOIL}). As far as we know, there is no statistical guarantee on the quality of such candidate models, no matter the parameter dimension is fixed or diverges.

\begin{figure}[!tb]
\centering
\begin{tikzpicture}[font=\small, very thick, node distance=4cm]
\node [rectangle, draw, minimum width=15em, minimum height=10.7em, text centered, rounded corners, fill=teal!20] (pop) {};
\node [above=0em of pop.90] {A {\em population} of size $5$};
\node [rectangle split, rectangle split horizontal, rectangle split parts=10, rectangle split part fill=white, draw, minimum height=1em, below=1em of pop.90] (sol1) {1\nodepart{two}1\nodepart{three}1\nodepart{four}0\nodepart{five}0\nodepart{six}0\nodepart{seven}0\nodepart{eight}0\nodepart{nine}0\nodepart{ten}0};
\node [rectangle split, rectangle split horizontal, rectangle split parts=10, rectangle split part fill=white, draw, minimum height=1em, below=0.5em of sol1] (sol2) {1\nodepart{two}1\nodepart{three}0\nodepart{four}0\nodepart{five}1\nodepart{six}0\nodepart{seven}0\nodepart{eight}0\nodepart{nine}1\nodepart{ten}0};
\node [rectangle split, rectangle split horizontal, rectangle split parts=10, rectangle split part fill=white, draw, minimum height=1em, below=0.5em of sol2] (sol3) {1\nodepart{two}1\nodepart{three}0\nodepart{four}1\nodepart{five}0\nodepart{six}0\nodepart{seven}0\nodepart{eight}0\nodepart{nine}1\nodepart{ten}1};
\node [rectangle split, rectangle split horizontal, rectangle split parts=10, rectangle split part fill=white, draw, minimum height=1em, below=0.5em of sol3] (sol4) {0\nodepart{two}1\nodepart{three}1\nodepart{four}0\nodepart{five}0\nodepart{six}0\nodepart{seven}1\nodepart{eight}0\nodepart{nine}0\nodepart{ten}0};
\node [rectangle split, rectangle split horizontal, rectangle split parts=10, rectangle split part fill=white, draw, minimum height=1em, below=0.5em of sol4] (sol5) {1\nodepart{two}1\nodepart{three}1\nodepart{four}0\nodepart{five}0\nodepart{six}0\nodepart{seven}0\nodepart{eight}0\nodepart{nine}0\nodepart{ten}0};
\node [rectangle split, rectangle split horizontal, rectangle split parts=10, rectangle split part fill=white, draw, minimum height=1em, right=6em of sol1] (sol) {1\nodepart{two}1\nodepart{three}1\nodepart{four}0\nodepart{five}0\nodepart{six}0\nodepart{seven}0\nodepart{eight}0\nodepart{nine}0\nodepart{ten}0};
\node [above=0em of sol] {$1$-st {\em model}/{\em solution} of size $10$};
\node [below=4em of sol.315] (sol_word) {$6$-th {\em variable}/{\em position}/{\em gene} of the $1$-st model};
\draw [draw, latex'-, dashed] (sol) -- (sol1);
\draw [draw, latex'-, dashed] (sol_word.90) -- +(0em, 4em);
\end{tikzpicture}
\caption{\label{fig:term}
An example of GA terminology. Note that the term {\em population} in GA is different from what a ``population'' means in statistics.}
\end{figure}
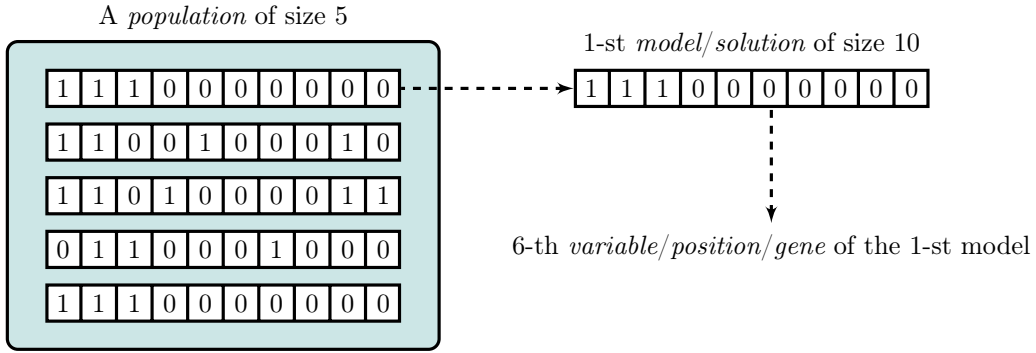

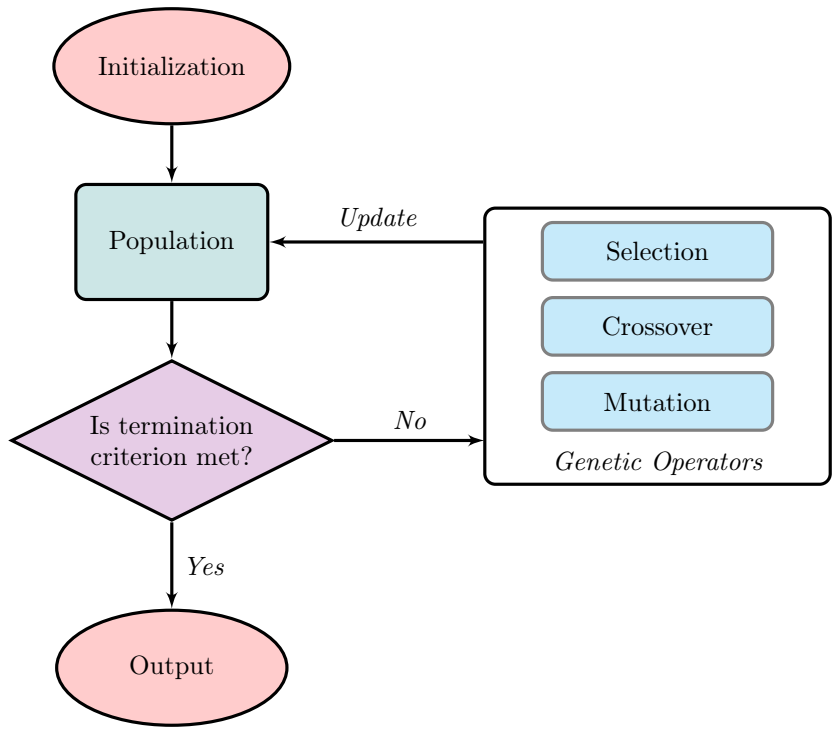
\begin{figure}[!tb]
\centering
\begin{tikzpicture}[font=\small,very thick,node distance=4cm]
\node [frame, text width=6em, minimum height=4em, fill=teal!20] (pop) {Population};
\node [diamond, draw, aspect=2, below=2em of pop, align=center, fill=violet!20] (test) {Is termination\\criterion met?};
\node [rectangle, draw, minimum width=12em, minimum height=9.6em, text centered, rounded corners, above right=-3em and 8em of test] (operators) {};
\node [above=0em of operators.270] {\em Genetic Operators};
\node [frame, draw=gray, minimum width=8em, minimum height=2em, below=0.5em of operators.90, fill=cyan!20] (selection) {Selection};
\node [frame, draw=gray, minimum width=8em, minimum height=2em, below=0.5em of selection, fill=cyan!20] (crossover) {Crossover};
\node [frame, draw=gray, minimum width=8em, minimum height=2em, below=0.5em of crossover, fill=cyan!20] (mutation) {Mutation};
\node [cloud, above=2em of pop] (init) {Initialization};
\node [cloud, below=3em of test] (term) {Output};
\path [line] (init) -- (pop);
\path [line] (pop) -- (test);
\path [line] (test) -- node[right, align=left] {\em Yes} (term);
\path [line] (test.east) -- node[above, align=center] {\em No} +(2.0, 0.0);
\draw [draw, latex'-] (pop) -- node[above, align=center] {\em Update} +(4.14, 0.0);
\end{tikzpicture}
\caption{\label{fig:GA}
A flowchart of a generic GA. It starts with an initial population and is updated with genetic operations until a termination criterion is met.}
\end{figure}

In this paper, we study genetic algorithm (GA, \citet{H75,GA:Intro,GA:Tutorial}) in order to search for high-quality candidate models over the whole model space. GA is a class of iterative algorithms inspired by the process of natural selection, and often used for global optimization or search problems; see Figure~\ref{fig:term}. There are two key elements of GA: a genetic representation of the solution domain, i.e., a binary sequence, and a fitness function to evaluate the candidate solutions such as all kinds of information criteria. A GA begins with an initial population of a given size that is improved through iterative application of genetic operations, such as selection, crossover and mutation, until convergence; see Figure~\ref{fig:GA}.

Specifically, we employ three basic genetic operations, i.e., selection, crossover and mutation, for the GA.
In each generation (the population in each iteration), we adopt elitism and proportional selection so that the fittest model is kept into the next generation, and that fitter models are more likely to be chosen as the ``parent'' models to breed the next generation, respectively.
Uniform crossover is then performed to generate one ``child'' model by recombining the genes from each pair of parent models.
Finally, a mutation operator is applied to randomly alter chosen child genes.
Besides the uniform mutation, we propose a new adaptive mutation strategy using the variable association strength to enhance the variable selection performance.
The genetic operations are iteratively performed until the size of the new generation reaches that of the previous one; see Figure~\ref{fig:GA:evo}.
It is worth noting that the crossover operator generates new models similar to their parents (i.e., local search), while the mutation operator increases the population diversity to prevent GAs from being trapped in local optimum (thus resulting in global search). See Section~\ref{sec:method} for more details.

In theory, we investigate the convergence properties of the GA in Theorem~\ref{thm:mc} based on the Markov chain theory.
A practical consequence is to design an adaptive termination strategy that significantly reduces the computational cost. Furthermore, we prove that a fitter schema (a collection of solutions with specific structures; see Definition~\ref{def:schema}) is more likely to survive and be expanded in the next generation, using the schema theory (Theorem~\ref{thm:schema} and Corollary~\ref{cor:schema:lower}). This implies that the average fitness of the subsequent population gets improved, which entitled the ``survival of the fittest'' phenomenon of the natural selection.

Our results are applied to variable selection and model confidence set (MCS).
In the former, the GA generates a manageable number of models (that is much smaller than all models up to some pre-determined size), over which the true model is found; see Proposition~\ref{prop:vs}. As for the latter, the collected models in the model confidence sets constructed by the GA are shown to be not statistically worse than the true model with a certain level of confidence; see Proposition~\ref{prop:sms}.

As far as we are aware, two other methods can also be used to prepare candidate models:
(i) collecting distinct models on regularization paths of penalized estimation methods (e.g., Lasso \citep{Lasso}, SCAD \citep{SCAD} and MCP \citep{MCP}), called as ``regularization paths (RP)'' method;
(ii) a simulated annealing (SA) algorithm recently proposed by \citet{NR17}.
The former has no rigorous control on the quality of candidate models since model evaluation is not taken into account, and the latter needs a pre-determined model size and an efficiency threshold to filter out bad models. In comparison, the GA uses information-criterion based fitness function to search for good models, and produces models of various sizes.
As a result, the candidate models produced by the GA lead to much improved multi-model inference results, as demonstrated in Sections~\ref{sec:sim} and \ref{sec:realdata}. \citet{AL14} and \citet{SMA} proposed approaches to prepare candidate models that do not work for general multi-model inference applications.
Best subset selection and forward stepwise regression can generate solution paths similar to the Lasso \citep{T15,HTT17}.
However, the former imposes intractable computational burden and the latter lacks of comprehensive theoretical investigation.

Extensive simulation studies are carried out in Section~\ref{sec:sim} to demonstrate the power of the GA in comparison with the RP and the SA in terms of computation time, quality of the candidate model set, and performance of multi-model inference applications. In particular, the GA-best model exhibits the best variable selection performance in terms of the high positive selection rate and low false positive rate. For model averaging and variable importance, the GA results in at least comparable performance to the RP and the SA, but exhibits greater robustness than the SA.
Additionally, the GA is also shown to possess better applicability than the RP in optimal high-dimensional model averaging.

Two real data examples are next carried out to illustrate the practical utility of the GA.
For the riboflavin dataset \citep{BKM14}, the GA-best model finds an informative gene which has not stood out in the literature \citep[e.g.,][]{BKM14,JM14:dLasso,LM15,CLW16,HY17}.
For the residential building dataset \citep{RA16,RA18}, we identify factors, such as preliminary estimated construction cost, duration of construction, and $1$-year delayed land price index and exchange rate, relevant to construction costs. These findings are further confirmed by the variable importance results using the SOIL \citep{SOIL}. Moreover, compared with the aforementioned competing methods, we again find that the GA generates the best candidate model set and results in the best model averaging performance on both datasets.

The rest of this paper is organized as follows.
In Section~\ref{sec:method} we present the GA for global model search, and list several possible ways for improving the implementation.
In Section~\ref{sec:theory} the GA is analyzed using the Markov chain and schema theories.
In Section~\ref{sec:appl} we illustrate how the GA assists multi-model inference tools such as variable selection and model confidence set. Sections~\ref{sec:sim} and ~\ref{sec:realdata} present extensive simulation studies and two real data analysis. In Section~\ref{sec:discussion}, we discuss future works. All proofs are presented in the supplementary materials.

\section{Methodology}
\label{sec:method}

Consider a linear regression model
\begin{align}
\bY=\bX\bbeta^{0}+\beps, \label{eq:LM}
\end{align}
where $\bY=(Y_{1},\ldots,Y_{n})^{\top}$ is the response vector, $\bX=[\bX_{1},\ldots,\bX_{d}]$ is the design matrix with $\bX_{j}$ representing the $j$-th column for $j=1,\ldots,d$, and $\beps=(\varepsilon_{1},\ldots,\varepsilon_{n})^{\top}$ is the noise vector with $\E[\varepsilon_{i}]=0$ and $\Var(\varepsilon_{i})=\sigma^{2}$.
Suppose $\bbeta^{0}=(\beta^{0}_{1},\ldots,\beta^{0}_{d})^{\top}$ is $s$-sparse (i.e., $\|\bbeta^{0}\|_{0}=s$) with $s\ll\min(n,d)$.
Throughout this paper, $s$ and $d$ are allowed to grow with $n$.

\paragraph{Genetic representation for variable selection.}
The genetic representation of a model is defined as a binary sequence of length $d$, say $u=(u_{1},\ldots,u_{d})$, and variable $j$ is said to be active (inactive) if $u_{j}=1$ ($u_{j}=0$).
For example, $u=(1,1,1,0,0,0,0,0,0,0)$ denotes the model with $d=10$ variables but only the first three variables being active. Note that $|u|=\sum_{j=1}^{d}u_{j}$ denote the model size. Denote $\bX_{u}$ as the submatrix of $\bX$ subject to $u$, and  $\cM=\{0,1\}^{d}$ as the model space.

\paragraph{Fitness function.}
Let $\Psi(t)$ denote the $t$-th generation of population, and $\overline{\Psi}(t)=\cup_{t'=0}^{t}\Psi(t')$ the collection of all models that have appeared up to the $t$-th generation.
For any model $u\in\Psi(t)$, the fitness function is then defined as
\begin{align}
f(u)=\begin{cases}-\GIC(u)&\text{if}~|u|<n\\\displaystyle\min_{v\in\overline{\Psi}(t),|v|<n}-\GIC(v)&\text{if}~|u|\geq n\end{cases}, \label{eq:fitness}
\end{align}
where
\begin{align}
\GIC(u)=n\log\what{\sigma}_{u}^{2}+\kappa_{n}|u|, \label{eq:GIC}
\end{align}
is the generalized information criterion (GIC, \citet{N84,S97}) and $\what{\sigma}_{u}^{2}=\bY^{\top}\big[\bI_{n}-\bX_{u}(\bX_{u}^{\top}\bX_{u})^{-1}\bX_{u}^{\top}\big]\bY/n$ is the mean squared error evaluated by the model $u$. GIC covers many types of information criteria (e.g., AIC \citep{AIC} with $\kappa_{n}=2$, BIC \citep{BIC} with $\kappa_{n}=\log n$, modified BIC \citep{mBIC} with $\kappa_{n}=\log\log|u|\log n$ with $d<n$ and extended BIC \citep{eBIC} with $\kappa_{n}\asymp\log n+2\log d$ with $p\geq n$). Since GIC cannot be computed for $|u|\geq n$, we define it as the worst fitness value up to the current generation. The rational is that any model with size larger than $n$ should be unfavorable to all models with size smaller than $n$ given the assumption that $s\ll\min(n,d)$.
This definition warrants an unconstrained optimization, which is convenient for subsequent theoretical analysis. This is different from other ways to deal with the ``infeasible solutions'' in the GA literature, e.g., the ``death penalty'' in \citet{CYPI16} and \citet{ZOMKO14}, which lead to constrained optimization.

\subsection{A Genetic Algorithm for Candidate Model Search}

We propose a genetic algorithm to search for good candidate models in Algorithm~\ref{alg:GA}.
Specifically, we use the RP method to generate an initial population, and then adopt proportional selection, uniform crossover and mutation operators to constitute the evolutionary process. 
Besides uniform mutation, we propose another mutation strategy based on the strength of variable association for improving empirical performances. An adaptive termination strategy is also proposed to enhance the computational efficiency. See Algorithm~\ref{alg:GA} for the overview of the GA.

\begin{algorithm}[!tb]
\caption{A Genetic Algorithm for Model Search}
\label{alg:GA}
\begin{algorithmic}[1]
\Require{Population size $K$ and mutation rate $\pi_{m}$}
\State Generate initial population $\Psi(0)=\big\{u^{1}(0),\ldots,u^{K}(0)\big\}$
\State $t\leftarrow 0$
\State $\it{Converge}\leftarrow\it{False}$
\Do
\State $t\leftarrow t+1$
\State {\em (Fitness evaluation)} Compute fitness values $f\big(u^{k}(t-1)\big),k=1,\ldots,K$
\State {\em (Elitism selection)} Set $u^{1}(t)=\argmax_{u\in\Psi(t-1)}f(u)$
\For{$k=2,\ldots,K$}
\State {\em (Proportional selection)} Randomly select two models from $\Psi(t-1)$ using $w_{k}$ in \eqref{eq:modelweight}
\State {\em (Uniform crossover)} Breed a child model using \eqref{eq:UX}
\State {\em (Mutation)} Mutate the child genes using \eqref{eq:mut:unif} or \eqref{eq:mut:varimp}
\EndFor
\State Set $\Psi(t)=\big\{u^{1}(t),\ldots,u^{K}(t)\big\}$
\If{Convergence criterion \eqref{eq:GA:terminate} is met}
\State $T\leftarrow t$
\State $\it{Converge}\leftarrow\it{True}$
\EndIf
\doWhile{$\it{Converge}$ is $\it{False}$}
\State \Return{$\Psi(T)=\big\{u^{1}(T),\ldots,u^{K}(T)\big\}$}
\end{algorithmic}
\end{algorithm}

\begin{figure}[!tb]
\centering
\includegraphics[scale=0.7]{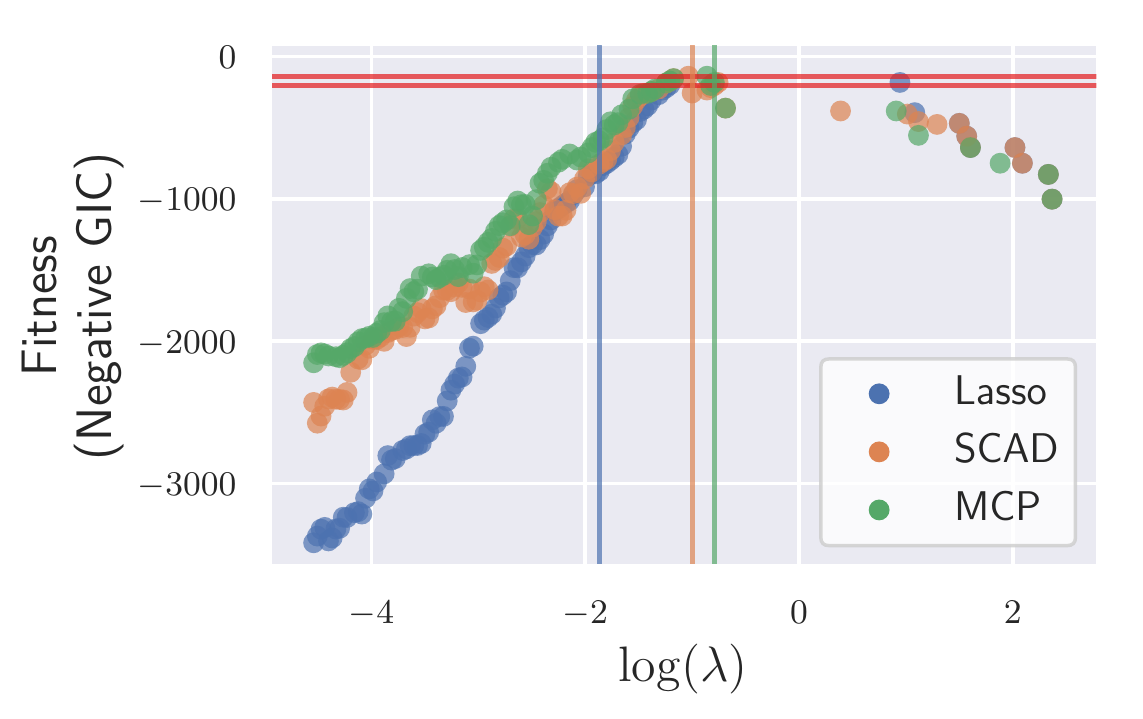}
\includegraphics[scale=0.7]{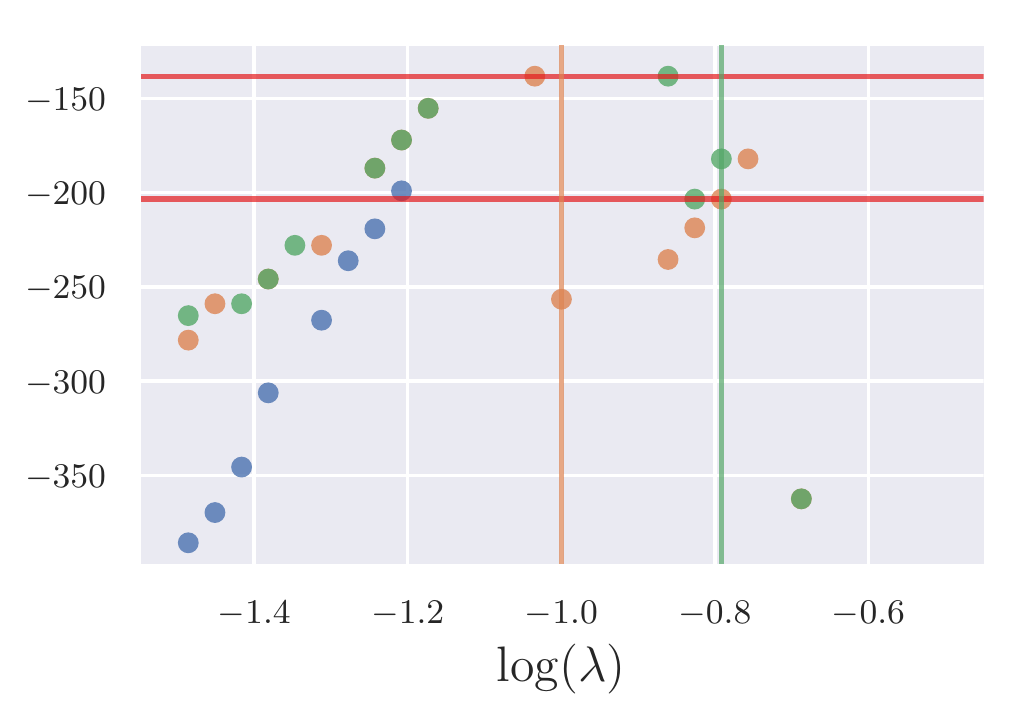}
\caption{\label{fig:fitness:RP}
Fitness values of the models obtained from the regularization paths of Lasso, SCAD and MCP.
The two red horizontal lines indicate the best and worst fitness values of the GA models, and the vertical lines locate the best $\lambda$ selected by $10$-fold cross-validation. The right panel is a zoomed view of the left panel around the selected $\lambda$. Among the $284$ RP models, only $10$ of them have fitness values not smaller than the GA-worst model. The result is obtained from the first dataset under the simulation Case 1 with $(n,d,s,\rho)=(200,400,6,0.5)$.}
\end{figure}

\paragraph{Initialization:}
The initial population $\Psi(0)=\big\{u^{1}(0),\ldots,u^{K}(0)\big\}$ only has very minimal requirement as follows: (i) $K\geq 2$ and (ii) $|u^{k}(0)|<n$ for some $k=1,\ldots,K$ (i.e., at least one model with commutable GIC). 
The condition (i) allows the GA to explore through the model space $\cM$; see Section~\ref{sec:mc}, and (ii) ensures $f\big(u^{k}(0)\big),k=1,\ldots,K$, are all available. The choice of $K$ will be discussed in Section~\ref{sec:popsize}. For fast convergence of the GA, we recommend the RP method to generate initial population. Please see Figure~\ref{fig:fitness:RP} for how models produced by RP are improved by the GA in terms of GIC.

\begin{figure}[!tb]
\centering
\begin{tikzpicture}[font=\small, node distance=8em, very thick]
\node [frame, text width=5em, label={$\Psi(t-1)$}, fill=teal!20, minimum height=10em] (pop_old) {$u^{1}(t-1)$\\[1em]$u^{2}(t-1)$\\[1em]$\vdots$\\[1em]$u^{K}(t-1)$};
\node [frame, fill=cyan!20, text width=8em, minimum height=2em, right=18em of pop_old.50] (best) {$u^{1}(t)=u^{*}(t-1)$};
\node [frame, text width=2.5em, fill=cyan!20, right=9em of pop_old.320, minimum height=4em] (parents) {$u^{p1,k}$\\[0.5em]$u^{p2,k}$};
\node [frame, text width=2.5em, minimum height=2em, fill=cyan!20, right=5em of parents] (child0) {$u^{c,k}$};
\node [frame, text width=2.5em, minimum height=2em, fill=cyan!20, right=5em of child0] (child) {$u^{k}(t)$};
\node [frame, text width=5em, label={$\Psi(t)$}, right=31em of pop_old, fill=teal!20, minimum height=10em] (pop_new) {$u^{1}(t)$\\[1em]$u^{2}(t)$\\[1em]$\vdots$\\[1em]$u^{K}(t)$};
\path [line] (pop_old.50) -- node[above, align=center, pos=0.5] {\em Elitism Selection} (best);
\path [line] (best.east) -- node {} (pop_new.130);
\path [line] (pop_old.320) -- node[above, align=center, pos=0.5] {\em Proportional\\\em Selection} (parents.west);
\path [line] (parents.east) -- node[above, align=center, pos=0.5] {\em Uniform\\\em Crossover} (child0.west);
\path [line] (child0.east) -- node[above, align=center, pos=0.5] {\em Mutation} (child.west);
\path [line] (child.east) -- node {} (pop_new.220);
\path (parents.west |- parents.north)+(-3.0,0.5) node (a1) {};
\path (child.east |- child.south)+(+0.5,-1.0) node (a2) {};
\path[rounded corners, draw=black!50, dashed] (a1) rectangle (a2);
\path (a1.east |- a1.south)+(1.0,-2.5) node (u1) [above, text width=6em, text centered] {$k=2,\ldots,K$};
\end{tikzpicture}
\caption{\label{fig:GA:evo}
Illustration of the evolution process of the GA.
$u^{*}(t-1)$ denotes the best model in $\Psi(t-1)$.
The candidate pool of the proportional selection is the entire $\Psi(t-1)$, which still includes $u^{*}(t-1)$.
For each $k=2,\ldots,K$, a pair of parent models are selected according to the probability $w_k$ and one child model $u^{c,k}$ is generated through uniform crossover \eqref{eq:UX}.
Finally, $u^{c,k}$ is processed by the mutation \eqref{eq:mut:unif} or \eqref{eq:mut:varimp} to produce $u^{k}(t)$.}
\end{figure}
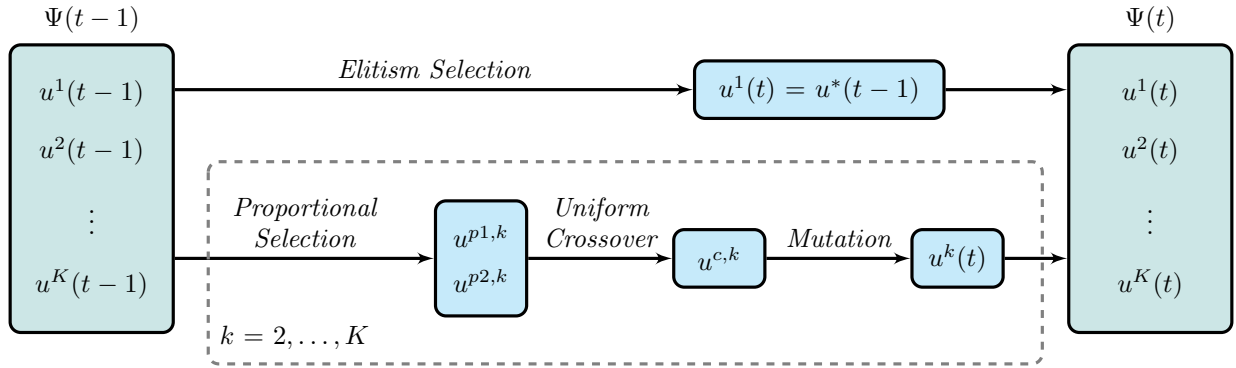

Given the $(t-1)$-th generation $\Psi(t-1)=\big\{u^{1}(t-1),\ldots,u^{K}(t-1)\big\}$, the GA produces the next generation $\Psi(t)=\big\{u^{1}(t),\ldots,u^{K}(t)\big\}$ through proportional selection, uniform crossover and mutation operations. See Figure~\ref{fig:GA:evo} to visualize the evolution process. In what follows, we give details for each step in our main algorithm. 

In the {\bf elitism selection} step, we choose $u^{*}(t-1):=\argmax_{u\in\Psi(t-1)}f(u)$, i.e., the best model in $\Psi(t-1)$ is kept into $\Psi(t)$, and define it as $u^{1}(t)$ for simplicity. The {\bf proportional selection} step chooses parent models from $\Psi(t-1)$ (including $u^{*}(t-1)$) based on the exponentially scaled fitness as follows. Define the fitness $f_{k}=f\big(u^{k}(t-1)\big)$ according to \eqref{eq:fitness}. For $k=1,\ldots,K$, first compute the weight $w_{k}$ for $u^{k}(t-1)$ as
\begin{align}
w_{k}=\frac{\exp\left(f_{k}/2\right)}{\sum_{l=1}^{K}\exp\left(f_{l}/2\right)},\quad k=1,\ldots,K. \label{eq:modelweight}
\end{align}
Then $(K-1)$ pairs of models are randomly selected with replacement from $\Psi(t-1)$, where the probability of selecting $u^{k}(t-1)$ is $w_k$. Note that the exponentially scaled information criteria are often used for model weighting in multi-model inference \citep[e.g.,][]{BA04:book,AICw,BMA}.

Each pair of parent models produces a child model by performing {\bf uniform crossover} with equal mixing rate (i.e., each child position has equal chance to be passed from the two parents).
That is, let $u^{p1,k}=(u^{p1,k}_{1},\ldots,u^{p1,k}_{d})$ and $u^{p2,k}=(u^{p2,k}_{1},\ldots,u^{p2,k}_{d})$ be the chosen parent models, and then the genes in the child model $u^{c,k}=(u^{c,k}_{1},\ldots,u^{c,k}_{d})$ is determined by
\begin{align}
u^{c,k}_{j}=\begin{cases}u^{p1,k}_{j}&\text{with probability}~1/2\\u^{p2,k}_{j}&\text{otherwise}\end{cases},\quad j=1,\ldots,d. \label{eq:UX}
\end{align}

In the last step, we apply {\bf mutation} to the child model  $u^{c,k}$. Given a mutation probability $\pi_{m}$ (usually low, such as $\pi_{m}=0.01$ or $1/d$), we consider the following two mutation schemes.
Denote by $u^{k}(t)=\big(u^{1}_{d}(t),\ldots,u^{k}_{d}(t)\big)$ the resulting model after mutation being applied to $u^{c,k}$.
\begin{itemize}
\item {\bf Uniform mutation:}
Genes in $u^{c,k}$ are randomly flipped with probability $\pi_{m}$, i.e.,
\begin{align}
u^{k}_{j}(t)=\begin{cases}1-u^{c,k}_{j}&\text{with probability}~\pi_{m}\\u^{c,k}_{j}&\text{otherwise}\end{cases},\quad j=1,\ldots,d. \label{eq:mut:unif}
\end{align}
\item {\bf Adaptive mutation:}
We propose a data-dependent mutation operator based on the variable association measures $\gamma_{j}$. For example, $\gamma_j$ can be either the marginal correlation learning $\big|\what{\Cor}(\bX_{j},\bY)\big|$ \citep{SIS} or the high-dimensional ordinary least-squares projection $\big|\bX_{j}(\bX\bX^{\top})^{-1}\bY\big|$ \citep[][available only for $d\geq n$]{HOLP}.
Let $V^{k}_{+}=\{j:u^{c,k}_{j}=1\}$ and $V^{k}_{-}=\{j:u^{c,k}_{j}=0\}$.
Define the mutation probability for the $u^{c,k}_{j}$ as
\begin{align*}
\bar{\pi}^{k}_{m,j}=\begin{cases}\dfrac{\gamma_{j}^{-1}}{\sum_{l\in V^{k}_{+}}\gamma_{l}^{-1}}|V^{k}_{+}|\pi_{m}&\text{if}~j\in V^{k}_{+}\\\dfrac{\gamma_{j}}{\sum_{l\in V^{k}_{-}}\gamma_{l}}|V^{k}_{-}|\pi_{m}&\text{if}~j\in V^{k}_{-}\end{cases}.
\end{align*}
Then the proposed mutation operation is performed by
\begin{align}
u^{k}_{j}(t)=\begin{cases}1-u^{c,k}_{j}&\text{with probability}~\bar{\pi}^{k}_{m,j}\\u^{c,k}_{j}&\text{otherwise}\end{cases},\quad j=1,\ldots,d. \label{eq:mut:varimp}
\end{align}
By defining $\bar{\pi}^{k}_{m,j}$ this way, unimportant active variables are more likely to be deactivated, and important inactive variables are more likely to be activated.
Also, it can be easily seen that this mutation operation results in the same expected number of deactivated and activated genes as those of uniform mutation operation. As far as we are aware, this is the first data dependent mutation method in the GA literature.
\end{itemize}
In numerical experiments, we note that the adaptive mutation performs slightly better than the uniform mutation. For space constraint, we just focus on the adaptive mutation with $\pi_{m}=1/d$.

As for {\bf termination}, we propose an adaptive criterion by testing whether the average fitness becomes stabilized; see Section~\ref{sec:adaterm} for more details. This is very different from the user specified criteria used in GA literature such as the largest number of generations \citep[e.g.,][]{MMM16} or the minimal change of the best solution \citep[e.g.,][]{ACLZ14}.

\begin{figure}[!tb]
\centering
\includegraphics[width=0.95\textwidth]{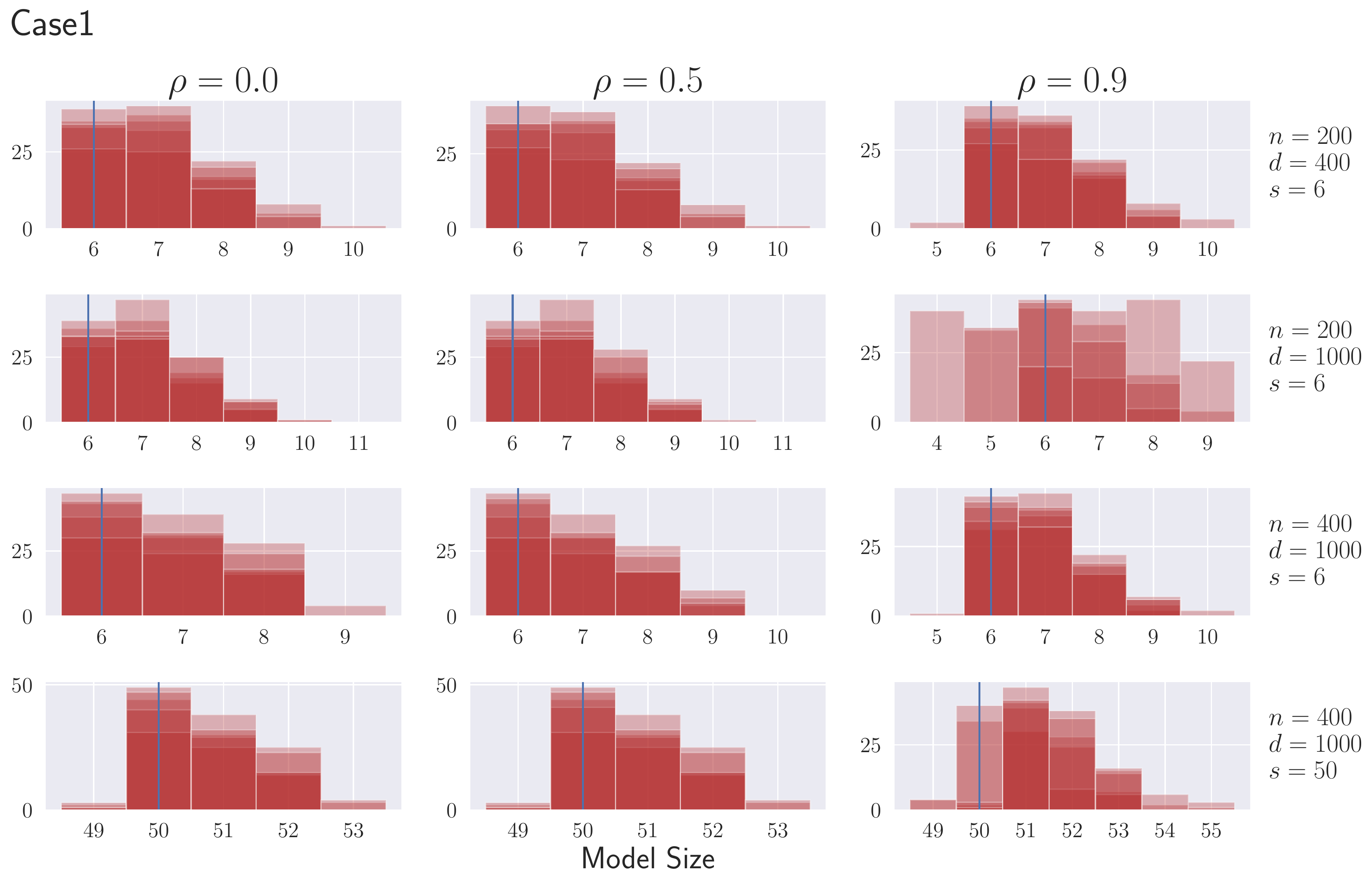}
\caption{\label{fig:sparse}
The overlapped distributions of the sizes of the final candidate models collected by the GA, with the blue vertical line indicating the true model size.
The results are obtained under the simulation Case 1 (see Section~\ref{sec:sim} for more details) and other cases exhibit similar patterns.}
\end{figure}

\begin{remark}\label{rmk:sparse}
We note that the models collected by the GA are in nature sparse since their sizes are around the true model size $s$; see Figure~\ref{fig:sparse} for example.
This empirically appealing feature allows us to construct GA-based sparse model confidence sets in the later sections. 
\end{remark}

\subsection{Computational Considerations}

Computational concern has been the major critiques that prevent GAs from being popular over other optimization methods such as gradient descent in machine learning and statistical communities.
In our experience, the most computational cost is taken by the calculation of the fitness evaluation, which could be alleviated by reducing the population size (Section~\ref{sec:popsize}) and the number of generations (Section~\ref{sec:adaterm}).

\subsubsection{Population Sizing}
\label{sec:popsize}

The population size $K$ plays an important role in GAs.
It is obvious that larger population makes GAs computationally more expensive.
On the other hand, empirical results indicate that small population would jeopardize performance \citep[e.g.,][]{KK06,PS06,LL05}.
We found that the minimum population size suggested in \citet{R93} makes a good balance.
The idea is to have a population such that every possible solution in the search space should be reachable from an randomly generated initial population by crossover only.
In binary gene coding cases, it means that the solutions in the initial population cannot be all $0$ or $1$ for any position.
For any $K$, the probability of such an event can be found by
\begin{align*}
P^{*}&=(1-1/2^{K-1})^{d}=\exp\big[d\log(1-1/2^{K-1})\big]\approx\exp(-d/2^{K-1}).
\end{align*}
Accordingly, for every given $P^{*}$, we can calculate the minimum population size
\begin{align*}
K^{*}\approx\big\lceil 1+\log(-d/\log P^{*})/\log 2\big\rceil,
\end{align*}
where $\lceil a\rceil$ is the smallest integer larger than $a\in\RR$.
For example, a population of size $K=25$ is enough to ensure that the required probability exceeds $99.99{\%}$ when $d=1{,}000$.

In our implementation, we conservatively use
\begin{align*}
K=4\big\lceil 1+\log(-d/\log P^{*})/\log 2\big\rceil
\end{align*}
with $P^{*}=0.9999$, to specify the population size according to model dimension.

\subsubsection{Adaptive Termination}
\label{sec:adaterm}

To adaptively terminate, we perform an independent two-sample $t$-test on whether the average fitness of $\Psi(t)$ and $\Psi(t-10)$ are the same at a significance level $0.05$:
\begin{align}\label{eq:GA:terminate}
H^{t}_{0}:~\bar{f}\big(\Psi(t)\big)=\bar{f}\big(\Psi(t-10)\big)\quad\text{v.s.}\quad H^{t}_{1}:~\bar{f}\big(\Psi(t)\big)\not=\bar{f}\big(\Psi(t-10)\big),
\end{align}
where $\bar{f}\big(\Psi(t)\big)$ is the average fitness of the $t$-th generation.
The $T$-th generation is set to be the final generation if $T$ is the smallest $t\geq 10$ such that $H^{t}_{0}$ is rejected. The generation gap $10$ is meant to weaken the correlation between the two generations being tested.
Note that the GA can be regarded as a Markov chain (see Section~\ref{sec:mc} for details) and therefore there exists dependence among generations. Hence, it is not appropriate to perform two-sample $t$-test of the average fitness from two consecutive generations.

\begin{remark}
This termination criterion is constructed based on the limiting distribution derived for the associated Markov chain (see the discussion below Theorem~\ref{thm:mc} for more details) and results in huge computational efficiency. In the literature, the GA iteration is often terminated at a fixed, predetermined number of generations, say $T_{\max}$, which is usually large such as $50,100$ or even larger \citep[e.g.][]{H17,JIM18}. Our termination criterion, on the other hand, entitles a scientific check for the convergence.
With the RP used for generating the initial population, the GA enters the stationary distribution (as the average fitness is tested to be stabilized) in just a few generations (say, around $20$ generations). In addition, we note that the computational cost incurred by the independent two-sample $t$-test \eqref{eq:GA:terminate} is negligible, as the fitness values are computed as the models are generated.
\end{remark}

\section{Theoretical Properties}
\label{sec:theory}

In this section, we study the theoretical properties of the GA, which belongs to the so-called canonical GA (CGA) family\footnote{CGAs are also called as simple or standard GAs in the literature. CGA uses binary sequence for solution representation, and updates a fixed-sized population via selection, crossover and mutation operators.} \citep{H75}. The proposed GA is a CGA that specifically employs elitism and proportional selection, uniform crossover and uniform or adaptive mutation as described in Section~\ref{sec:method}. We first investigate the convergence properties for a general CGA family based on Markov chain theory, i.e., Theorem~\ref{thm:mc}. Furthermore, Theorem~\ref{thm:mcs} presents a brand new theoretical framework to construct MCSs for the globally best model. We next establish a new schema theory (Theorem~\ref{thm:schema} and Corollary~\ref{cor:schema:lower}) to elicit the evolutionary mechanism for the GA. It is worthy noting that the theoretical results established in this section apply to the general CGA framework and hence not restricted to the specific variable selection problem.

\subsection{Convergence Analysis}
\label{sec:mc}

In this section, we show that the Markov chain associated with a CGA class has a unique stationary distribution from which the asymptotic inclusion of the globally best model, i.e., global convergence, can be deduced. Such a result justifies the adaptive termination rule in Section~\ref{sec:adaterm}, and can be used to reduce the search space for variable selection problems; see Proposition~\ref{prop:vs}. 
Note that the theoretical results obtained in this section hold for any finite sample size.

Recall that $\Psi(t)=\big\{u^{1}(t),\ldots,u^{K}(t)\big\}$ represents the $t$-th generation of the population, and denote by $\big\{\Psi(t)\big\}_{t\geq 0}$ the associated Markov chain with values on the finite state (population) space $\cM^{K}$. The corresponding transition matrix is defined as $\bP=\big[P_{\bu\bv}\big]_{\bu,\bv\in\cM^{K}}$, where $P_{\bu\bv}=P\big(\Psi(t+1)=\bv\big|\Psi(t)=\bu\big)$. We need the following definitions for our subsequent analysis.
\begin{definition}
A square matrix $\bA=[a_{ij}]\in\RR^{K\times K}$ is said to be {\em non-negative (positive)} if $a_{ij}\geq 0$ ($a_{ij}>0$) for all $i,j\in\{1,\ldots,d\}$.
A non-negative square matrix $\bA$ is said to be
\begin{enumerate}[label=(\alph*)]
\item {\em primitive} if there exists a positive integer $k$ such that $\bA^{k}$ is positive;
\item {\em reducible} if there exists two square matrices $\bA_{11}$ and $\bA_{22}$ and a matrix $\bA_{21}$ with suitable dimensions such that $\bA$ can be expressed as the form
\begin{align*}
\bA=\begin{bmatrix}\bA_{11}&\bO\\\bA_{21}&\bA_{22}\end{bmatrix},
\end{align*}
where $\bO$ denotes a zero matrix with suitable dimensions, by applying the same permutations to rows and columns;
\item {\em irreducible} if it is not reducible;
\item {\em stochastic} if $\sum_{j=1}^{K}a_{ij}=1$ for all $i=1,\ldots,K$.
\end{enumerate}
\end{definition}

Let
\begin{align*}
u^{*}:=\argmax_{u\in\cM}f(u)
\end{align*}
denote the best model in $\cM$ and suppose it is unique, i.e., $f(u^{*})>f(u)$ for all $u\in\cM-\{u^{*}\}$.
Moreover, denote the collection of states that contains $u^{*}$ by
\begin{align}\label{eq:Mmax}
\cM_{\max}=\Big\{\bu=\big\{u^{1},\ldots,u^{K}\big\}\in\cM^{K}:u^{*}\in\bu\Big\}.
\end{align}

The following theorem describes two important convergence properties.
\begin{theorem}\label{thm:mc}
Let $\bP$ denote the transition probability matrix of the Markov chain associated with a CGA with elitism selection, population size $K\geq 2$ and mutation rate $\pi_{m}\in(0,1)$.
\begin{enumerate}[label=(\alph*)]
\item There exists a unique stationary distribution $\bpi=\big(\pi(\bu):\bu\in\cM^{K}\big)^{\top}$ that satisfies $\bpi^{\top}=\bpi^{\top}\bP$ and $\pi(\bu)=\lim_{t\to\infty}P\big(\Psi(t)=\bu\big)$ with $\pi(\bu)>0$ for $\bu\in\cM_{\max}$ and $\pi(\bu)=0$ for $\bu\not\in\cM_{\max}$.
\item (Theorem~6 of \citet{R94})
We have
\begin{align}
\lim_{t\to\infty}P\big(u^{*}\in\Psi(t)\big)=1. \label{eq:conv:global}
\end{align}
\end{enumerate}
\end{theorem}

As far as we are aware, the existence of the stationary distribution stated in Theorem~\ref{thm:mc} (a) for CGAs {\em with elitism selection} is new, even though similar results for non-elitist CGAs has been presented for over decades \citep[e.g.,][]{R94,DGMP10}.
This is in contrast with the GA literature that typically concerns global convergence \citep[e.g.,][]{R94,A98,DGMP10} rather than the stationary distribution.
As for Theorem~\ref{thm:mc} (b), the elitism selection is a necessary condition \citep{R94,A98} and it is different from the path-consistency property of non-convex penalization approaches (e.g., \citet{KK12,CCCP}) in that the former captures the {\em best} model for any sample size $n$ as $t\to\infty$ and the latter targets at the {\em true} model as $n\to\infty$.
Later, Theorem~\ref{thm:mc} (b) is extended to a selection consistency result as $n\to\infty$; see Proposition~\ref{prop:vs}.

Part (a) of Theorem~\ref{thm:mc} has the following implication.
Recall that $\bar{f}(\bu)$ is the average fitness of any population $\bu$, and thus we have
\begin{align*}
\lim_{t\to\infty}\E\big[\bar{f}\big(\Psi(t)\big)\big]=\sum_{\bu\in\cM^{K}}\pi(\bu)\bar{f}(\bu),
\end{align*}
which is a constant given data $(\bX,\bY)$. This indicates that the average fitness oscillates around a constant in the long run, as $\Psi(t)$ becomes stabilized (i.e., the associated Markov chain converges). This justifies the termination check in \eqref{eq:GA:terminate}.

\begin{remark}
It is worth noting that Theorem~\ref{thm:mc} does not only apply to the GA but also any CGA with elitism selection.
The key reason is that the child solutions generated through selection, crossover and mutation operators always remain in the search space for unconstrained optimization or search problems.
Accordingly, instead of the ones mentioned in Section~\ref{sec:method}, Theorem~\ref{thm:mc} still holds for any other selection, crossover and mutation operations (e.g., rank-based or tournament selection \citep{SPM15} and the newly proposed crossover and mutation operations developed in \citet{HA18a} and \citet{HA18b}, respectively).
\end{remark}

In contrast to the asymptotic result in Theorem~\ref{thm:mc} (b) as $t\to\infty$, it is also of practical relevance to construct a $100(1-\alpha){\%}$ MCS that covers the best model $u^{*}$ after a finite number of generations.
A particularly appealing feature is that every model in this set is sparse.
This is conceptually different from the MCS constructed based on the debiased principle \citep{VGBRD14:dLasso,ZZ14:LDP,JM14a:dLasso}, which mostly contains dense models.

\begin{theorem}\label{thm:mcs}
Let $\Psi(t)$ denote the $t$-th population of a CGA with elitism selection, $K\geq 2$ and $\pi_{m}\in(0,1)$.
Then for any $\alpha\in(0,1)$ there exists a positive integer $T_{\alpha}$ such that
\begin{align}
P\big(u^{*}\in\Psi(t)\big)\geq 1-\alpha \label{eq:mcs}
\end{align}
for any $t\geq T_{\alpha}$.
\end{theorem}

The proof of Theorem~\ref{thm:mcs} implies the global convergence property described in Theorem~\ref{thm:mc} (b) by letting $\alpha=0$ and thus $T_{0}=\infty$. From the proof of Theorem~\ref{thm:mcs}, we note that obtaining the value of $T_{\alpha}$ requires the knowledge of the constant $\xi$ as defined in \eqref{eq:mcs:xi}, which is often unknown. By definition, $\xi$ can be obtained by estimating the submatrix $\bR$ in the transition matrix $\bP$.
That is, $\xi$ is the smallest of the row sums of $\bR$. Since $\bR$ has $\big|\cM^{K}\big|-\big|\cM_{\max}\big|=\binom{K+2^{d}-1}{K}-\binom{K+2^{d}-2}{K-1}$ rows and $\big|\cM_{\max}\big|=\binom{K+2^{d}-2}{K}$ columns, the size of $\bR$ is massive. For instance, when $(K,d)=(10,5)$, there are about $3\times 10^{17}$ elements in $\bR$. Albeit \citet{V93} provides a useful formula to compute the elements in the $\bP$, the computational cost is too large to be carried out in practice. Hence, we leave an accurate estimation or approximation of $\xi$ to future study.

\subsection{Evolvability Analysis}
\label{sec:schema}

In this section, we establish a schema theorem to study the evolution process of the GA. Specifically, it is proven that the average fitness gets improved over generations.
To the best of our knowledge, we are the first to develop a schema theorem for GAs with proportional selection, uniform crossover and uniform mutation at the same time in the GA literature.
The most closely related schema theorems are provided by \citet{P01a,P01b} for GAs with proportional selection and one-point crossover, and by \citet{MWC04} for GAs with uniform crossover alone.

In the following we give the definition of a schema with general GA terminology (i.e., using ``solutions'' instead of ``models''), followed by an example as illustration.

\begin{definition}\label{def:schema}
A {\em schema} $H=(H_{1},\ldots,H_{d})\in\{0,1,\ast\}^{d}$ is a ternary sequence of length $d$, where the ``$\ast$'' is a wildcard symbol, meaning that we do not care whether it is $0$ or $1$.
The indices where the schema has a $0$ or $1$ are called the {\em fixed positions}.
We say a solution $u=(u_{1},\ldots,u_{d})$ matches $H$ if all fixed positions of $H$ are the same as the corresponding positions in $u$.
The {\em order of a schema $H$}, denoted by $\ord(H)$, is defined by the number of fixed positions in $H$.
Moreover, by adopting the notations used in the order theory (e.g., \citet{FW17}), for any schema $H$ we define the {\em expansion operator ${\uparrow}(H)$} to map $H$ to the set of all possible solutions that match $H$, i.e.,
\begin{align*}
{\uparrow}(H)=\big\{u\in\cM:u_{j}=H_{j}~\text{or}~H_{j}=\ast~\text{for each}~j=1,\ldots,p\big\}.
\end{align*}
\end{definition}

\begin{example}
Suppose a schema $H=(1,0,\ast,0,\ast)$. In this case, $\ord(H)=3$, and   ${\uparrow}(H)=\big\{(1,0,0,0,0),(1,0,0,0,1),(1,0,1,0,0),(1,0,1,0,1)\big\}$.
\end{example}

Let $m(H,t)$ denote the number of solutions that match a schema $H$ in the $t$-th generation, and $\alpha(H,t)$ the probability that the schema $H$ survives or is created after the $t$-th generation.
\citet{PLO98} noted that $m(H,t+1)$ follows a binomial distribution with the number of trials $K$ and success probability $\alpha(H,t)$, i.e., ($K$ is the population size)
\begin{align}
m(H,t+1)\sim\Bin\big(K,\alpha(H,t)\big). \label{eq:schema:prob}
\end{align}
Hence, we have $\E\big[m(H,t+1)\big]=K\alpha(H,t)$.
Accordingly, higher $\alpha(H,t)$ leads to higher $\E\big[m(H,t+1)\big]$ and thus tends to result in more solutions of $H$ in the next generation.
Since the population size is fixed, more solutions of a fitter schema imply higher average fitness in the subsequent generation. Hence, we will show that $\alpha(H^{1},t)$ is larger than $\alpha(H^{2},t)$ if the average fitness of $H^{1}$ is larger than that of $H^{2}$.

To prove the above result, we need to define the following different notions of Hamming distance. The first concerns two models $u=(u_{1},\ldots,u_{p})$ and $v=(v_{1},\ldots,v_{p})$, i.e., $\delta(u,v)=\sum_{j=1}^{p}\dsone(u_{j}\not=v_{j})$, while the second type of Hamming distance is between a model and a schema $H$ on the fixed positions: $\delta(u,H)=\sum_{j:H_{j}\not=\ast}\dsone(u_{j}\not=H_{j})$. The last one is Hamming distance between models $u$ and $v$ with respect to the fixed positions of any schema $H$: $\delta_{H}(u,v)=\sum_{j:H_{j}\not=\ast}\dsone(u_{j}\not=v_{j})$.

We are now ready to characterize $\alpha(H,t)$ explicitly for the GA with uniform mutation. Recall from \eqref{eq:modelweight} that $w_{k}$ denotes the probability that model $u^{k}$ is selected as a parent model.

\begin{theorem}\label{thm:schema}
Given the $t$-th generation $\Psi(t)=\{u^{1},\ldots,u^{K}\}$ and a schema $H$, define the probability that a solution matching $H$ is selected by the proportional selection operator as
\begin{align*}
\alpha_{sel}(H,t)=\sum_{k:u^{k}\in{\uparrow}(H)}w_{k}.
\end{align*}
For the GA with uniform mutation, we have
\begin{align}
\alpha(H,t)&=\alpha_{sel}^{2}(H,t)(1-\pi_{m})^{\ord(H)}+\alpha_{sel}(H,t)\sum_{l:u^{l}\not\in{\uparrow}(H)}w_{l}\frac{(1-\pi_{m})^{\ord(H)}}{\big[2(1-\pi_{m})\big]^{\delta(u^{l},H)}} \notag\\
&\qquad+\sum_{k,l:u^{k},u^{l}\not\in{\uparrow}(H)}w_{k}w_{l}\frac{(2\pi_{m})^{h_{kl}}(1-\pi_{m})^{\ord(H)}}{\big[2(1-\pi_{m})\big]^{\delta_{H}(u^{k},u^{l})}}, \label{eq:schema}
\end{align}
where $h_{kl}=\big|\{j:H_{j}\not=\ast,u^{k}_{j}=u^{l}_{j}\not=H_{j}\}\big|$.
\end{theorem}

The general result in Theorem~\ref{thm:schema} provides an exact form of $\alpha(H,t)$, which is quite difficult to interpret and analyze. Accordingly, we derive a simple-to-analyze lower bound for $\alpha(H,t)$.

\begin{corollary}\label{cor:schema:lower}
Suppose conditions in Theorem~\ref{thm:schema} hold.
For $\pi_{m}\leq 0.5$, we have
\begin{align}
\alpha(H,t)\geq(1-\pi_{m})^{\ord(H)}\alpha_{sel}(H,t)^{2}+2^{-\ord(H)}\alpha_{sel}(H,t)\big[1-\alpha_{sel}(H,t)\big]+\big[1-\alpha_{sel}(H,t)\big]^{2}\pi_{m}^{\ord(H)}. \label{eq:schema:lower}
\end{align}
\end{corollary}

It can be seen from \eqref{eq:schema:lower} that the lower bound of $\alpha(H,t)$ gets larger when the schema selection probability $\alpha_{sel}(H,t)$ increases or the schema $H$ has lower order (i.e., $\ord(H)$ is small).
By definition, fitter schema $H$ leads to larger $\alpha_{sel}(H,t)$ and therefore higher $\alpha(H,t)$ and $\E\big[m(H,t+1)\big]$.
Since an expansion of the fitter schema $H$ is expected in a fixed-size population, fitter models matching $H$ are more likely to be generated in place of weaker models; see Section~\ref{sec:sim:schema} for a numerical verification.
Accordingly, the subsequent generation is anticipated to have higher average fitness.
This entitles the ``survival of the fittest'' phenomenon of the natural selection and acknowledges the evolvability of the GA.

\section{GA-assisted Multi-model Inference}
\label{sec:appl}

In this section, we describe how the GA helps multi-model inferences. Note that existing information-criteria based variable selection \citep[e.g.,][]{eBIC,KKC12,HBIC} and MCS procedures \citep[e.g.,][]{HLN11,MSCS:FTest,MSCS:LRT} typically concern the true model rather than the globally best model, which is the target of the GA. To bridge this gap, we first present a lemma suggesting that the true model indeed possess the lowest GIC value and therefore become the globally best model in large samples.

The following regularity condition is needed. 
\begin{assumption}
\begin{enumerate}
\item[(A1)] There exists a positive constant $C_{1}$ such that $\lambda_{\min}\big(\bX_{u^{0}}^{\top}\bX_{u^{0}}/n\big)>C_{1}$ for all $n$, where $u^{0}$ denotes the true model;
\item[(A2)] There is a positive constant $C_{2}$ such that
\begin{align*}
\inf_{u\not=u^{0},|u|<\tilde{s}}\bmu^{\top}(\bI-\bH_{u})\bmu\geq C_{2}n,
\end{align*}
where $\bmu=\bX\bbeta^{0}$ and $\bH_{u}=\bX_{u}(\bX_{u}^{\top}\bX_{u})^{-1}\bX_{u}^{\top}$ denotes the hat matrix of the model $u$, for some positive integer $\tilde{s}$ with $s\leq\tilde{s}<n$.
\end{enumerate}
\end{assumption}
Condition (A1) ensures the design matrix of the true model is well-posed and Condition (A2) is the asymptotic identifiability condition used in \citet{eBIC}, indicating that the model is identifiable if no model with comparable size can predict as well as the true model.

Recall that $\kappa_{n}$ is defined in the GIC formulation \eqref{eq:GIC}.
\begin{lemma}\label{lemma:GIC}
Suppose Assumption 4.1 holds, $\beps\sim\Norm(\bm{0},\sigma^{2}\bI)$, $\log d=O(n^{\tau})$, $\kappa_{n}=O(n^{\tau})$ and $\kappa_{n}\to\infty$ for some positive constant $\tau<1$.
Then for any positive integer $\tilde{s}$ satisfying $\tilde{s}\geq s$ and $\tilde{s}\log d=o(n)$, we have, as $n\to\infty$,
\begin{align}
\min_{u\in\cM_{\tilde{s}}-\{u^{0}\}}\GIC(u)-\GIC(u^{0})>0, \label{eq:selconsistent}
\end{align}\label{eq:lemma:GIC}
where $\cM_{\tilde{s}}=\big\{u\in\cM:|u|\leq\tilde{s}\big\}$.
\end{lemma}

\subsection{Variable Selection}
\label{sec:appl:vs}

The GA offers a practical way to perform variable selection by only searching the models generated by the GA instead of the whole model space. The existing information-criterion based selection methods search a  constrained model space $\cM_{\tilde{s}}$ for some $s\leq\tilde{s}\ll n$. However, by using the GA, we only need to evaluate at most $K\times T$ models (recall that $K$ and $T$ are the population size and the number of generations to convergence, respectively).
For example, under the simulation Case 1 with $(n,d,s,\rho)=(200,400,6,0.5)$ (see Section~\ref{sec:sim:settings}), it is nearly impossible to go through $\binom{400}{6}\approx 5.5\times 10^{12}$ models with the true size $6$, not to mention to compare all the models with sizes at most $\tilde{s}$ for some $6\leq\tilde{s}\ll n$.
On the other hand, the GA searches for the true model in all $500$ simulation runs, each with less than $1{,}750$ models evaluated ($K=92$ and $T\leq 19$ generations to convergence).

By combining Theorem~\ref{thm:mc} (b) and Lemma~\ref{lemma:GIC}, Proposition~\ref{prop:vs} shows that the true model becomes the best model in large samples and is eventually captured by the GA. Let $\Psi_{\tilde{s}}(t)$ denote the $t$-th generation of a GA population on the constrained model space $\cM_{\tilde{s}}$. The fitness function \eqref{eq:fitness} makes models of sizes at least $n$ nearly impossible to be generated. Accordingly, it is equivalent to setting $\tilde{s}=n-1$.

\begin{proposition}\label{prop:vs}
Suppose conditions in Lemma~\ref{lemma:GIC} hold and $\Psi_{\tilde{s}}(t)$ satisfies the conditions in Theorem~\ref{thm:mc}. Define
\begin{align*}
\what{u}(t)=\argmin_{u\in\Psi_{\tilde{s}}(t)}\GIC(u).
\end{align*}
Then we have
\begin{align}
\lim_{t\to\infty}\lim_{n\to\infty}P\left(\what{u}(t)=u^{0}\right)=1. \label{eq:GA:vs}
\end{align}
\end{proposition}

\subsection{Model Confidence Set}
\label{sec:appl:sms}

In this section, we construct practically feasible model confidence sets with the aid of the GA, in comparison with the one based on Theorem~\ref{thm:mcs}. The main idea is to employ the two-step testing procedure of \citet{V89}, given that the candidate models are produced by GA. 

Given a candidate model set $\Psi=\big\{u^{1},\ldots,u^{K}\big\}$, let $u^{\#}=\argmin_{u\in\Psi}\GIC(u)$ denote the best candidate model in $\Psi$.
Collect
\begin{align}
\cA_{\alpha}=\big\{u\in\Psi:H_{0,u}~\text{is not rejected at a significance level}~\alpha\big\} \label{eq:sms}
\end{align}
by performing the hypothesis testing
\begin{align}
H_{0,u}:~\text{Model}~u~\text{is not worse than}~u^{\#}\quad\text{vs.}\quad H_{1,u}:~\text{Model}~u~\text{is worse than}~u^{\#}. \label{eq:sms:test}
\end{align}
for every $u\in\Psi-\{u^{\#}\}$ at significance level $\alpha$.
We name the model confidence set $\cA_{\alpha}$ as survival model set (SMS) since the models therein survive the elimination testing \eqref{eq:sms:test}. Recall from Section~\ref{sec:schema} that the GA models, even after the globally best model is found, keep being improved until convergence. Accordingly, a manageable number of good (and sparse) models are included in the SMS when the GA is used to provide candidate models. Later, we use the relative size $|\cA_{\alpha}|/|\Psi|$ to measure the quality of the candidate model set in Section~\ref{sec:sim:quality}.

To perform the hypothesis testing \eqref{eq:sms:test} where $u$ and $u^{\#}$ may not be nested, we employ the two-step procedure of \citet{V89} by decomposing (\ref{eq:sms:test}) as first model distinguishability test
\begin{align}
H^{dis}_{0,u}:~u~\text{and}~u^{\#}~\text{are indistinguishable}\quad\text{vs.}\quad H^{dis}_{1,u}:~u~\text{and}~u^{\#}~\text{are distinguishable} \label{eq:sms:test:dis}
\end{align}
and if $H^{dis}_{0,u}$ is rejected, then a superiority test
\begin{align}
H^{sup}_{0,u}:~\E\big[\GIC(u)\big]\leq\E\big[\GIC(u^{\#})\big]\quad\text{vs.}\quad H^{sup}_{1,u}:~\E\big[\GIC(u)\big]>\E\big[\GIC(u^{\#})\big]. \label{eq:sms:test:sup}
\end{align}
The rejection of $H_{0,u}$ at significance level $\alpha$ is equivalent to that $H^{dis}_{0,u}$ and $H^{sup}_{0,u}$ are both rejected at significance level $\alpha$.
We note that the original superiority test in \citet{V89} is based on likelihood ratio, and therefore certain adjustment is needed for our case; see Section~\ref{supp:V89:sup} for detailed description. The {\sf R} package {\bf nonnest2} \citep{nonnest2} is used to test \eqref{eq:sms:test:dis} and extract necessary quantities for the GIC-based superiority test \eqref{eq:sms:test:sup}.

The following proposition justifies the asymptotic validity of the constructed SMS.
\begin{proposition}\label{prop:sms}
Suppose conditions in Proposition~\ref{prop:vs} hold and $\Psi_{\tilde{s}}(t)$ satisfies the conditions in Theorem~\ref{thm:mc}. 
Let $\cA_{\alpha}(t)$ denote a $100(1-\alpha){\%}$ SMS with $\Psi_{\tilde{s}}(t)$ serving as the candidate model set.
Then we have
\begin{align*}
\lim_{t\to\infty}\lim_{n\to\infty}P\big(u\in\cA_{\alpha}(t)\big)\geq 1-\alpha
\end{align*}
for all $u\in\Psi_{\tilde{s}}(t)-\{u^{0}\}$ such that $H_{0,u}$ is not rejected at significance level $\alpha$.
\end{proposition}

\section{Simulation Studies}
\label{sec:sim}

In this section, we conduct extensive simulation studies to provide numerical support for the new schema theory supplied in Section~\ref{sec:schema} and show that the GA outperforms the RP method and the SA algorithm of \citet{NR17}.
The simulated data were generated based on the linear model \eqref{eq:LM} with $\beps\sim\Norm_{n}(\bm{0},\bI)$.
Each row of the design matrix $\bX$ was generated independently from $\Norm_{d}(\bm{0},\bSig)$, where $\bSig$ is a Toeplitz matrix with the $(k,l)$-th entry $\Sigma_{kl}=\rho^{|k-l|}$ for $\rho=0,0.5$ and $0.9$. Results were obtained based on $500$ simulation replicates.

\subsection{Simulation Settings}
\label{sec:sim:settings}

We consider six simulation cases below with both high-dimensional (i.e., $d\geq n$; Cases 1--4) and low-dimensional (Cases 5 and 6) settings. Cases 1 and 2 with $\rho=0$ refers to the first two simulation cases used in \citet{SOIL}. Case 3 is inspired from the simulations settings used in \citet{NR17}, but our settings ensure $X_{s+1}$ and $X_{s+2}$ are always marginally distributed as $\Norm(0,1)$ for any $\rho\in[0,1)$.
Case 4 is similar to Case 3 but with weak signals.
Cases 5 and 6 refers to the simulation example 2 of \citet{mBIC}, with weak signals in Case 6.

Let $\bm{1}_{p}$ and $\bm{0}_{p}$ denote the $p$-dimensional vectors of $1$'s and $0$'s, respectively.
\begin{description}
\item[Case 1:] $\bbeta^{0}=(4\bm{1}_{s-2}^{\top},-6\sqrt{2},\frac{4}{3},\bm{0}_{d-s}^{\top})^{\top}$.
\item[Case 2:] $\bbeta^{0}$ as in Case 1.
Re-define $X_{s+1}=0.5X_{1}+2X_{s-2}+\eta_{1}$, where $\eta_{1}\sim\Norm(0,0.01)$.
\item[Case 3:] $\bbeta^{0}=(3\bm{1}_{s}^{\top},\bm{0}_{d-s}^{\top})^{\top}$.
Re-define $X_{s+1}=\frac{2}{3\sqrt{(1+\rho)}}(X_{1}+X_{2})+\eta_{2}$ and $X_{s+2}=\frac{2}{3\sqrt{(1+\rho)}}(X_{3}+X_{4})+\eta_{3}$, where $\eta_{1},\eta_{2}\overset{iid}{\sim}\Norm(0,1/9)$.
\item[Case 4:] $\bbeta^{0}=(3\log(n)/\sqrt{n}\bm{1}_{s}^{\top},\bm{0}_{d-s}^{\top})^{\top}$.
Re-define $X_{s+1}$ and $X_{s+2}$ as in Case 3.
\end{description}
Two cases are set up for moderate dimensional (i.e., $d<n$) scenarios:
\begin{description}
\item[Case 5:] $\beta^{0}_{1}\geq\cdots\geq\beta^{0}_{s}$ are iid $\Unif(0.5,1.5)$, sorted decreasingly, and $\beta_{j}=0$ for $j>s$.
\item[Case 6:] $\bbeta^{0}=(3\log(n)/\sqrt{n}\bm{1}_{s}^{\top},\bm{0}_{d-s}^{\top})^{\top}$.
\end{description}

For the GA implementation, we use
\begin{align}
\GIC(u)=n\log\what{\sigma}_{u}^{2}+3.5|u|\log d, \label{eq:PLIC}
\end{align}
to evaluate models.
This choice of $\kappa_{n}=3.5\log d$ makes the GIC coincide with the pseudo-likelihood information criterion \citep{PLIC} and the high-dimensional BIC \citep{HBIC}. The penalization constant $3.5$ is specifically used due to the superior performance shown in the simulation studies in \citet{HBIC}.
It should be mentioned that \eqref{eq:PLIC} works well regardless the relationship between $n$ and $d$ \citep[e.g.,][]{PLIC,HBIC}.
Our {\sf Python} implementation of the GA, the RP and the SA is publicly available in the Github repository  \href{https://github.com/aks43725/cand}{\url{https://github.com/aks43725/cand}}.

\subsection{Schema Evolution}
\label{sec:sim:schema}

The discussion followed by Corollary~\ref{cor:schema:lower} concludes that fitter schema $H$ leads to larger $\alpha_{sel}(H,t)$ (the probability of selecting a model matching $H$ in the $t$-th generation) and hence larger $\E\big[m(H,t+1)\big]$ (the expected number of models matching $H$ in the $(t+1)$-th generation).
In the following we provide empirical evidence by observing that $m(H,t+1)$ aligns with $\alpha_{sel}(H,t)$ for three schemata:
\begin{align*}
H^{1}=(\bm{1}_{s},\bm{0}_{2s},\ast,\ldots,\ast),\quad H^{2}=(\bm{1}_{s+2},\ast,\ldots,\ast)\quad\text{and}\quad H^{3}=(\bm{1}_{s-1},0,\ast,\ldots,\ast),
\end{align*}
which represent good, fair and bad performing schemata, respectively.
In particular, $H^{1}$ is expected to perform the best by covering good models such as the true model.
The $2s$ $0$'s are placed to deteriorate its overall performance through ruling out some models that are too good to observe the evolution of $m(\cdot,t)$ and $\alpha_{sel}(\cdot,t)$.
$H^{2}$ is expected to be slightly worse than $H^{1}$ because models matching it are all overfitting by having at least two false discoveries.
We anticipate $H^{3}$ to have the worst performance due to missing one true signal.
Note that ${\uparrow}(H^{1})\cap{\uparrow}(H^{2})\cap{\uparrow}(H^{3})$ does not cover the whole model space $\cM$.
For implementation, we used uniform mutation as needed in the theoretical conditions.
Moreover, since the GA with initial population provided by the RP is too good to observe the evolution process, we used an approach proposed in Section~\ref{supp:initial} to randomly generate an initial population.

\begin{figure}[!tb]
\centering
\includegraphics[width=0.9\textwidth]{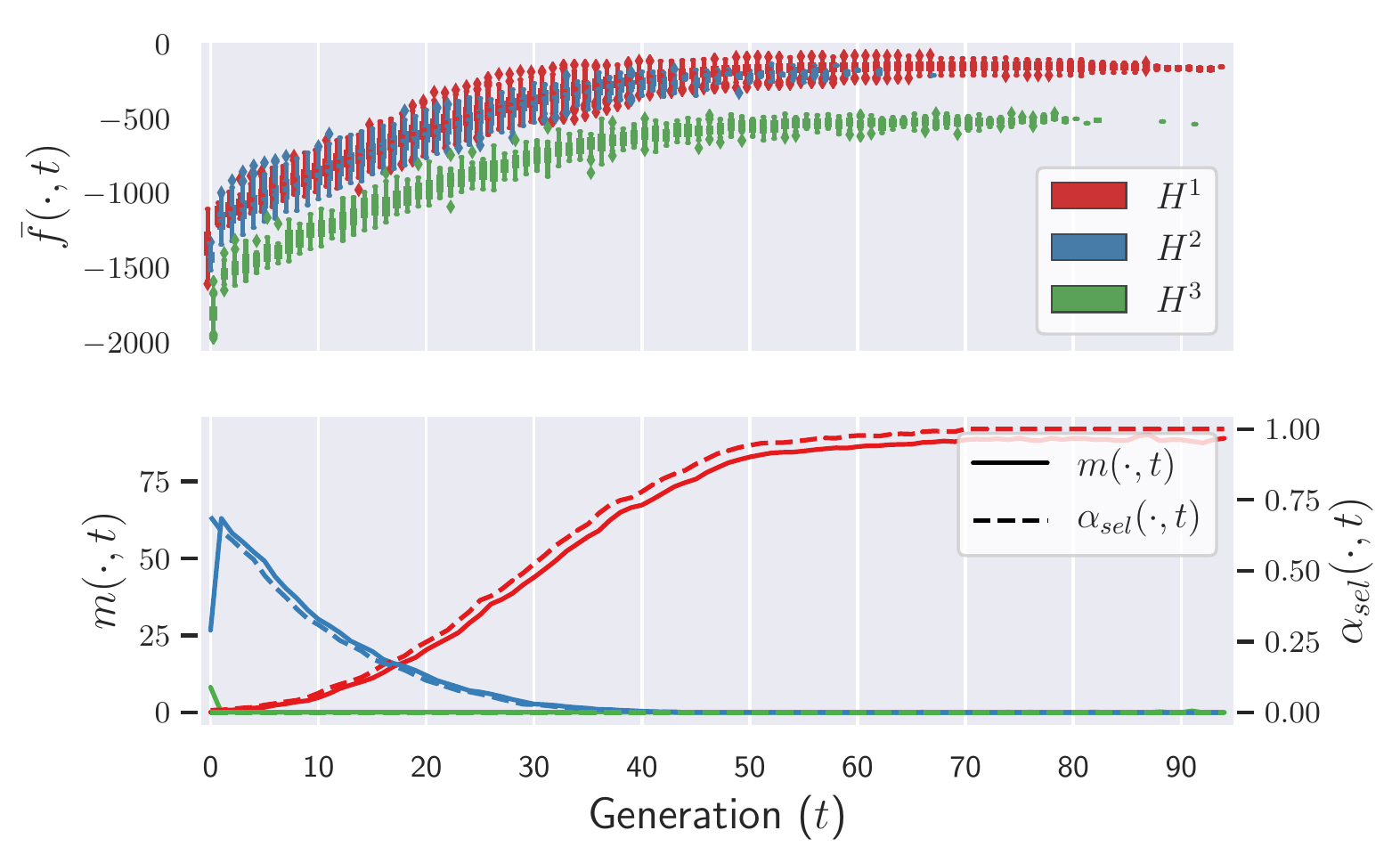}
\caption{\label{fig:schema}
Schema performance (upper panel) and evolution (lower panel) under Case 3 with $(n,d,s,\rho)=(200,400,6,0.5)$.}
\end{figure}

Figure~\ref{fig:schema} is obtained under Case 3 with $(n,d,s,\rho)=(200,400,6,0.5)$ and serves as a representative example since other cases (included in supplementary, Section~\ref{supp:sim:schema}) exhibits similar patterns.
The upper panel confirms our performance assertion on the overall schema performance, i.e., $H^{1}$ is slightly better than $H^{2}$ and $H^{3}$ is the worst.
From the lower panel, it is evident that the pattern of $m(H^{1},t+1)$ aligns with that of $\alpha_{sel}(H^{1},t)$ in all cases.
In addition, the strong schema $H^{1}$ evolves to take over the whole population eventually even it is a minority at the beginning, and vice versa for the weaker schema $H^{2}$.
On the other hand, the evolution process of $H^{3}$ illustrates a typical example that a particularly weak schema extincts soon and never rises again.
In summary, a good schema expands and a weak one diminishes over generations, resulting in an improved average fitness until convergence.

\subsection{Comparison with Existing Methods}

In this section, we compare the GA with the RP and the SA in terms of computation time, quality of candidate model sets, and performance of multi-model inference applications such as variable importance, model confidence set and model averaging.
For the RP, we collect the unique models on the regularization paths of Lasso, SCAD and MCP using the \textsf{Python} package \textbf{pycasso}.
Recall that the GA takes the RP for the initial population.
The SA is implemented to search for models of sizes appeared in the last GA generation, and the best $K$ models are kept as the final candidate model set.
Other tuning parameters are settled according to the simulation settings in \citet{NR17}.

In the following, we show that the GA evidently improve the models generated by the RP in reasonable computation time, and that the SA takes a long time to implement but produces at most comparable results to those of the GA.
In particular, the GA exhibits the best performance in all cases in terms of variable selection and quality of candidate model set.
In terms of model averaging and variable importance, the GA performs at least comparably to the RP and the SA in high-dimensional cases, while just comparably under low-dimensional settings.

\subsubsection{Computation Time}
\label{sec:sim:time}

The averaged computation time for the three methods are depicted in Figure~\ref{fig:time}.
It is obvious that the GA is a bit slower than the RP but way much (like more than $10$ times) faster than the SA.

\subsubsection{Variable Selection}
\label{sec:sim:vs}

\begin{figure}[!tb]
\centering
\includegraphics[width=\textwidth]{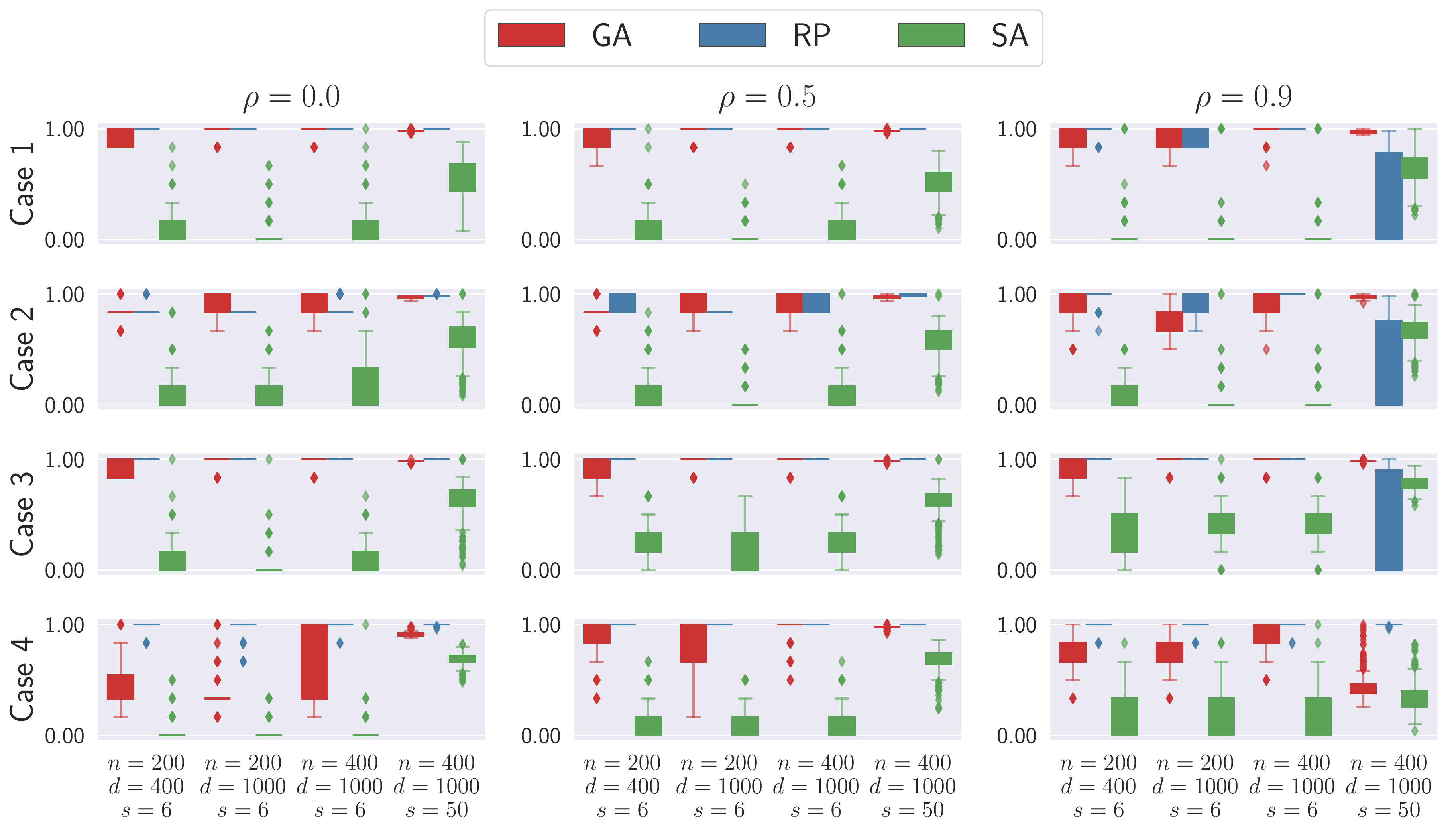}
\includegraphics[width=\textwidth]{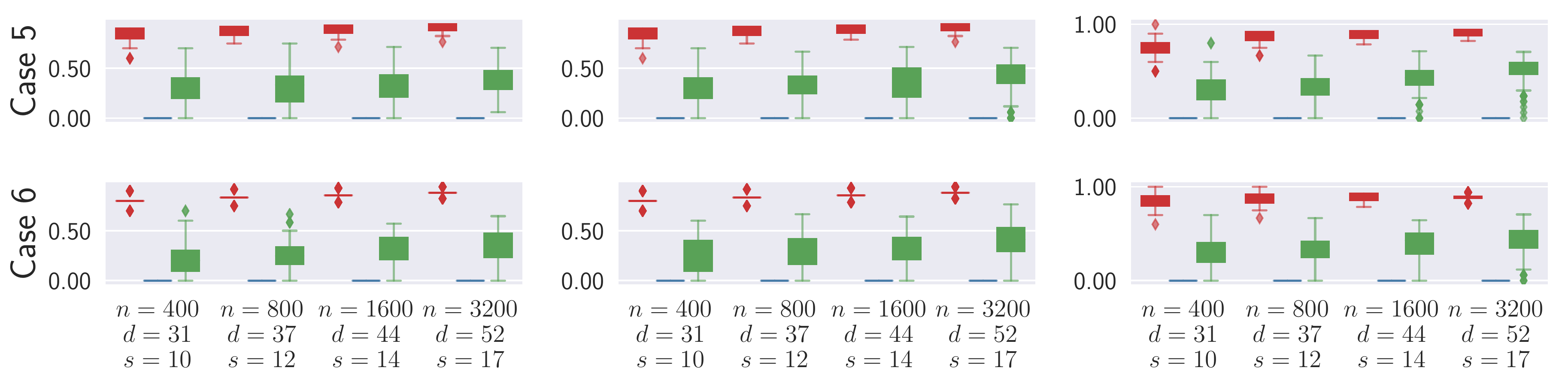}
\caption{\label{fig:mselect:PSR}
Positive selection rate (PSR) of the best model.}
\end{figure}

\begin{figure}[!tb]
\centering
\includegraphics[width=\textwidth]{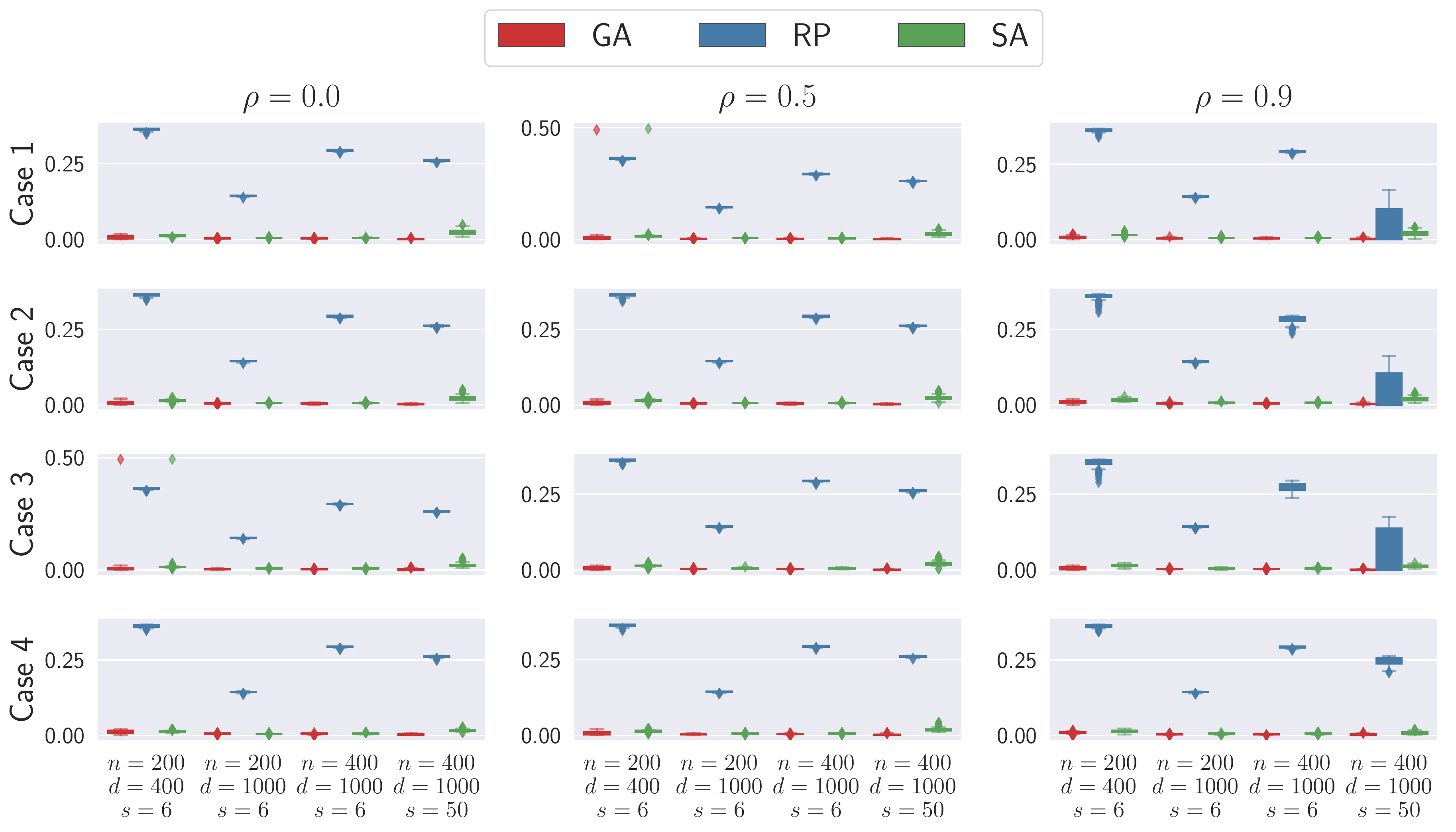}
\includegraphics[width=\textwidth]{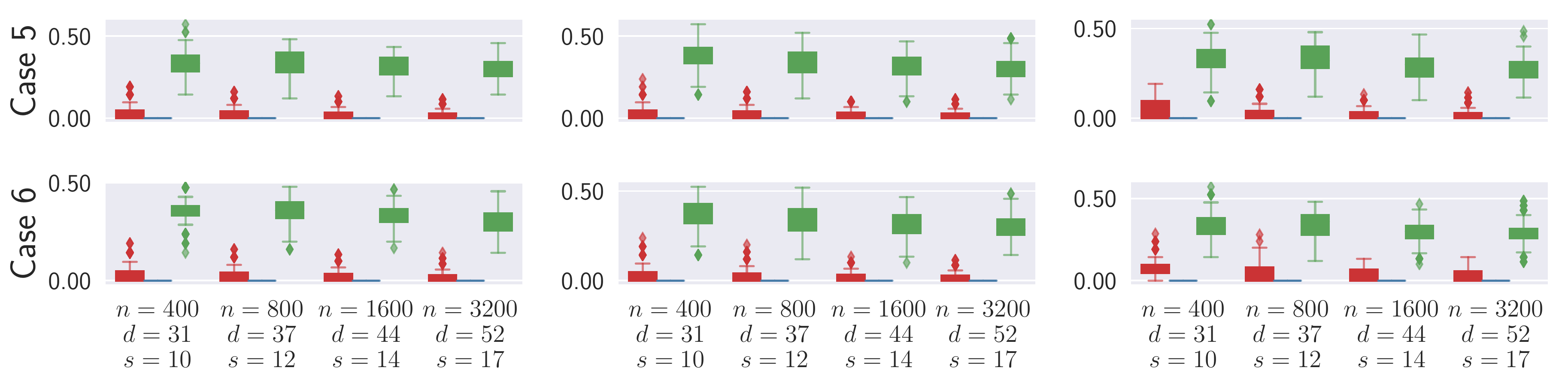}
\caption{\label{fig:mselect:FDR}
False discovery rate (FDR) of the best model.}
\end{figure}

To evaluate the performance of variable selection, the boxplots of the positive selection rate (PSR, the proportion of true signals that are active in the best model) and the false discovery rate (FDR, the proportion of false signals that are active in the best model) are drawn in Figure~\ref{fig:mselect:PSR} and Figure~\ref{fig:mselect:FDR}, respectively.
We see that the GA-best model gives fairly high PSR and low FDR in all cases, demonstrating excellent variable selection performance.
Under high-dimensional settings (Cases 1--4), the RP produces high PSR but also high FDR, while the SA results in the opposite (PSR and FDR are both low).
For moderate dimensional cases (Cases 5 and 6), both of the RP and the SA give low PSR and FDR.
In summary, the GA-best model possesses much better variable selection performance than those from the RP and the SA.

\subsubsection{Quality of Candidate Models}
\label{sec:sim:quality}

\begin{figure}[!tb]
\centering
\includegraphics[width=\textwidth]{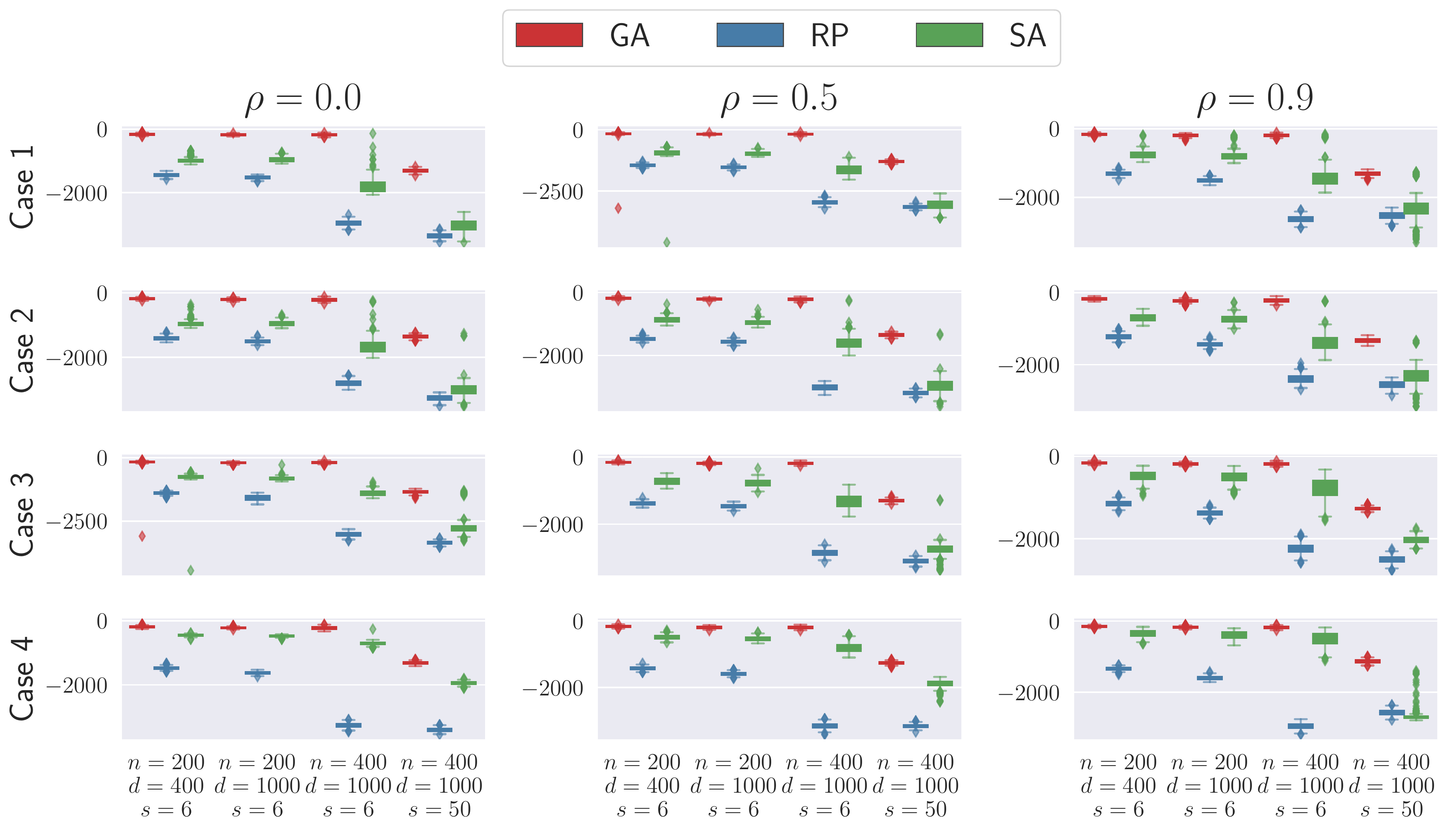}
\includegraphics[width=\textwidth]{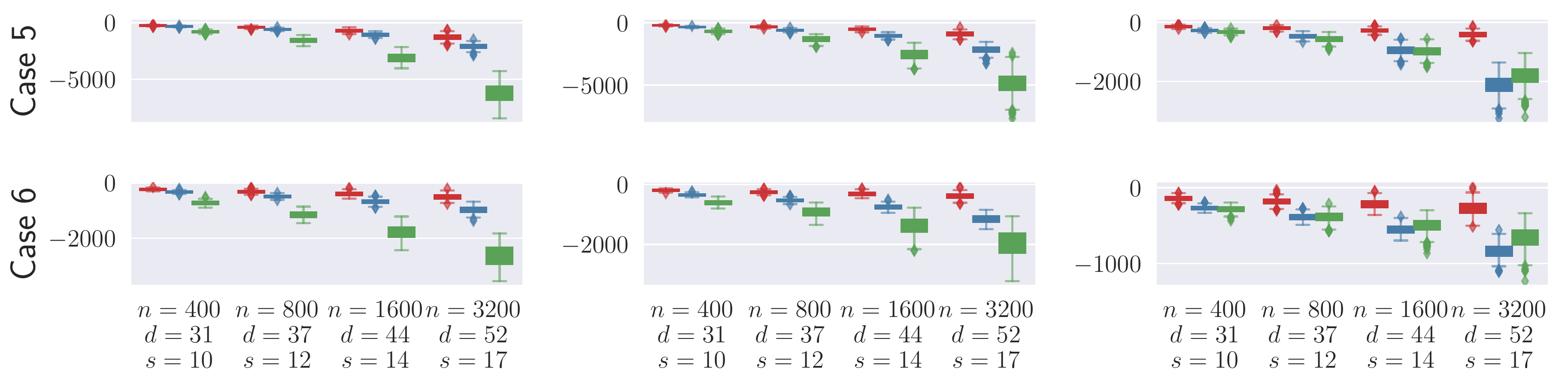}
\caption{\label{fig:avgfit}
Boxplots of the average fitness of the candidate model sets.}
\end{figure}

We evaluate the quality of candidate model sets using two criteria: (i) the average fitness and (ii) the relative size of $95{\%}$ SMSs (see Section~\ref{sec:appl:sms} for the SMS construction) to the original candidate model set.
Figure~\ref{fig:avgfit} exhibits the boxplots of average fitness and suggests that the GAs produce the fittest models in all cases.
The SA takes the second place in high-dimensional cases (Cases 1--4), yet is outperformed by the RP in moderate dimensional cases (Cases 5 and 6) with $\rho=0$ and $0.5$, where the covariates are not strongly correlated.
To conclude, the candidate model set produced by the GA possesses the best quality among the three approaches.

\begin{figure}[!tb]
\centering
\includegraphics[width=\textwidth]{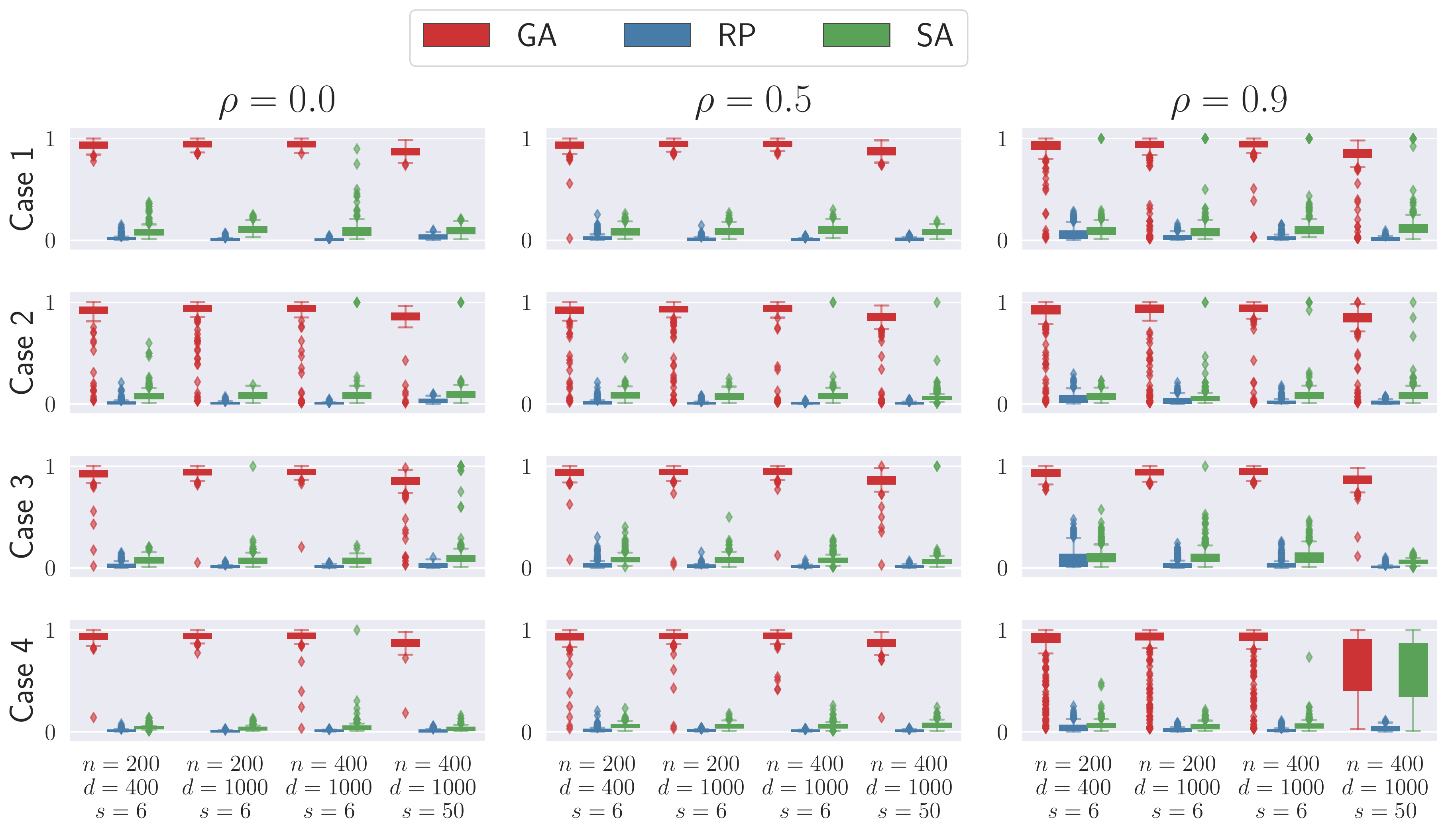}
\includegraphics[width=\textwidth]{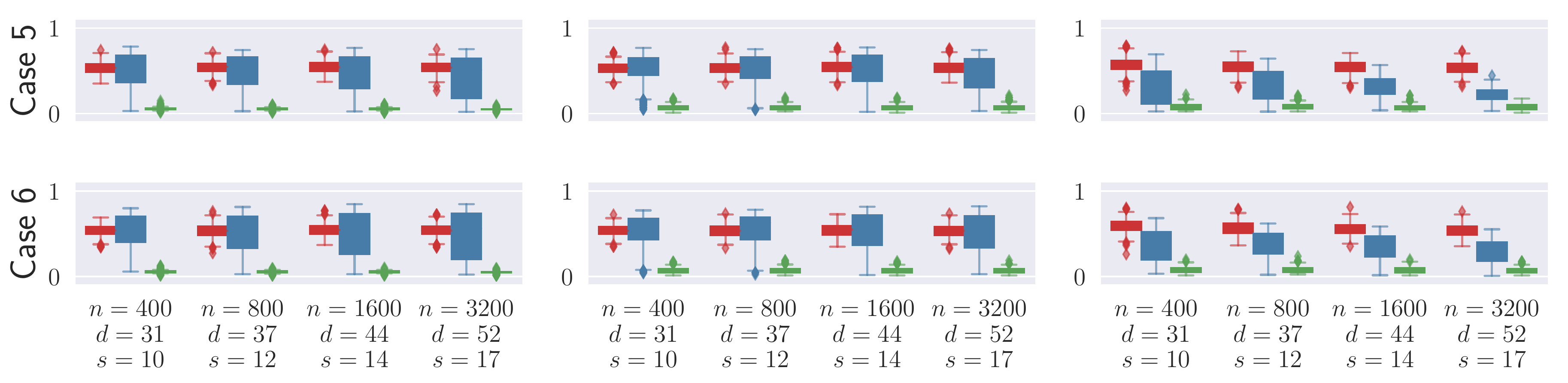}
\caption{\label{fig:sms}
Relative size of $95{\%}$ SMS over the original candidate model set.}
\end{figure}

Figure~\ref{fig:sms} displays the boxplots of the relative size of $95{\%}$ SMSs $\cA_{0.05}$ against the original candidate model set $\Psi$, i.e., $|\cA_{0.05}|/|\Psi|$, where larger values indicate better quality of $\Psi$.
We see that the relative sizes for the GA are typically higher than those from the RP and SA in all cases, and are close to $1$ in high-dimensional settings (e.g., Cases 1--4).
This supports the conclusion we made about the quality of candidate models in the previous paragraph.

\subsubsection{Model Averaging}
\label{sec:sim:ma}

Model averaging, especially in high-dimensional predictive analysis, is a prominent application of multi-model inference.
The GA does not perform significantly better than the RP and the SA in model averaging, but exhibits better applicability than the RP, and greater robustness than the SA.

Given a candidate model set $\Psi=\big\{u^{1},\ldots,u^{K}\big\}$, the model averaging predictor is defined by
\begin{align}
\what{\bY}=\sum_{k=1}^{K}w_{k}\what{\bY}_{u^{k}}, \label{eq:ma}
\end{align}
where $\what{\bY}_{u^{k}}=\bX_{u}(\bX_{u}^{\top}\bX_{u})^{-1}\bX_{u}^{\top}\bY$ are the least-squares predictors and $w_{k}$ with $0\leq w_{k}\leq 1$ denote the model weights of $u^{k}$ for $k=1,\ldots,K$.
We use the root mean squared error (RMSE) defined by
\begin{align*}
\sqrt{n^{-1}(\bY-\what{\bY})^{\top}(\bY-\what{\bY})}
\end{align*}
to assess the performance of model averaging.

Two model weighting schemes are considered to obtain the model weights $w_{k}$: (i) GIC-based weights as in \eqref{eq:modelweight} with $f_{k}$ replaced by $-\GIC(u^{k})$, and (ii) the weighting approach proposed by \citet{AL14}, which we called it the ``AL weighting'' hereafter (see Section~\ref{supp:AL14} for detailed construction).
We note that (i) is the the most commonly used model weighting scheme in multi-model inference (e.g., Akaike weights \citep{A78,B87,BA04:book,AICw} and Bayesian model averaging \citep{BMA}), and (ii) is developed for optimal predictive performance in high-dimensional model averaging.

\begin{figure}[!tb]
\centering
\includegraphics[width=\textwidth]{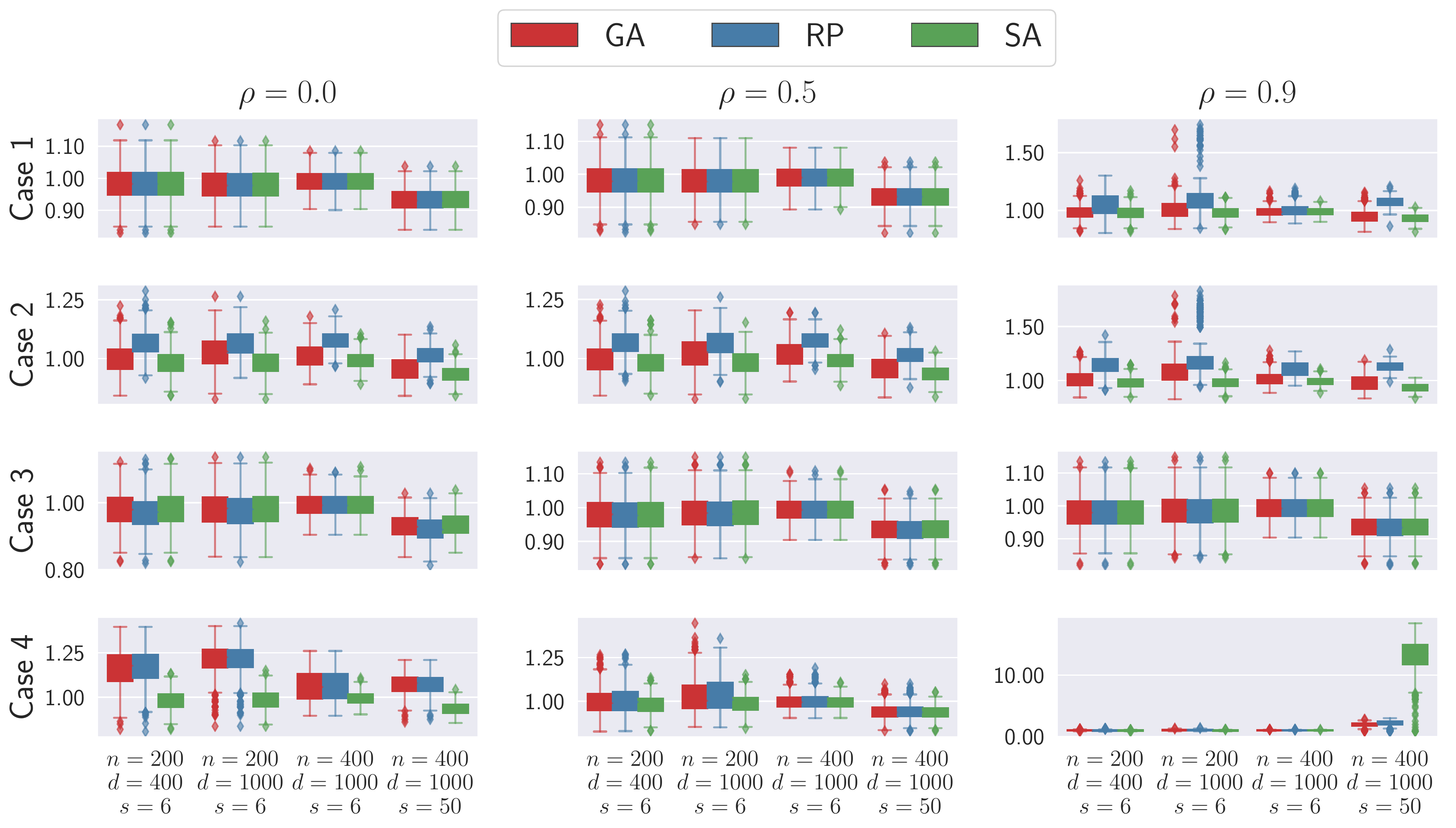}
\includegraphics[width=\textwidth]{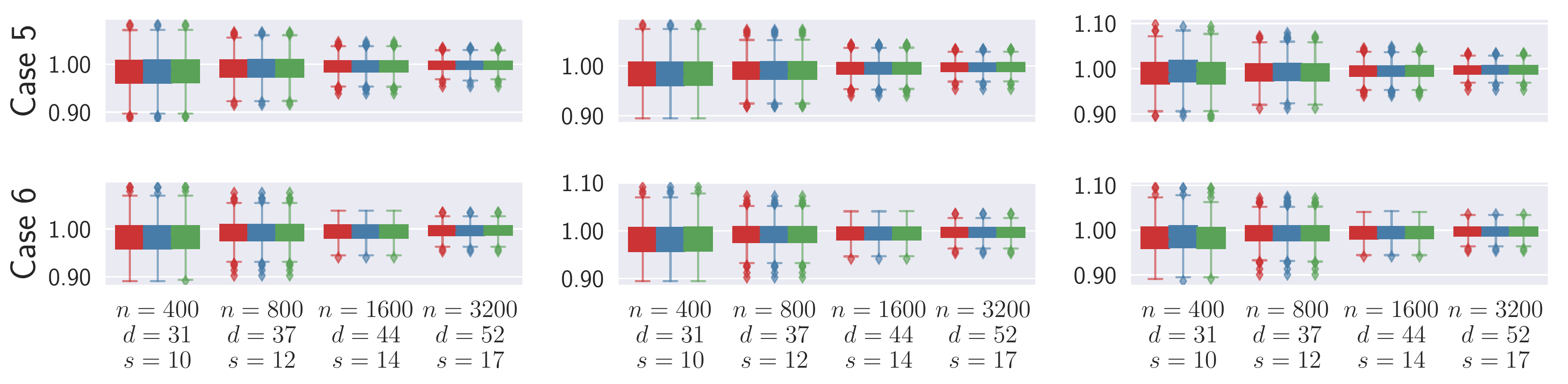}
\caption{\label{fig:ma:PLIC}
Boxplots of the RMSE obtained by model averaging using the GIC-based weighting.}
\end{figure}

Figure~\ref{fig:ma:PLIC} displays the boxplots of the RMSE using the GIC-based model weighting, showing that the GA exhibits good and robust (in contrast to the wildly high RMSE by SA in Case 4 with $(n,d,s,\rho)=(400,1000,50,0.9)$; see Remark~\ref{rmk:SA} for more details) results over all cases.
On the other hand, the RP is obviously worse than the GA in Case 2, and the SA's performance is just comparable to that of the GA.
The three methods perform similarly in the rest cases (i.e., Cases 1, 3, 5 and 6).

\begin{figure}[!tb]
\centering
\includegraphics[width=\textwidth]{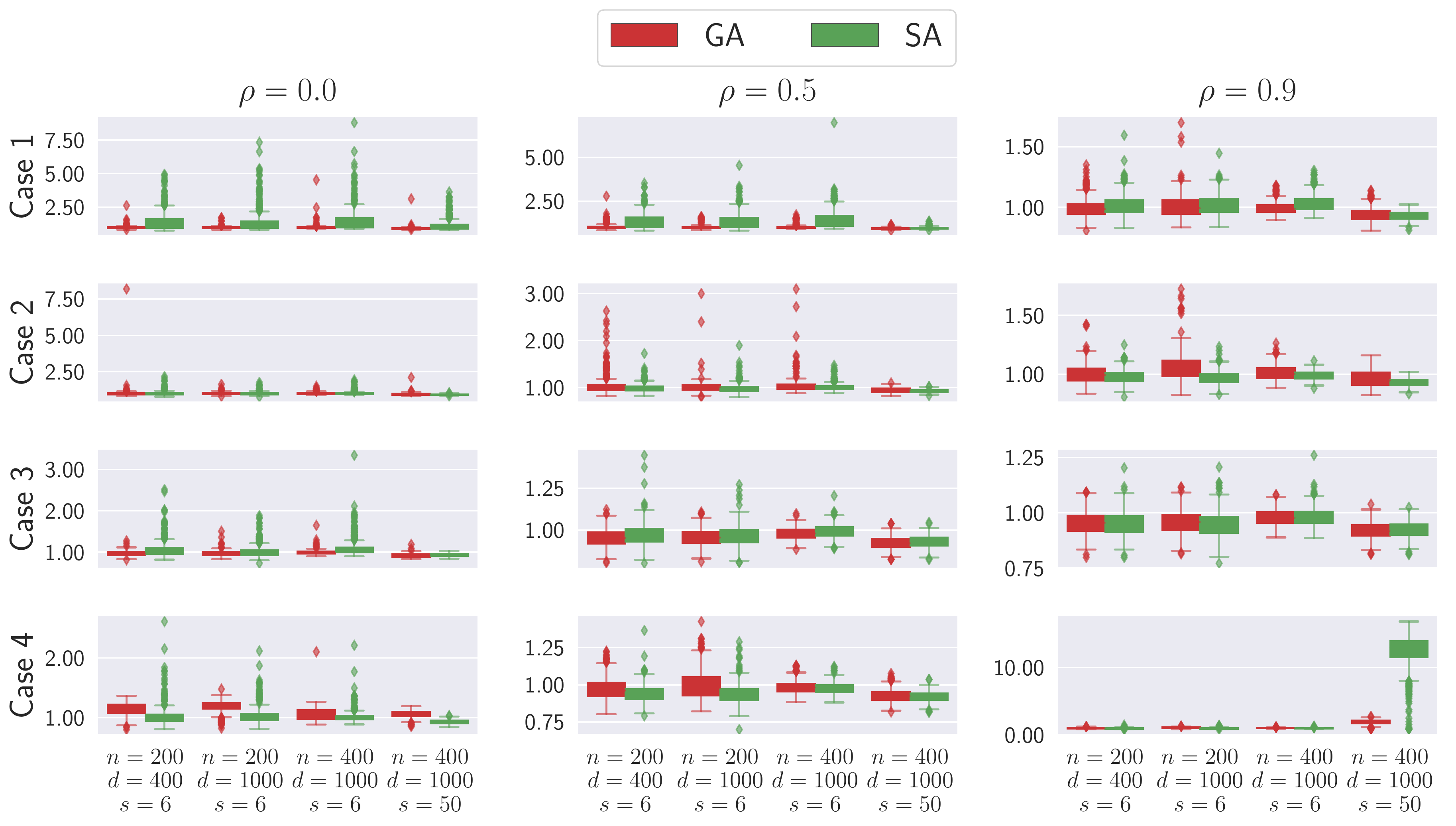}
\includegraphics[width=\textwidth]{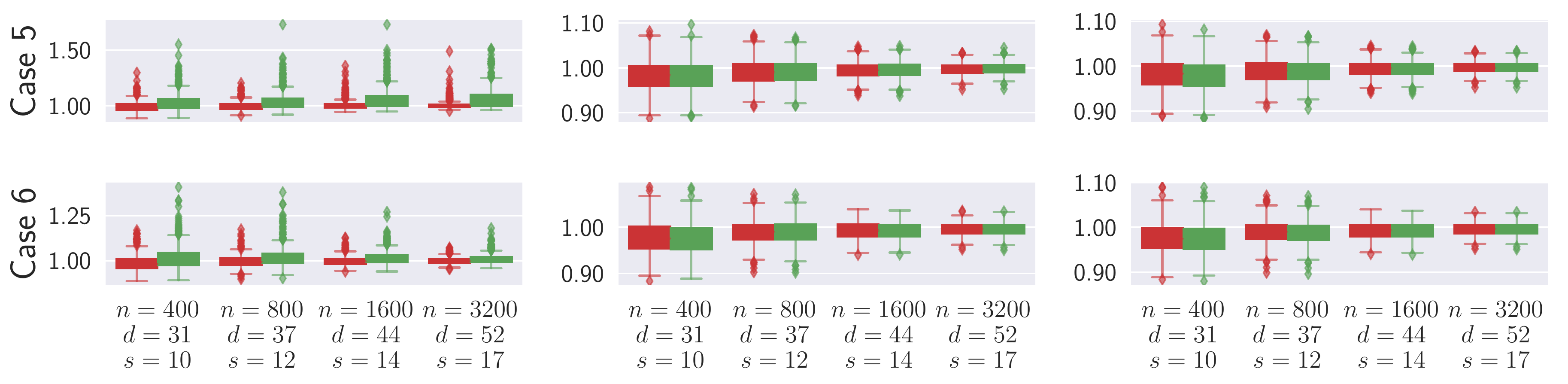}
\caption{\label{fig:ma:AL14}
Boxplots of the root mean squared error obtained by high-dimensional model averaging approach of \citet{AL14}.
The RP method fails to perform in all cases and thus is not shown.}
\end{figure}

The RMSEs obtained by the AL weighting are shown in Figure~\ref{fig:ma:AL14}.
Different from the results using the GIC-based model weights, the GA behaves slightly better than SA in some cases (e.g., Case 1 with $\rho=0.0$ and $0.5$ and Case 3 with $\rho=0.0$) and comparably in the rest.
Yet similarly, the GA performs robustly and the SA has wildly high RMSE in Case 4 with $(n,d,s,\rho)=(400,1000,50,0.9)$.
On the other hand, the results for the RP are omitted due to the computational infeasibility (inverting a singular matrix) in generating the AL weights.
Accordingly, the GA is shown to possess better applicability in optimal high-dimensional model averaging.

\subsubsection{Variable Importance}
\label{sec:sim:varimp}

To evaluate the performance of high-dimensional variable importance, we employ the sparsity oriented importance learning \citep[SOIL;][]{SOIL} defined by
\begin{align*}
\SOIL_{j}\equiv\SOIL(j;\bw,\Psi)=\sum_{k=1}^{K}w_{k}\dsone(u^{k}_{j}=1)
\end{align*}
with the GIC-based model weights $w_{k}$ given in \eqref{eq:modelweight}.
It can well separate the variables in the true model from the rest in the sense that $\SOIL_{j}$ rarely gives $0$ ($1$) if the variable $j$ is (not) in the true model.
Moreover, it rarely gives variables not in the true model significantly higher values than those in the true model even if the signal is weak.
In the original work \citep{SOIL}, the candidate models were generated using the RP method.
Our results indicate that the GA performs at least comparably to the SA and the RP in separating the true signals from the rest.

\begin{figure}[!tb]
\centering
\includegraphics[width=\textwidth]{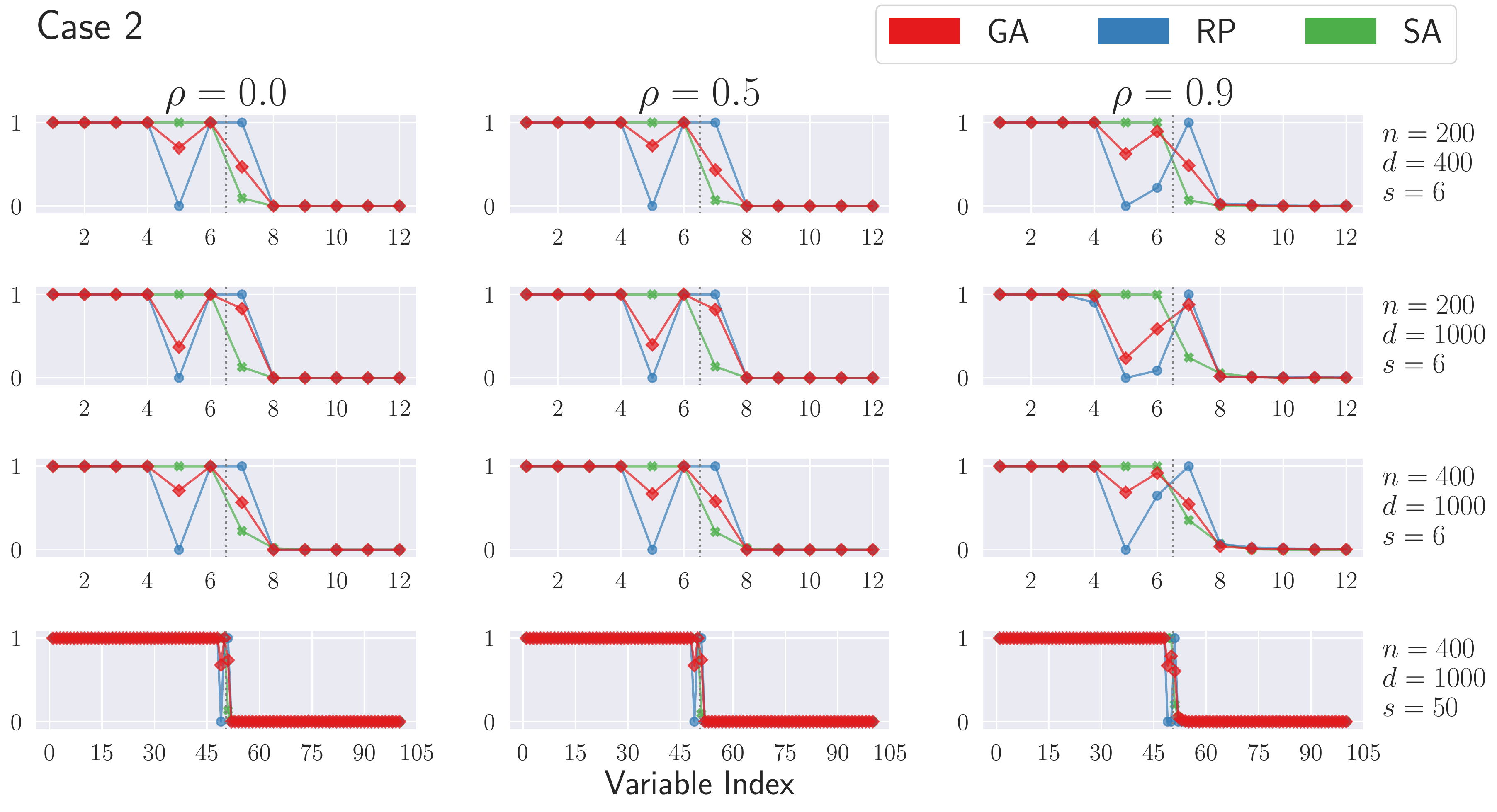}
\caption{\label{fig:soil:case2}
(Case 2)
Averaged SOIL measures.}
\end{figure}

\begin{figure}[!tb]
\centering
\includegraphics[width=\textwidth]{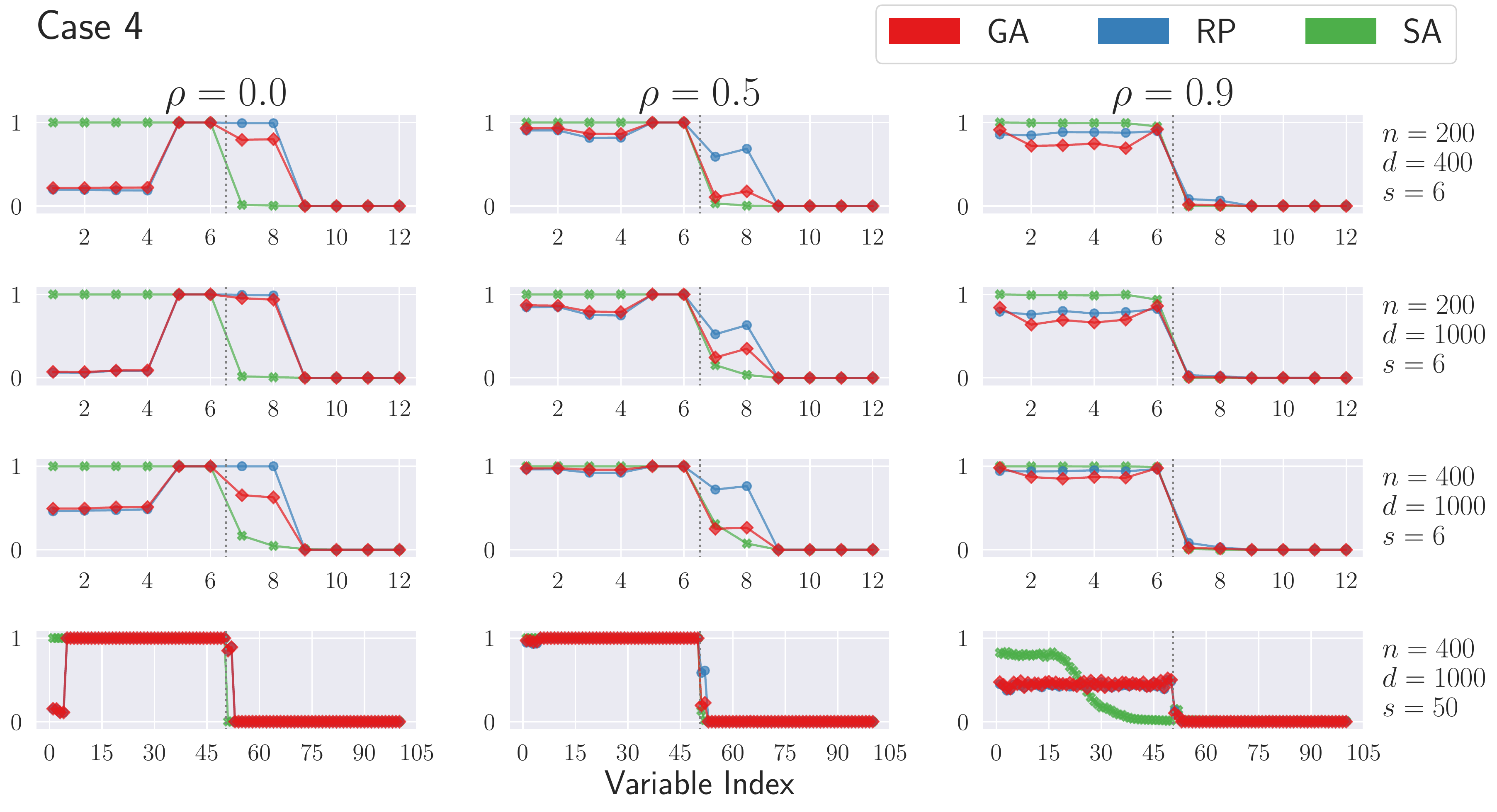}
\caption{\label{fig:soil:case4}
(Case 4)
Averaged SOIL measures.}
\end{figure}

Figure~\ref{fig:soil:case2} and Figure~\ref{fig:soil:case4} depict the averaged SOIL values for the first $2s$ variables for Cases 2 and 4, respectively, where the active ones are before the vertical gray line and the rest are not shown due to $\SOIL_{j}\approx 0$ for $j>2s$ no matter which method was used for candidate model preparation.
Results for Cases 1, 3, 5 and 6 are presented in supplementary (Section~\ref{supp:sim:varimp}) due to high similarity among the three methods. the GA exhibits the best performance that separate the true signals from the rest. Specifically, the resulting SOIL values are by no means close to $0$ and $1$ for truly active and inactive variables, respectively.
On the other hand, in Case 2 the RP results in $\SOIL_{s-2}\equiv 1$ and $\SOIL_{s+1}=0$, where $X_{s-2}$ is a true signal and $X_{s+1}$ is not.
Moreover, in Case 4 with $(n,d,s,\rho)=(200,1000,50,0.9)$, since the SA results in $\SOIL_{j}\leq 0.03$ for $j=38,\ldots,50$, these $13$ true signals may easily be regarded as not important.

\begin{remark}\label{rmk:SA}
From Case 4 with $(n,d,s,\rho)=(200,1000,50,0.9)$, we note that the SA's performance in model averaging and variable importance critically depends on the model size specification.
Recall that the SA only searches for the models with sizes resulting from the GA candidate models.
For this simulation case, the GA model sizes are around a half of the number of strong signals, i.e., $s/2$.
Such model size misspecification causes the SA to perform poorly in model averaging and variable importance. On the other hand, the GA still behaves well even when all of its resulting candidate models miss certain number of true signals.
\end{remark}

\section{Real Data Example}
\label{sec:realdata}

In this section, we present two real data examples to exhibit the usefulness of the proposed GA. Additionally, hypothesis testing \eqref{eq:sms:test} was conducted to compare models in terms of the GIC.

\subsection{The Riboflavin Dataset}
\label{sec:realdata:riboflavin}

We first introduce the riboflavin (vitamin B) production dataset that was widely studied in high-dimensional variable selection literature \citep[e.g.,][]{BKM14,JM14:dLasso,LM15,CLW16,HY17}.
The response variable is the logarithm of the riboflavin production rate in Bacillus subtilis for $n=71$ samples and the covariates are the logarithm of the expression level of $d=4{,}088$ genes.
Please see more details in Supplementary Section A.1 of \citet{BKM14}. 

The proposed GA delivers new insights by yielding better variable selection results than the existing works. From Table~\ref{tab:riboflavin:compare}, the GA-best model contains only one active gene {\sf XHLA-at} which was not identified by previous approaches. However, it turns out to be the fittest model (with all $p$-values $<0.0001$) among those listed.
Moreover, the importance of the gene {\sf XHLA-at} is confirmed by having $\SOIL_{\textsf{XHLA-at}}=1$ and all other SOIL values less than $0.01$.
Accordingly, we suggest a further investigation on the gene {\sf XHLA-at} is needed from scientists.

\begin{table}[!tb]
\centering\small
\caption{\label{tab:riboflavin:compare}
Variable selection results and GIC values of the selected models for the riboflavin dataset.}
\begin{tabular}{llr}
\hline\hline
Method & Active Covariates & \multicolumn{1}{r}{GIC} \\
\hline
Proposed GA                                   & {\sf XHLA-at}                & $-20.520$ \\
Multisplit procedure \citep{MMB09}$^\dagger$  & {\sf YXLD-at}                & $-14.357$ \\
Stability selection \citep{StabSel}$^\dagger$ & {\sf YXLD-at}, {\sf YOAB-at}, {\sf LYSC-at} & $-1.431$ \\
Debiased Lasso \citep{JM14:dLasso}            & {\sf YXLD-at}, {\sf YXLE-at} &  $15.643$ \\
B-TREX \citep{LM15}                           & {\sf YXLD-at}, {\sf YOAB-at}, {\sf YXLE-at} & $10.624$ \\
$\mathrm{AV}_{\infty}$ \citep{CLW16}          & {\sf YXLD-at}, {\sf YOAB-at}, {\sf YEBC-at}, & $-5.681$ \\
                                              & {\sf ARGF-at}, {\sf XHLB-at} & \\
RP, SA, and Ridge-type projection \citep{B13}$^\dagger$ & None               & $-11.775$ \\
\hline\hline
\noalign{\vskip 0.5ex}
\multicolumn{3}{p{0.9\textwidth}}{$^\dagger$Obtained by \citet{BKM14} using the {\sf R} package \textbf{hdi}.}
\end{tabular}
\end{table}

Table~\ref{tab:riboflavin:mumin} summarizes the results of $95{\%}$ SMSs (see Section~\ref{sec:appl:sms}), and shows the GA outperforms the RP and the SA in terms of the quality of candidate model set and model averaging. For the former, besides the much fittest (i.e., lowest GIC) model, the GA also gives the highest relative size of $95{\%}$ SMSs of $56/67=83.58{\%}$ (compared to $1/54=1.85{\%}$ for the RP and $11/16=68.75{\%}$ for the SA). For model averaging, the GA results in the smallest RMSE using the GIC-based weighting. Moreover, as the only method leading to successful AL weighting (see Section~\ref{sec:sim:ma}) computation, the GA is shown to possess better applicability in optimal high-dimensional model averaging.

\begin{table}[!tb]
\centering\small
\caption{\label{tab:riboflavin:mumin}
Results of the relative size of $95{\%}$ SMSs and model averaging for the riboflavin dataset.}
\begin{tabular}{lrrrr}
\hline\hline
&&& \multicolumn{2}{r}{RMSE of Model Averaging} \\
\cline{4-5}
Method & $\#(\text{Candidate Models})$ & $\#(\text{Models in}~95{\%}~\text{SMS})$ & GIC-based & AL \\
\hline
GA & $67$ & $56$ & $0.6941$ & $0.6162$ \\
RP & $54$ &  $1$ & $0.9139$ &      N/A \\
SA & $16$ & $11$ & $0.9139$ &      N/A \\
\hline\hline
\end{tabular}
\end{table}

\subsection{Residential Building Dataset}
\label{sec:realdata:building}

The second dataset was used to study $n=372$ residential condominiums from as many 3- to 9-story buildings constructed between 1993 and 2008 in Tehran, Iran \citep{RA16,RA18}.
Construction cost, sale price, $8$ project physical and financial (PF) variables and $19$ economic variables and indices (EVI) with up to $5$ time lags before the construction were collected on the quarterly basis.
Similar to the analysis in \citet{RA18}, we study how construction cost is influenced by the PF and delayed EVI factors, but exclude the only categorical PF variable, project locality.
Accordingly, we have $d=7+19\times 5=102$ covariates.
We define the variable coding in Table~\ref{tab:building:coding} for the ease of presentation.

Table~\ref{tab:building:compare} and Table~\ref{tab:building:soil} respectively summarize the variable selection and variable importance results of the GA, the RP and the SA.
From the former, we see that the GA-best model gives the best performance (i.e., lowest GIC), and its variable structure agrees with the findings by \citet{RA18}, which suggest that PF and EVI factors (especially 4-quarter delayed ones) be informative.
Moreover, the second column in Table~\ref{tab:building:soil} confirms the relevance of PF-5, PF-7, $1$-quarter delayed EVI-05, $4$-quarter delayed EVI-07 and EVI-13, and $5$-quarter delayed EVI-12.
We also note that the RP- and SA-best models do not consist of sensible variable structures and are significantly worse than the GA-best model ($p$-values $<0.0001$).

\begin{table}[!tb]
\centering\small
\caption{\label{tab:building:compare}
Summary of the best models for the residential building dataset.}
\begin{tabular}{llr}
\hline\hline
Method & Active Variables & \multicolumn{1}{c}{GIC} \\
\hline
GA & PF-5, PF-7, EVI-05-Lag1, EVI-07-Lag4, EVI-12-Lag5, EVI-13-Lag4 & $2571.49$ \\
RP & (None)                                                         & $3788.05$ \\
SA & PF-2, PF-3, PF-4, PF-5, PF-6, PF-7                             & $2699.63$ \\
\hline\hline
\end{tabular}
\end{table}

\begin{table}[!tb]
\centering\small
\caption{\label{tab:building:soil}
SOIL values of the important variables for the residential building dataset.
SOIL values less than $0.05$ are not listed.}
\begin{tabular}{lccc}
\hline\hline
              & \multicolumn{3}{c}{SOIL} \\
\cline{2-4}
Variable Code & GA      & RP      & SA \\
\hline
PF-2          &         &         & $1.000$ \\
PF-3          &         &         & $1.000$ \\
PF-4          &         &         & $1.000$ \\
PF-5          & $1.000$ &         & $1.000$ \\
PF-6          &         &         & $1.000$ \\
PF-7          & $1.000$ &         & $1.000$ \\
EVI-05-Lag1   & $1.000$ &         &         \\
EVI-07-Lag4   & $1.000$ &         &         \\
EVI-12-Lag5   & $1.000$ &         &         \\
EVI-13-Lag4   & $1.000$ &         &         \\
EVI-19-Lag1   &         & $1.000$ &         \\
\hline\hline
\end{tabular}
\end{table}

\begin{figure}[!tb]
\centering
\includegraphics[width=0.5\textwidth]{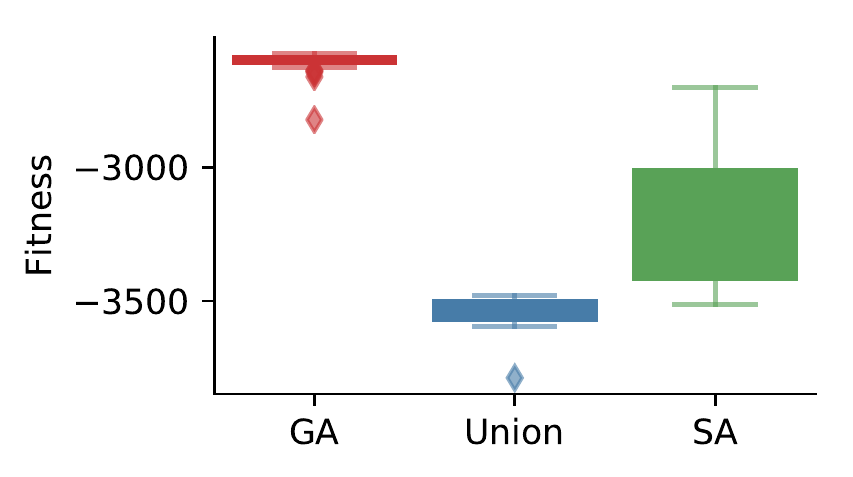}
\caption{\label{fig:building:fitness}
Boxplots of the fitness values of the candidate models for the residential building dataset.}
\end{figure}

Figure~\ref{fig:building:fitness} and Table~\ref{tab:building:mumin} respectively display the boxplots of the fitness values of the candidate models and the multi-model analysis results to evaluate the quality of candidate model sets and model averaging.
The former (Figure~\ref{fig:building:fitness}) suggests that the GA models generally possess higher fitness (i.e., lower GIC) values.
Again, the GA is shown to produce the best candidate model set by having the fittest best model (all $p$-values $<0.0001$) and the highest relative size of $95{\%}$ SMS of $41/48=85.41{\%}$ (compared to approximately $14{\%}$ for the RP and the SA).
In addition to generating the best candidate model set, the GA also results in the lowest RMSE of model averaging using both the GIC-based and AL weighting methods.
These results suggest that good candidate models be helpful in enhancing the performance of multi-model inference.

\begin{table}[!tb]
\centering\small
\caption{\label{tab:building:mumin}
Results of relative size of $95{\%}$ SMSs and model averaging for the residential building dataset.}
\begin{tabular}{lrrrr}
\hline\hline
&&& \multicolumn{2}{r}{RMSE of Model Averaging} \\
\cline{4-5}
Method & $\#(\text{Candidate Models})$ & $\#(\text{Models in}~95{\%}~\text{SMS})$ & GIC-based & AL \\
\hline
GA & $48$ & $41$ &  $27.5553$ & $28.4411$ \\
RP & $11$ &  $3$ & $104.9914$ &       N/A \\
SA & $84$ & $13$ &  $32.7367$ & $32.2841$ \\
\hline\hline
\end{tabular}
\end{table}

\begin{figure}[!tb]
\centering\small
\includegraphics[width=0.95\textwidth]{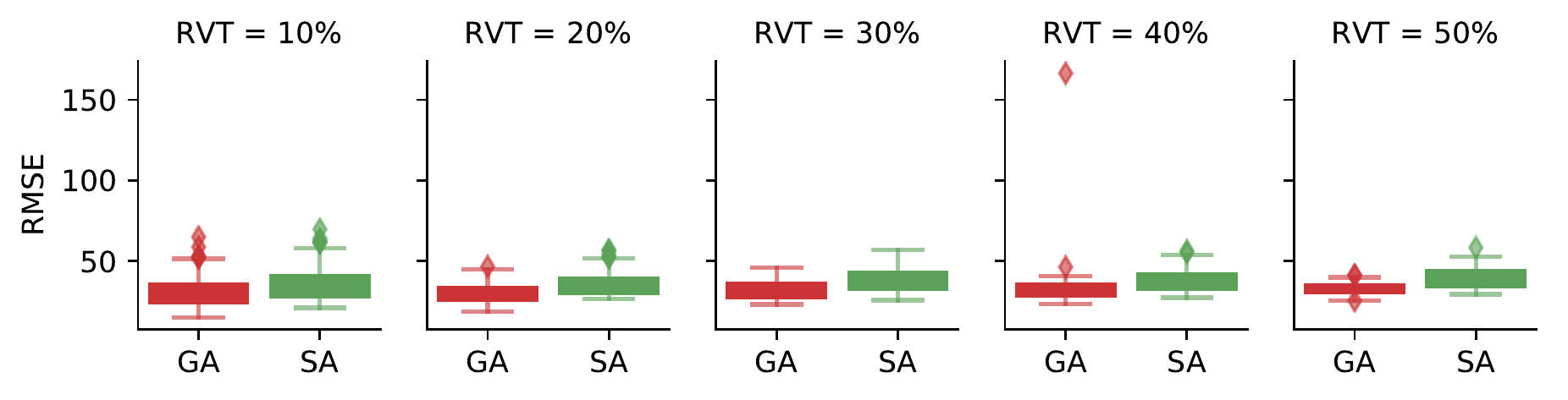}
\hfill
\caption{\label{fig:building:ma1}
Boxplots of RMSE of model averaging using the AL weighting for the residential building dataset.
The RP method failed in weight calculations in all cases and therefore is not shown.}
\end{figure}

To further investigate the predictive performance via model averaging with the AL weighting, we randomly split the dataset using five ratios of validation to training (RVTs) of $10{\%},20{\%},30{\%},40{\%}$ and $50{\%}$.
For each RVT, $100$ randomly selected validation and training datasets were generated by splitting the original dataset, and the boxplots of RMSE are drawn in Figure~\ref{fig:building:ma1}.
In summary, the GA generally results in lower RMSE, suggesting its superior predictive performance.

\section{Discussion}
\label{sec:discussion}

In the end, we propose three future directions.
Firstly, we are interested in developing more implementable algorithms for Theorem~\ref{thm:mcs} to construct the proposed MCS procedure.
Secondly, we believe that incorporating GAs into modern computational tools such as neural networks may produce more powerful statistical inference procedures.
For instance, the {\em deep neuroevolution} developed by the Uber AI Labs uses GAs to train deep reinforcement learning (DRL) models and demonstrates amazing performance on hard DRL benchmarks such as Atari and Humanoid Locomotion \citep[e.g.,][]{Uber1,Uber2,Uber3}; see \url{https://eng.uber.com/deep-neuroevolution/} for a comprehensive introduction.
Lastly, we want to investigate more advanced GA variants (e.g., adaptive GAs \citep[e.g.,][]{T12:IAGA,SX13:IAGA,RG13,LBC15}, the immune GAs \citep[e.g.,][]{JW00,YZ08,ZOMKO14} or the hybrid GAs \citep[e.g.,][]{CCW05,KZ08,CS09,ZSLGC11}) from statistical and machine learning perspectives.

{\linespread{1}\selectfont
\bibliographystyle{asa}
\bibliography{refs}}

\clearpage
\pagenumbering{arabic}
\setcounter{section}{0}
\renewcommand*{\thepage}{S.\arabic{page}}
\renewcommand*{\thesection}{S.\arabic{section}}
\setcounter{table}{0}
\renewcommand{\thetable}{S.\arabic{table}}%
\setcounter{figure}{0}
\renewcommand{\thefigure}{S.\arabic{figure}}%

\makeatletter
\def\@thanks{}
\makeatother

\title{Supplementary material for:\protect\\
  Enhancing Multi-model Inference with Natural Selection}

\date{Purdue University}

\maketitle

The supplementary material is organized as follows:
\begin{itemize}
\item In Section~\ref{supp:theory}, the proofs for Section~\ref{sec:theory} are presented.
\item In Section~\ref{supp:appl}, the proofs for Section~\ref{sec:appl} are presented.
\item In Section~\ref{supp:details}, we present the details of
\begin{itemize}
\item the GIC-based superiority test \eqref{eq:sms:test:sup},
\item the model averaging approach of \citet{AL14}, and
\item the approach of random initial population generation used in Section~\ref{sec:sim:schema}.
\end{itemize}
\item In Section~\ref{supp:sim}, we present additional simulation results for computation time, schema evolution, and variable importance.
\item In Section~\ref{supp:building}, we present the variable coding for the residential building dataset used in Section~\ref{sec:realdata:building}.
\item In Section~\ref{supp:lemmas}, we present technical lemmas.
\end{itemize}

\section{Proofs for Section~\ref{sec:theory}}
\label{supp:theory}

\subsection{Proof of Theorem~\ref{thm:mc}}
\label{supp:proof:thm:mc}

To prove (a), first note that given $u^{*}\in\Psi(t)$ the subsequent generations cannot travel to any state that does not contain $u^{*}$ due to elitism selection.
This means $\bP$ is reducible, and $\cM_{\max}$ is closed (in the sense that $P\big(u^{*}\not\in\Psi(t')\big|u^{*}\in\Psi(t)\big)=0$ for all $t'>t$).

Without loss of generality, there exists square matrices $\bA$ and $\bT$, and a matrix $\bR$ with suitable dimensions such that
\begin{align*}
\bP=\begin{bmatrix}\bA&\bO\\\bR&\bT\end{bmatrix},
\end{align*}
where $\bA$ is a $|\cM_{\max}|\times|\cM_{\max}|$ transition probability submatrix corresponding to the states in $\cM_{\max}$.
According to Lemma~\ref{lemma:R94:thm2} (Theorem~2 of \citet{R94}), it suffices to show that $\bA=[a_{\bu\bv}]_{\bu,\bv\in\cM_{\max}}$ is stochastic and primitive, and $\bR$ and $\bT$ are not zero matrices.

To show $\bA$ is stochastic and primitive, first note that $\bA$ corresponds to the transition probability matrix for the states $\bu\in\cM_{\max}$.
Since any $P\big(\Psi(t+1)\not\in\cM_{\max}\big|\Psi(t)\in\cM_{\max}\big)=0$ for any $t\geq 0$, we must have $\sum_{\bv\in\cM_{\max}}a_{\bu\bv}=1$.
This indicates that $\bA$ is stochastic.

For any fixed-size population $\bu$, the child models generated by selection and crossover operations still belong to $\cM$, and they can be transformed to any other models through the mutation operator with $\pi_{m}\in(0,1)$.
In other words, any model $u\in\bu$ with $u\not=u^{*}$ can be mapped to any $v\in\cM$.
This implies any state in $\cM_{\max}$ can travel to any other state in $\cM_{\max}$ with positive probability.
Accordingly, $\bA$ is positive and thus primitive.

Similar argument yields that $P_{\bu\bu}=P\big(\Psi(t+1)=\bu\big|\Psi(t)=\bu\big)>0$ for all $\bu\in\cM^{K}$, and therefore $\bT$, the transition probability matrix corresponding to the states not in $\cM_{\max}$, is not zero.
Moreover, since the generational best model can only be improved, any model $u$ can be transformed to $u^{*}$ with positive probability due to the mutation operator with $p_{m}\in(0,1)$.
Hence for any $t\geq 0$ we have
\begin{align}
P\big(\Psi(t+1)=\bv\big|\Psi(t)=\bu\big)>0\quad\text{for all}~\bu\not\in\cM_{\max}~\text{and}~\bv\in\cM_{\max}. \label{eq:Rmtx}
\end{align}
Note that the entries of $\bR$ collects all such transition probabilities.
Consequently, it is a positive, and thus nonzero matrix.

The result of (b) is a straightforward consequence of (a).
That is, since $\bpi$ is a distribution over $\cM^{K}$ and $\pi(\bu)=0$ for all $\bu\not\in\cM_{\max}$, we have $\sum_{\bu\in\cM_{\max}}\pi(\bu)=1$.
By the definition of $\cM_{\max}$, it further implies the asymptotic inclusion of the best model as $t\to\infty$.

\subsection{Proof of Theorem~\ref{thm:mcs}}
\label{supp:proof:thm:mcs}

It suffices to show that
\begin{align}
P\big(\Psi(T_{\alpha})\in\cM_{\max}\big)\geq 1-\alpha. \label{eq:mcs:1}
\end{align}
Since the GA with elitism selection satisfies
\begin{align*}
\Big\{\Psi(t)\in\cM_{\max}\Big\}\subset\Big\{\Psi(t+1)\in\cM_{\max}\Big\}\quad\text{for all}~t\geq 0,
\end{align*}
it suffices to show that there exists a positive integer $T_{\alpha}$ such that
\begin{align}
P\left(\bigcup_{t=1}^{T_{\alpha}}\Big\{\Psi(t)\in\cM_{\max}\Big\}\,\middle|\,\Psi(0)=\bu\right)\geq 1-\alpha\quad\text{for any}~\bu\in\cM^{K}. \label{eq:proof:mcs}
\end{align}

Let
\begin{align*}
P_{\bu\cM_{\max}}=\sum_{\bv\in\cM_{\max}}P\big(\Psi(t+1)=\bv\big|\Psi(t)=\bu\big)
\end{align*}
denotes the total probability that a population $\bu$ is transmitted into any population with the best solution in one iteration.
According to \eqref{eq:Rmtx}, define
\begin{align}
\xi:=\inf_{\bu\in\cM^{K}}P_{\bu\cM_{\max}}=\inf_{\bu\in\cM^{K}}\sum_{\bv\in\cM_{\max}}P\big(\Psi(t+1)=\bv\big|\Psi(t)=\bu\big)>0. \label{eq:mcs:xi}
\end{align}
Note that, for all $\bu\in\cM^{K}$ and positive integer $t$,
\begin{align*}
1-\xi\geq P\Big(\Psi(t)\not\in\cM_{\max}\,\Big|\,\Psi(0)=\bu\Big)=\E\Big[\dsone\big(\Psi(t)\not\in\cM_{\max}\,\big|\,\Psi(0)=\bu\big)\Big].
\end{align*}

It then holds, for any $\bu\in\cM^{K}$ and positive integer $T$,
\begin{align*}
&P\left(\bigcap_{t=1}^{T}\Big\{\Psi(t)\not\in\cM_{\max}\Big\}\,\middle|\,\Psi(0)=\bu\right)\\
&\qquad=\E\left[\dsone\left(\bigcap_{t=1}^{T}\Big\{\Psi(t)\not\in\cM_{\max}\Big\}\right)\middle|\,\Psi(0)=\bu\right]\\
&\qquad=\E\left[\prod_{t=1}^{T}\dsone\big(\Psi(t)\not\in\cM_{\max}\big)\middle|\,\Psi(0)=\bu\right]\\
&\qquad=\E\left[\E\Big[\dsone\big(\Psi(T)\not\in\cM_{\max}\big)\,\Big|\,\Psi(T-1)\Big]\prod_{t=1}^{T-1}\dsone\big(\Psi(t)\not\in\cM_{\max}\big)\middle|\,\Psi(0)=\bu\right]\\
&\qquad\leq(1-\xi)\E\left[\prod_{t=1}^{T-1}\dsone\big(\Psi(t)\not\in\cM_{\max}\big)\middle|\,\Psi(0)=\bu\right],
\end{align*}
where the third equality is due to the Markov property.
By keeping doing this we obtain
\begin{align*}
P\left(\bigcap_{t=1}^{T}\Big\{\Psi(t)\not\in\cM_{\max}\Big\}\,\middle|\,\Psi(0)=\bu\right)\leq(1-\xi)^{T}.
\end{align*}
Since $\xi\in(0,1)$, there exists a positive integer $T_{\alpha}$ such that $(1-\xi)^{T_{\alpha}}\leq\alpha<(1-\xi)^{T_{\alpha}-1}$.
Accordingly, the desired confidence statement \eqref{eq:proof:mcs} follows.

\subsection{Proof of Theorem~\ref{thm:schema}}
\label{supp:proof:thm:schema}

We first characterize the individual probabilities caused by the selection, crossover and mutation operations.
Firstly, it is obvious that the probability that models $u^{k}$ and $u^{l}$ are selected is $w_{k}w_{l}$.
Secondly, the probability that the uniform mutation operation transforms a given model $v$ into a solution that matches $H$ is $\pi_{m}^{\delta(v,H)}(1-\pi_{m})^{\ord(H)-\delta(v,H)}$.

Finally, we discuss the effect of the uniform crossover operation, given two parent models $u^{k}$ and $u^{l}$ are selected.
Due to the mechanism of the uniform crossover, all possible child models has equal probabilities to be generated.
This allows us to focus on the fixed positions of $H$.
Note that it is possible that $u^{k}$ and $u^{l}$ can never generate a child model that is a solution that matches $H$.
Therefore, we define
\begin{align*}
h_{kl}=\big|\{j:H_{j}\not=\ast,u^{k}_{j}=u^{l}_{j}\not=H_{j}\}\big|
\end{align*}
as the minimum $\delta(v,H)$ among all the child models $v$ produced by the uniform crossover with parent models $u^{k}$ and $u^{l}$.
Now, suppose $v$ is a model generated through uniform crossover with $u^{k}$ and $u^{l}$, we have
\begin{align*}
P\big(\delta_{H}(v,H)=h+h_{kl}\,\big|\,\text{parent models}~u^{k},u^{l}\big)=\frac{\binom{\delta_{H}(u^{k},u^{l})-h_{kl}}{h}}{2^{\delta_{H}(u^{k},u^{l})}}\quad\text{for}~h=0,1,\ldots,\delta_{H}(u^{k},u^{l})-h_{kl}.
\end{align*}

Accordingly, a general form of $\alpha(H,t)$ can be written by
\begin{align}
\alpha(H,t)&=\sum_{k,l:u^{k},u^{l}\in\Psi(t)}w_{k}w_{l}\left[\sum_{h=0}^{\delta_{H}(u^{k},u^{l})-h_{kl}}\frac{\binom{\delta_{H}(u^{k},u^{l})-h_{kl}}{h}}{2^{\delta_{H}(u^{k},u^{l})-h_{kl}}}\pi_{m}^{h+h_{kl}}(1-\pi_{m})^{\ord(H)-h-h_{kl}}\right] \notag\\
&=\sum_{k,l:u^{k},u^{l}\in\Psi(t)}w_{k}w_{l}\frac{\pi_{m}^{h_{kl}}(1-\pi_{m})^{\ord(H)-\delta_{H}(u^{k},u^{l})}}{2^{\delta_{H}(u^{k},u^{l})-h_{kl}}} \notag\\
&\hspace{8em}\times\left[\sum_{h=0}^{\delta_{H}(u^{k},u^{l})-h_{kl}}\binom{\delta_{H}(u^{k},u^{l})-h_{kl}}{h}\pi_{m}^{h+h_{kl}}(1-\pi_{m})^{\delta_{H}(u^{k},u^{l})-h-h_{kl}}\right] \notag\\
&=\sum_{k,l:u^{k},u^{l}\in\Psi(t)}w_{k}w_{l}\frac{\pi_{m}^{h_{kl}}(1-\pi_{m})^{\ord(H)-\delta_{H}(u^{k},u^{l})}}{2^{\delta_{H}(u^{k},u^{l})-h_{kl}}} \label{eq:schema:all}
\end{align}
Note that on the right hand side of \eqref{eq:schema:all}, the summation can be tore apart to three cases based on whether the parents are solutions that match $H$.
That is,
\begin{align}
\alpha(H,t)=P\big(\text{Case 1}\big)+P\big(\text{Case 2}\big)+P\big(\text{Case 3}\big), \label{eq:schema:cases}
\end{align}
where Cases 1, 2 and 3 refer to the events that the final child model after crossover and mutation is a solution that matches $H$ given that
\begin{enumerate}
\item both parents match $H$ (i.e., $k,$ such that $u^{k},u^{l}\in{\uparrow}(H)$),
\item only one of the parents matches $H$ (i.e., $k$ such that $u^{k}\in{\uparrow}(H)$ and $l:u^{l}\not\in{\uparrow}(H)$), and
\item neither of the parents matches $H$ (i.e., $k,l$ such that $u^{k},u^{l}\not\in{\uparrow}(H)$),
\end{enumerate}
respectively.

For Case 1, since both parents belong to ${\uparrow}(H)$, it follows that $\delta_{H}(u^{k},u^{l})=0$ and $h_{kl}=0$, and hence
\begin{align}
P\big(\text{Case 1}\big)&=\sum_{k,l:u^{k},u^{l}\in{\uparrow}(H)}w_{k}w_{l}(1-\pi_{m})^{\ord(H)} \notag\\
&=\left(\sum_{k:u^{k}\in{\uparrow}(H)}w_{k}\right)^{2}(1-\pi_{m})^{\ord(H)} \notag\\
&=\alpha_{sel}(H,t)^{2}(1-\pi_{m})^{\ord(H)}. \label{eq:prob:case1}
\end{align}

For Case 2, since one of the parents matches $H$, it holds $h_{kl}=0$ and $\delta_{H}(u^{k},u^{l})=\delta(u^{l},H)$.
It then holds that
\begin{align}
P\big(\text{Case 2}\big)&=\sum_{\substack{k:u^{k}\in{\uparrow}(H)\\l:u^{l}\not\in{\uparrow}(H)}}w_{k}w_{l}\frac{(1-\pi_{m})^{\ord(H)}}{\big[2(1-\pi_{m})\big]^{\delta(u^{l},H)}}\notag\\
&=\alpha_{sel}(H,t)\sum_{l:u^{l}\not\in{\uparrow}(H)}w_{l}\frac{(1-\pi_{m})^{\ord(H)}}{\big[2(1-\pi_{m})\big]^{\delta(u^{l},H)}}. \label{eq:prob:case2}
\end{align}

For Case 3, there seems no simplification available, and therefore we have
\begin{align}
P\big(\text{Case 3}\big)=\sum_{k,l:u^{k},u^{l}\not\in{\uparrow}(H)}w_{k}w_{l}\frac{(2\pi_{m})^{h_{kl}}(1-\pi_{m})^{\ord(H)}}{\big[2(1-\pi_{m})\big]^{\delta_{H}(u^{k},u^{l})}}. \label{eq:prob:case3}
\end{align}
The proof is then complete by plugging \eqref{eq:prob:case1}, \eqref{eq:prob:case2} and \eqref{eq:prob:case3} into \eqref{eq:schema:cases}.

\subsection{Proof of Corollary~\ref{cor:schema:lower}}
\label{supp:proof:cor:schema:lower}

First note that $2(1-\pi_{m})>1$ since $\pi_{m}\leq 0.5$.
Since $\delta(u^{l},H)\leq\ord(H)$ for all $u^{l}\not\in{\uparrow}(H)$, it follows that
\begin{align*}
P(\text{Case 2})&=\alpha_{sel}(H,t)\sum_{l:u^{l}\not\in{\uparrow}(H)}w_{l}\frac{(1-\pi_{m})^{\ord(H)}}{\big[2(1-\pi_{m})\big]^{\delta(u^{l},H)}}\\
&\geq 2^{-\ord(H)}\alpha_{sel}(H,t)\sum_{l:u^{l}\not\in{\uparrow}(H)}w_{l}\\
&=2^{-\ord(H)}\alpha_{sel}(H,t)\big[1-\alpha_{sel}(H,t)\big].
\end{align*}
Similarly, since $2\pi_{m}<1$, $h_{kl}\leq\ord(H)$ and $\delta_{H}(u^{k},u^{l})\leq\ord(H)$ for all $u^{k},u^{l}\not\in{\uparrow}(H)$, we have
\begin{align*}
P(\text{Case 3})&=\sum_{k,l:u^{k},u^{l}\not\in{\uparrow}(H)}w_{k}w_{l}\frac{(2\pi_{m})^{h_{kl}}(1-\pi_{m})^{\ord(H)}}{\big[2(1-\pi_{m})\big]^{\delta_{H}(u^{k},u^{l})}}\\
&\geq\pi_{m}^{\ord(H)}\big[1-\alpha_{sel}(H,t)\big]^{2}.
\end{align*}
Accordingly, we have the desired result \eqref{eq:schema:lower}.

\section{Proof for Section~\ref{sec:appl}}
\label{supp:appl}

\subsection{Proof of Lemma~\ref{lemma:GIC}}

Without loss of generality, let $u^{0}$ denote the binary sequence with first $s$ genes active and the rest inactive and $\sigma^{2}=1$.
Recall that $\bX_{u}$ denotes the submatrix of $\bX$ subject to the active variable indices in $u$.
Let $\bH_{u}=\bX_{u}(\bX_{u}^{\top}\bX_{u})^{-1}\bX_{u}^{\top}$ the projection matrix of the submatrix $\bX_{u}$.

We first consider the case $u\not\supseteq u^{0}$, i.e., model $u$ misses at least one relevant variable.
We can write
\begin{align*}
\GIC(u)-\GIC(u^{0})&=n\log\left(1+\frac{\RSS(u)-\RSS(u^{0})}{\RSS(u^{0})}\right)+\kappa_{n}\big(|u|-s\big)\\
&\geq n\log\left(1+\frac{\RSS(u)-\RSS(u^{0})}{\RSS(u^{0})}\right)-\kappa_{n}s.
\end{align*}
Note that
\begin{align}
\RSS(u^{0})=\bY^{\top}(\bI-\bH_{u^{0}})\bY=\beps^{\top}(\bI-\bH_{u^{0}})\beps=\sum_{i=1}^{d-s}Z_{i}^{2}=n(1+o(1)), \label{eq:RSS0}
\end{align}
where the $Z_{i}$ are independent $\Norm(0,1)$ variables, and
\begin{align}
\RSS(u)-\RSS(u^{0})&=\bY^{\top}(\bI-\bH_{u})\bY-\beps^{\top}(\bI-\bH_{u^{0}})\beps\notag\\
&=\bmu^{\top}(\bI-\bH_{u^{0}})\bmu+2\bmu^{\top}(\bI-\bH_{u})\beps-\beps^{\top}\bH_{u}\beps+\beps^{\top}\bH_{u^{0}}\beps, \label{eq:RSSdiff}
\end{align}
where $\bmu=\bX_{u^{0}}\bbeta^{0}_{u^{0}}$.
By Condition~(A2), uniformly over $u$ with $|u|\leq\tilde{s}$, it holds
\begin{align}
\min_{u\in\cM_{\tilde{s}}-\{u^{0}\}}\bmu^{\top}(\bI-\bH_{u})\bmu\geq C_{2}n. \label{eq:RSSdiff:1}
\end{align}

Write
\begin{align*}
\bmu^{\top}(\bI-\bH_{u})\beps=\sqrt{\bmu^{\top}(\bI-\bH_{u})\bmu}Z_{u},\quad\text{where}~Z_{u}=\frac{\bmu^{\top}(\bI-\bH_{u})\beps}{\sqrt{\bmu^{\top}(\bI-\bH_{u})\bmu}}\sim\Norm(0,1).
\end{align*}
Note that for any model $u$ with $|u|\leq\tilde{s}$, there exists a positive constant $C_{3}$ such that
\begin{align*}
P\big(|Z_{u}|>t\big)=C_{3}\exp\left(-\frac{t^{2}}{2}\right).
\end{align*}
By the union bound, it follows that
\begin{align*}
P\left(\max_{u\in\cM_{\tilde{s}}-\{u^{0}\},u\not\supseteq u^{0}}|Z_{u}|>t\right)\leq\sum_{u\in\cM_{\tilde{s}}-\{u^{0}\},u\not\supseteq u^{0}}P\Big(|Z_{u}|>t\Big)\leq 2^{\tilde{s}}C_{3}\exp\left(-\frac{t^{2}}{2}\right).
\end{align*}
Let $t=\sqrt{2s\log d}$, we arrive at
\begin{align*}
P\left(\max_{u\in\cM_{\tilde{s}}-\{u^{0}\},u\not\supseteq u^{0}}|Z_{u}|>t\right)\leq C_{3}\left(\frac{2}{d}\right)^{\tilde{s}}\to 0
\end{align*}
as $n\to\infty$.
Accordingly,
\begin{align*}
\max_{u\in\cM_{\tilde{s}}-\{u^{0}\},u\not\supseteq u^{0}}|Z_{u}|=O_{P}\big(\sqrt{\tilde{s}\log d}\big)=o_{P}\big(\sqrt{n}\big),
\end{align*}
and therefore we have
\begin{align}
\max_{u\in\cM_{\tilde{s}}-\{u^{0}\},u\not\supseteq u^{0}}\bmu^{\top}(\bI-\bH_{u})\beps&\leq\sqrt{\bmu^{\top}(\bI-\bH_{u})\bmu}\max_{u\in\cM_{\tilde{s}}-\{u^{0}\},u\not\supseteq u^{0}}Z_{u}\notag\\
&=\sqrt{\bmu^{\top}(\bI-\bH_{u})\bmu}~o_{P}(\sqrt{n})\notag\\
&=o_{P}\left(\sqrt{\bmu^{\top}(\bI-\bH_{u})\bmu}\right). \label{eq:RSSdiff:2}
\end{align}

Now we deal with the last two terms in \eqref{eq:RSSdiff}.
Note that we can write
\begin{align*}
\beps^{\top}\bH_{u}\beps=\sum_{i=1}^{|u|}Z_{i}^{2}\sim\chi^{2}_{|u|},
\end{align*}
where $Z_{i}$ are some independent $\Norm(0,1)$ variables.
By the union bound, it then holds
\begin{align*}
P\left(\max_{u\in\cM_{\tilde{s}}-\{u^{0}\},u\not\supseteq u^{0}}\beps^{\top}\bH_{u}\beps>t\right)\leq\sum_{j=1}^{\tilde{s}}\binom{d}{j}P\big(\chi^{2}_{j}>t\big)\leq d^{\tilde{s}}P\big(\chi^{2}_{\tilde{s}}>t\big).
\end{align*}
It is east to see that (see, for example, \citet{Y99})
\begin{align*}
P\big(\chi^{2}_{\tilde{s}}>t\big)\leq\exp\left(-\frac{t-\tilde{s}}{2}\right)\left(\frac{t}{\tilde{s}}\right)^{\tilde{s}/2}.
\end{align*}
Let $t=3s\log d$, we arrive at
\begin{align*}
P\left(\max_{u\in\cM_{\tilde{s}}-\{u^{0}\},u\not\supseteq u^{0}}\beps^{\top}\bH_{u}\beps>t\right)&\leq\left(\frac{e\log d}{d}\right)^{\tilde{s}/2}\to 0
\end{align*}
as $n\to\infty$.
Consequently, we have
\begin{align}
\max_{u\in\cM_{\tilde{s}}-\{u^{0}\},u\not\supseteq u^{0}}\beps^{\top}\bH_{u}\beps=O_{P}\big(\tilde{s}\log d\big)=o_{P}(n). \label{eq:RSSdiff:3}
\end{align}
Similarly,
\begin{align}
\beps^{\top}\bH_{u^{0}}\beps=o_{P}(n). \label{eq:RSSdiff:4}
\end{align}

By \eqref{eq:RSSdiff:1}, \eqref{eq:RSSdiff:2}, \eqref{eq:RSSdiff:3} and \eqref{eq:RSSdiff:4}, it is easy to see that $\RSS(u)-\RSS(u^{0})$ is dominated by $\bmu^{\top}(\bI-\bH_{u^{0}})\bmu$.
Coupled with \eqref{eq:RSS0}, there is a positive constant $C_{4}$ such that
\begin{align*}
\log\left(1+\frac{\RSS(u)-\RSS(u^{0})}{\RSS(u^{0})}\right)\geq\log(1+C_{4})
\end{align*}
in probability.
Since $\kappa_{n}=o(n)$, we conclude that
\begin{align}
\min_{u\in\cM_{\tilde{s}}-\{u^{0}\},u\not\supseteq u^{0}}\GIC(u)-\GIC(u^{0})\geq n\log(1+C_{4})-\kappa_{n}s>0 \label{eq:selconsistent:1}
\end{align}
as $n\to\infty$.

Now we consider the case $u\supseteq u^{0}$ but $u\not=u^{0}$.
Since $(\bI-\bH_{u})\bX_{u^{0}}=\bO$, we have $\bY^{\top}(\bI-\bH_{u})\bY=\beps^{\top}(\bI-\bH_{u})\beps$ and
\begin{align*}
\RSS(u^{0})-\RSS(u)=\beps^{\top}(\bH_{u}-\bH_{u^{0}})\beps=\sum_{i=1}^{|u|-s}Z_{u,i}^{2}\sim\chi^{2}_{|u|-s},
\end{align*}
where $Z_{u,i}$ are some independent $\Norm(0,1)$ variables depending on $u$.
By the union bound we have
\begin{align*}
P\left(\min_{u\in\cM_{\tilde{s}}-\{u^{0}\},u\supseteq u^{0}}\GIC(u)-\GIC(u^{0})\leq 0\right)&\leq\sum_{u\in\cM_{\tilde{s}}-\{u^{0}\},u\supseteq u^{0}}P\big(\RSS(u^{0})-\RSS(u)\geq\kappa_{n}(|u|-s)\big)\\
&=\sum_{u\in\cM_{\tilde{s}}-\{u^{0}\},u\supseteq u^{0}}P\left(\beps^{\top}(\bH_{u}-\bH_{u^{0}})\beps\geq\kappa_{n}(|u|-s)\right)\\
&\leq\sum_{u\in\cM_{\tilde{s}}-\{u^{0}\},u\supseteq u^{0}}\big[\kappa_{n}\exp(1-\kappa_{n})\big]^{\frac{|u|-s}{2}}\\
&=\sum_{j=s+1}^{\tilde{s}}\binom{d-s}{j-s}\big[\kappa_{n}\exp(1-\kappa_{n})\big]^{\frac{j-s}{2}}\\
&\leq\sum_{m=0}^{d-s}\binom{d-s}{m}\big[\kappa_{n}\exp(1-\kappa_{n})\big]^{\frac{m}{2}}-1\\
&=\left(1+\sqrt{\frac{e\kappa_{n}}{\exp\kappa_{n}}}\right)^{d-s}-1\to 0\quad\text{as}~n\to\infty,
\end{align*}
where the second inequality follows from the sharp deviation bound on the $\chi^{2}$ distribution (see Lemma~3 of \citet{SIS}).
Hence we have
\begin{align}
\min_{u\in\cM_{\tilde{s}}-\{u^{0}\},u\supseteq u^{0}}\GIC(u)-\GIC(u^{0})>0 \label{eq:selconsistent:2}
\end{align}
with probability tending to $1$.
Accordingly, the desired result \eqref{eq:selconsistent} follows from \eqref{eq:selconsistent:1} and \eqref{eq:selconsistent:2}.

\subsection{Proof of Proposition~\ref{prop:vs}}
\label{supp:proof:prop:vs}

From Lemma~\ref{lemma:GIC} we know that the true model $u^{0}$ is the best model in the model space $\cM_{\tilde{s}}$ with probability tending to $1$.
Along with Theorem~\ref{thm:mc} (b) we have
\begin{align*}
\lim_{t\to\infty}\lim_{n\to\infty}P\big(u^{0}=u^{*}\in\Psi_{\tilde{s}}(t)\big)=1.
\end{align*}
By the definition of $\what{u}(t)$, we arrive at
\begin{align*}
\lim_{t\to\infty}\lim_{n\to\infty}P\left(\what{u}(t)=u^{0}\right)=1.
\end{align*}
This completes the proof.

\subsection{Proof of Proposition~\ref{prop:sms}}

By the construction of the $\cA_{\alpha}(t)$, we have
\begin{align*}
\lim_{n\to\infty}P\big(u\in\cA_{\alpha}(t)\big)\geq 1-\alpha
\end{align*}
for all $u\in\Psi_{\tilde{s}}(t)-\{\what{u}(t)\}$ with $H_{0,u}$ not rejected and any $t\geq 0$.
Along with Proposition~\ref{prop:vs}, which ensures that $\lim_{t\to\infty}\lim_{n\to\infty}P\left(\what{u}(t)=u^{0}\right)=1$, the desired result then holds.

\section{Details of the Auxiliary Methods}
\label{supp:details}

\subsection{GIC-Based Superiority Test}
\label{supp:V89:sup}

A natural test statistic for the GIC-based superiority test \eqref{eq:sms:test:sup} can be derived based on the difference of the GIC values of models $u$ and $u^{\#}$.
Note that the first term in the GIC \eqref{eq:GIC} comes from simplifying the log likelihood with Gaussian noise.
That is, the general form for GIC can be written as
\begin{align*}
\GIC(u)=-2\log L(\what{\bbeta}_{u};\bX,\bY)+\kappa_{n}|u|,
\end{align*}
where $L(\bbeta_{u};\bX,\bY)$ is the likelihood function of model $u$ evaluated at $\bbeta_{u}$ given data $(\bX,\bY)$, and $\what{\bbeta}_{u}=\big(\bX_{u}^{\top}\bX_{u}\big)^{-1}\bX_{u}\bY$ for any model $u\in\cM$ with $|u|<n$.
As a result, we write
\begin{align*}
\GIC(u)-\GIC(u^{\#})=\big(|u|-|u^{\#}|\big)\kappa_{n}-2\log\frac{L(\what{\bbeta}_{u};\bX,\bY)}{L(\what{\bbeta}_{u^{\#}};\bX,\bY)}.
\end{align*}
Note that the first term on the R.H.S. is merely a constant and the sampling variation comes only from the second term.
When $u$ and $u^{\#}$ are distinguishable (i.e., $H^{dis}_{0,u}$ in \eqref{eq:sms:test:dis} is rejected), \citet{V89} showed that the normalized log likelihood ratio
\begin{align}
n^{-1/2}\log\frac{L(\what{\bbeta}_{u};\bX,\bY)}{L(\what{\bbeta}_{u^{\#}};\bX,\bY)}\Longrightarrow\Norm(0,\omega_{u}^{2}), \label{eq:V89:LR}
\end{align}
where
\begin{align*}
\omega_{u}^{2}=\Var\left(\log\frac{L(\bbeta^{0}_{u};\bsfX,\sfY)}{L(\bbeta^{0}_{u^{\#}};\bsfX,\sfY)}\right),
\end{align*}
denotes the population variance of the log likelihood ratio of $u$ and $u^{\#}$, $\bbeta^{0}_{u}$ is the true regression coefficient under model $u$, and $\bsfX$ and $\sfY$ are the population counterparts of the design vector and the response scalar, respectively.
Accordingly, under $H^{sup}_{0,u}$, the result \eqref{eq:V89:LR} can be used to show that
\begin{align}
n^{-1/2}\big[\GIC(u)-\GIC(u^{\#})\big]&=n^{-1/2}\left[\big(|u|-|u^{\#}|\big)\kappa_{n}-2\log\frac{L(\what{\bbeta}_{u};\bX,\bY)}{L(\what{\bbeta}_{u^{\#}};\bX,\bY)}\right] \notag\\
&\Longrightarrow\Norm(0,4\omega_{u}^{2}). \label{eq:GICdiff:Norm}
\end{align}
In practice, we plug-in a consistent estimate of $\omega_{u}^{2}$, denoted by $\what{\omega}_{u}^{2}$ (see \citet{V89} for the formula), into \eqref{eq:GICdiff:Norm} to perform the test.
Accordingly, we reject $H^{\sup}_{0,u}$ if
\begin{align*}
\GIC(u)-\GIC(u^{\#})>2z_{1-\alpha}\what{\omega}_{u}\sqrt{n},
\end{align*}
where $z_{1-\alpha}$ is the $(1-\alpha)$-quantile of standard normal distribution, and the value of $\what{\omega}_{u}^{2}$ can be extracted from the {\sf R} package {\bf nonnest2} \citep{nonnest2} when implementing the distinguishability test \eqref{eq:sms:test:dis}.

\subsection{Model Averaging Approach of \texorpdfstring{\citet{AL14}}{Lg}}
\label{supp:AL14}

Given a candidate model set $\Psi=\{u^{1},\ldots,u^{K}\}$, let $\bD_{k}$ be a $n\times n$ diagonal matrix with the $l$-th element being $(1-h_{kl})^{-1}$, where $h_{kl}$ is the $l$-th diagonal element of the hat matrix $\bH_{u^{k}}=\bX_{u^{k}}\big(\bX_{u^{k}}^{\top}\bX_{u^{k}}\big)^{-1}\bX_{u^{k}}^{\top}$, and $\wtil{\bH}_{k}=\bD_{k}(\bH_{u^{k}}-\bI)+\bI$.
Following \citet{AL14}, the $K$-dimensional weight vector $\bw=(w_{1},\ldots,w_{K})^{\top}$ can be computed by
\begin{align}
\what{\bw}=\argmin_{\bw\in[0,1]^{K}}\left(\bY^{\top}\bY-2\bw^{\top}\ba+\bw^{\top}\bB\bw\right), \label{eq:AL14}
\end{align}
where $\ba=(a_{1},\ldots,a_{K})^{\top}$ with $a_{k}=\bY^{\top}\wtil{\bH}_{k}\bY$, and $\bB$ is a $K\times K$ matrix with the $(k,j)$-th element $B_{kl}=\bY^{\top}\wtil{\bH}_{k}^{\top}\wtil{\bH}_{l}\bY$.
Note that the common constraint $\sum_{k=1}^{K}w_{k}=1$ for model weights does not necessarily to be imposed.
In fact, \citet{AL14} show that their weighting approach leads to the smallest possible estimation error of the model averaging predictor \eqref{eq:ma} without the constraint.

\subsection{A Variable Association Measure Assisted Approach for Generating the Initial Population}
\label{supp:initial}

Given variable association measures $\gamma_{j},j=1,\ldots,d$ (e.g., the marginal correlation learning $\big|\what{\Cor}(\bX_{j},\bY)\big|$ \citep{SIS} or the HOLP $\big|\bX_{j}(\bX\bX^{\top})^{-1}\bY\big|$ \citep[][available only for $d\geq n$]{HOLP}), we introduce an approach to randomly generate the initial population $\{u^{0},\ldots,u^{0}_{K}\}$ for the GA as follows.
\begin{description}
\item[Step 1:] Assign the model sizes $|u^{0}_{k}|,k=1,\ldots,K$, by generating $K$ independent\\ $\HyperGeom\big(6\min(n,d),2\min(n,d),\min(n,d)\big)$ random variables, where $\HyperGeom(N,M,n)$ denotes the hypergeometric distribution with the probability mass function
\begin{align*}
P\big(\HyperGeom(N,M,n)=m\big)=\dfrac{\binom{M}{m}\binom{N-M}{n-m}}{\binom{N}{n}},m=\min(0,n+M-N),\ldots,\min(n,M).
\end{align*}
\item[Step 2:] For $k=1,\ldots,K$, the active positions of $u^{0}_{k}$ are determined by randomly selecting $|u^{0}_{k}|$ numbers from $[d]$ without replacement according to the probability distribution $\big\{\gamma_{j}/\sum_{l=1}^{d}\gamma_{l}\big\}_{j=1,\ldots,d}$.
\end{description}
This approach ensures the model sizes are around $\min(n,d)/3$ and never exceed $\min(n,d)$.
Moreover, by making use of the variable association measures $\gamma_{j}$, the resulting models are likely to contains the true signals so that their performance are by no means poor.

\section{Supplementary Simulation Results}
\label{supp:sim}

\subsection{Computation Time}
\label{supp:sim:time}

\begin{figure}[!tb]
\centering
\includegraphics[width=0.95\textwidth]{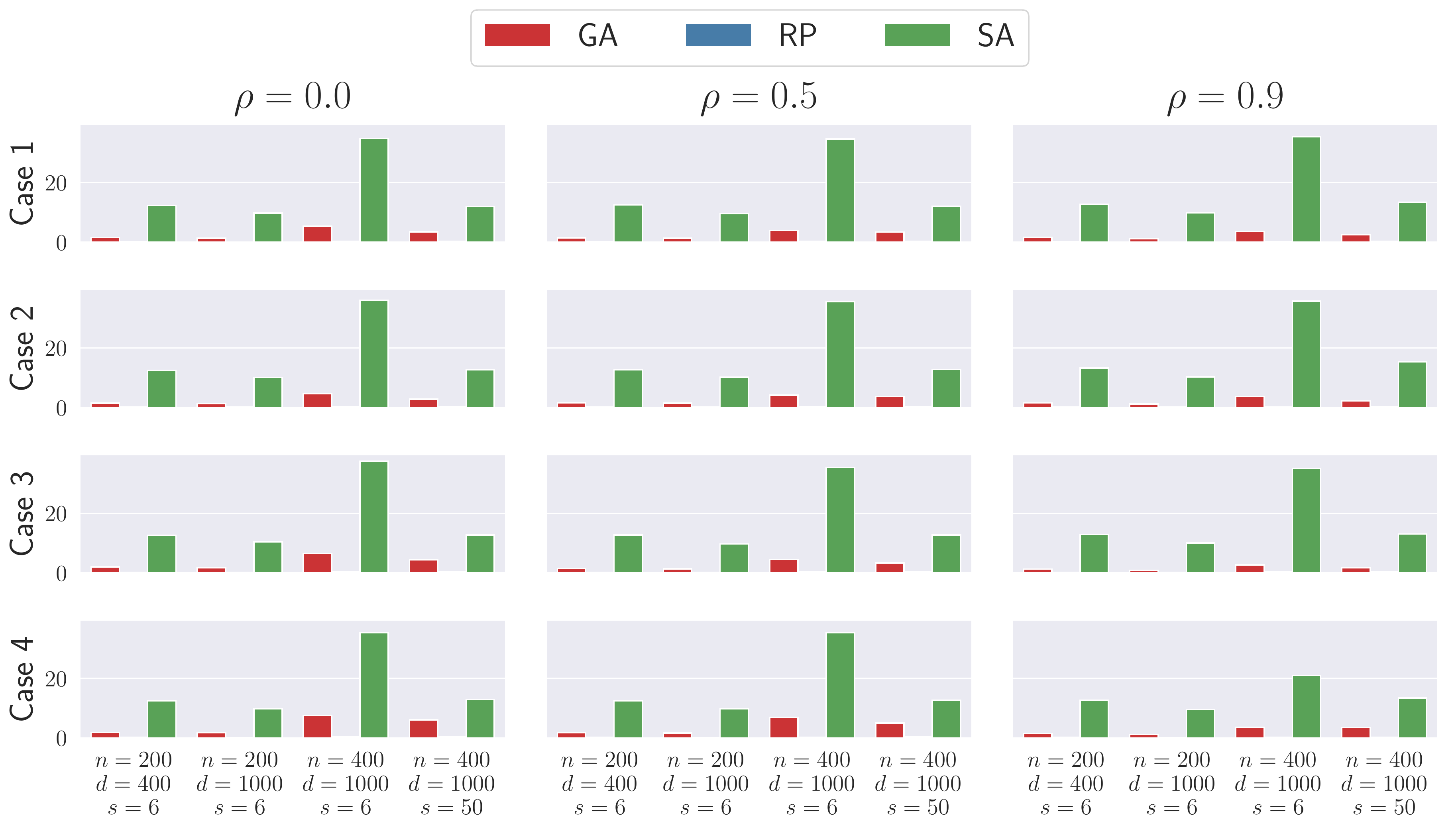}
\includegraphics[width=0.95\textwidth]{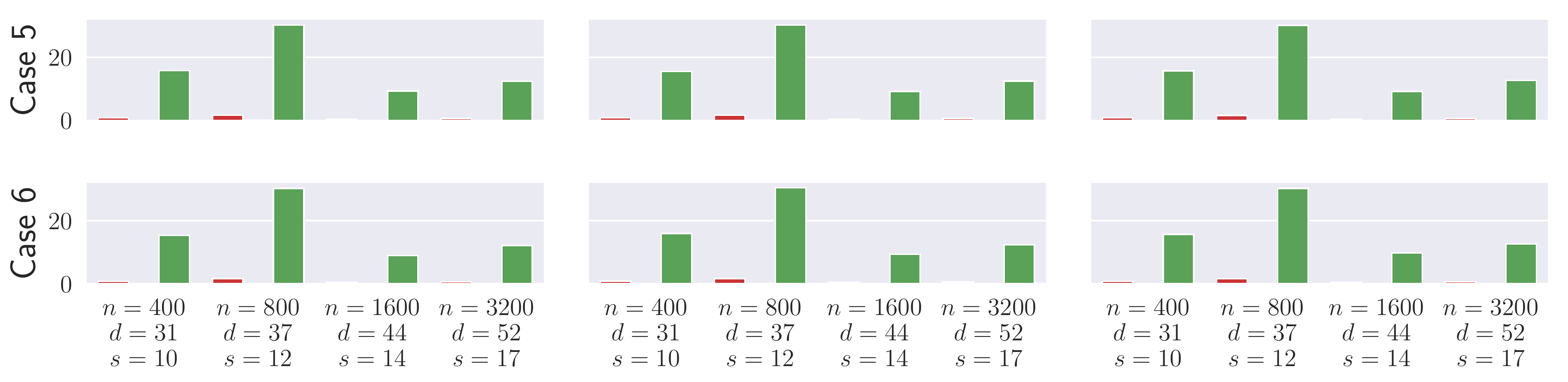}
\caption{\label{fig:time}
Computation time.
(The RP method is too fast to be visualized.)}
\end{figure}

Figure~\ref{fig:time} displays the bar graph for the averaged computation time for implementing the three methods.

\subsection{Schema Evolution}
\label{supp:sim:schema}

\begin{figure}[!tb]
\centering
\includegraphics[width=\textwidth,page=1]{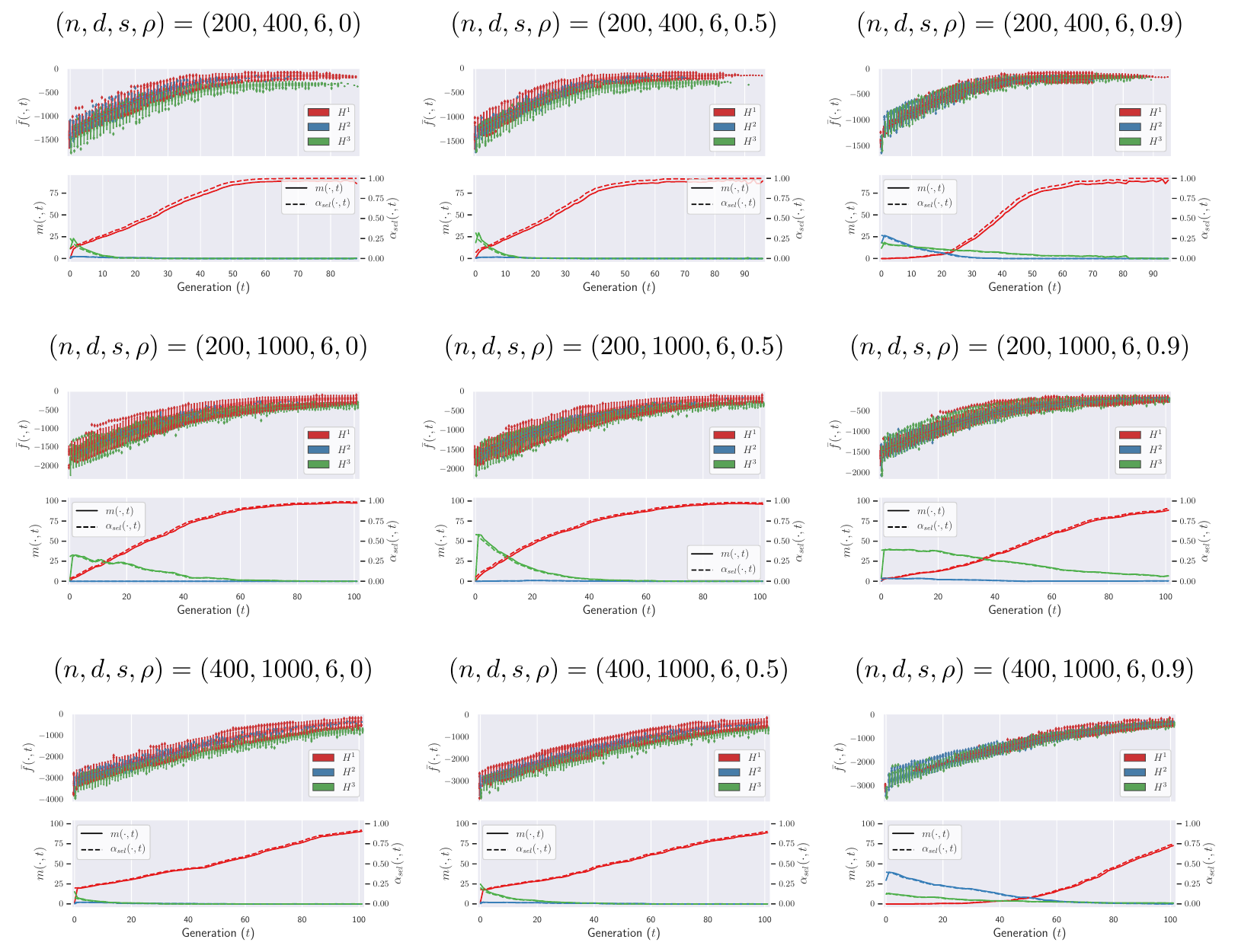}
\caption{\label{fig:schema:case1}
(Case 1) Schema performance (upper panel) and evolution (lower panel) obtained from $500$ simulation runs.}
\end{figure}

\begin{figure}[!tb]
\centering
\includegraphics[width=\textwidth,page=2]{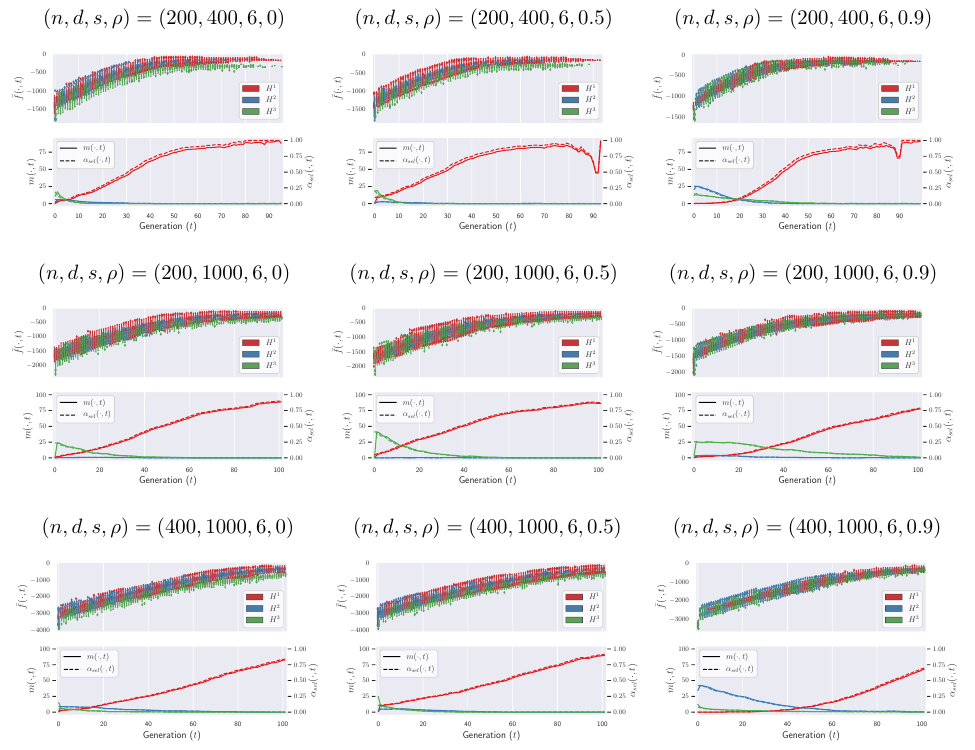}
\caption{\label{fig:schema:case2}
(Case 2) Schema performance (upper panel) and evolution (lower panel) obtained from $500$ simulation runs.}
\end{figure}

\begin{figure}[!tb]
\centering
\includegraphics[width=\textwidth,page=3]{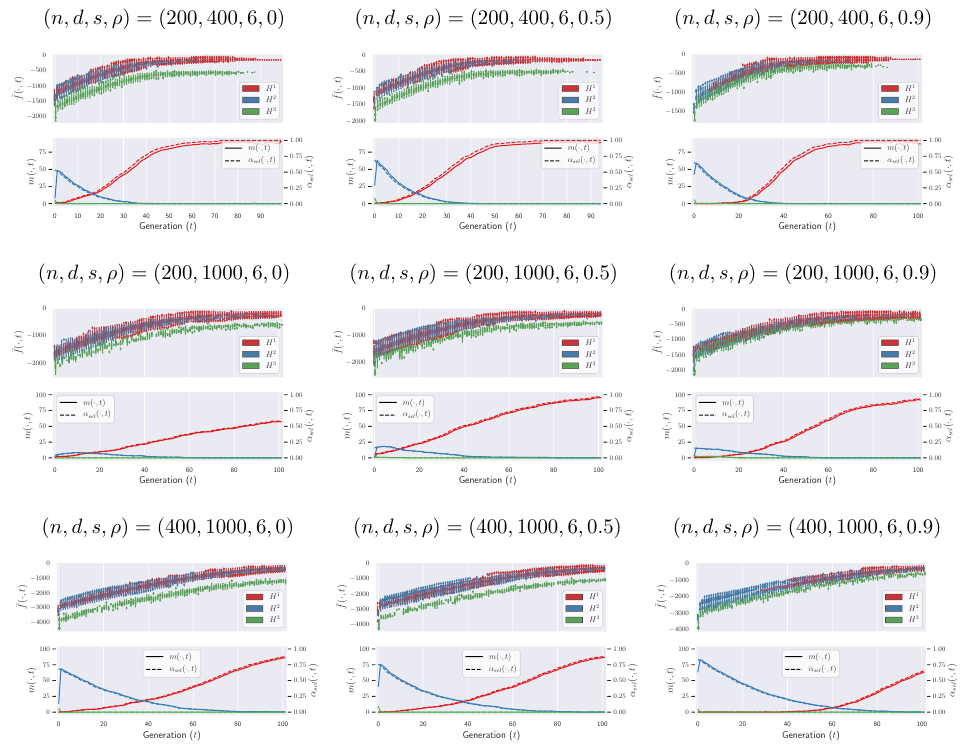}
\caption{\label{fig:schema:case3}
(Case 3) Schema performance (upper panel) and evolution (lower panel) obtained from $500$ simulation runs.}
\end{figure}

\begin{figure}[!tb]
\centering
\includegraphics[width=\textwidth,page=4]{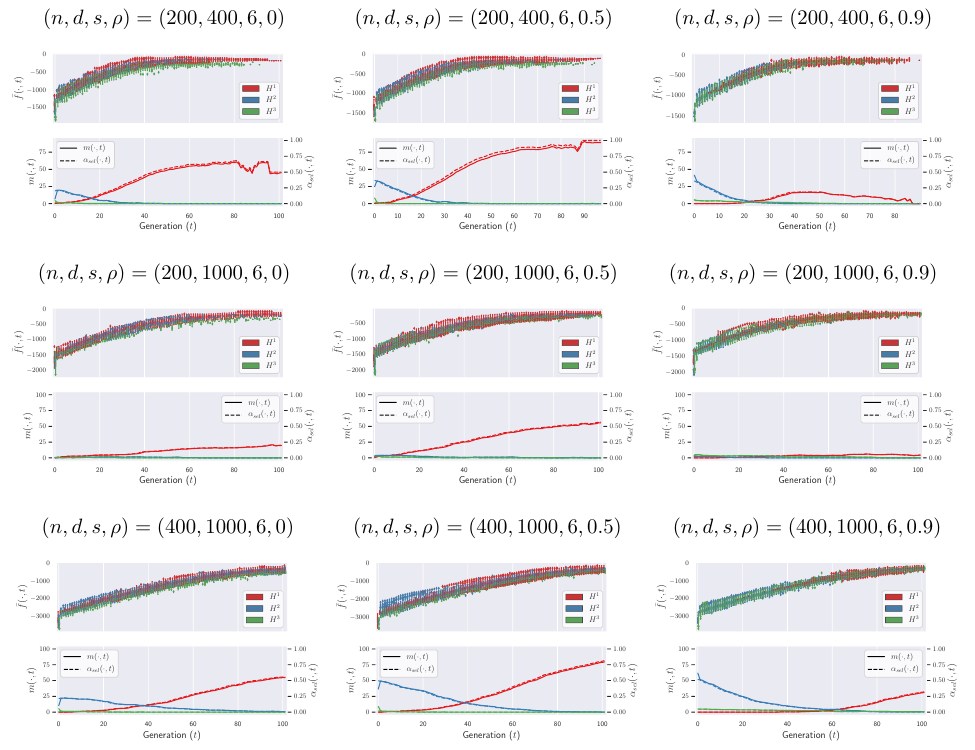}
\caption{\label{fig:schema:case4}
(Case 4) Schema performance (upper panel) and evolution (lower panel) obtained from $500$ simulation runs.}
\end{figure}

\begin{figure}[!tb]
\centering
\includegraphics[width=\textwidth,page=5]{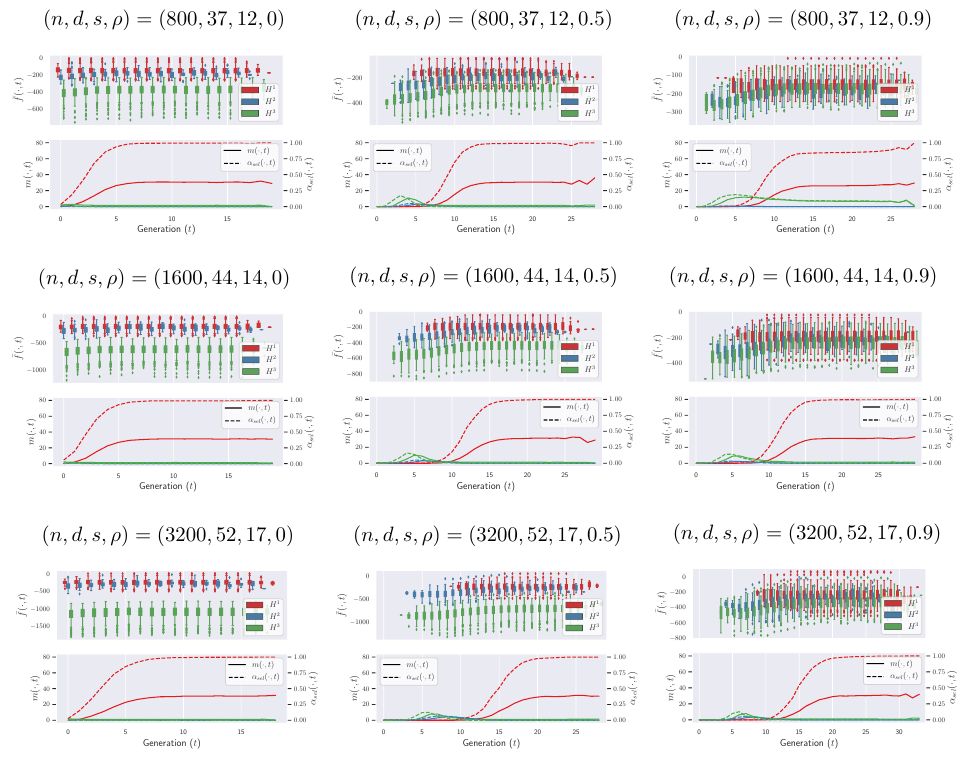}
\caption{\label{fig:schema:case5}
(Case 5) Schema performance (upper panel) and evolution (lower panel) obtained from $500$ simulation runs.}
\end{figure}

\begin{figure}[!tb]
\centering
\includegraphics[width=\textwidth,page=6]{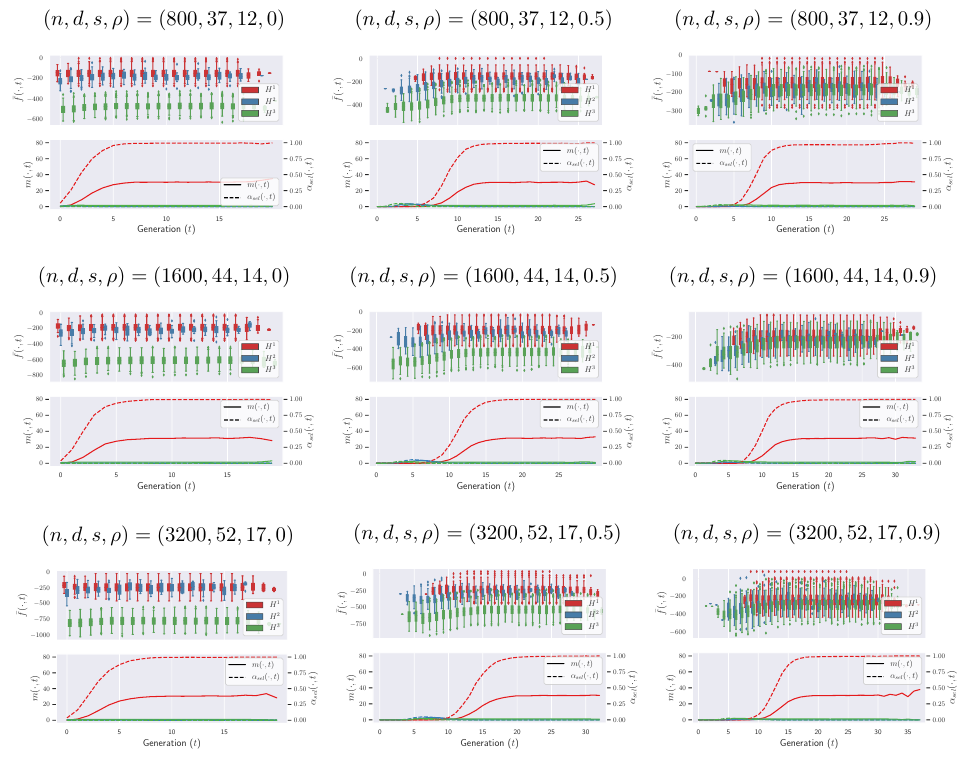}
\caption{\label{fig:schema:case6}
(Case 6) Schema performance (upper panel) and evolution (lower panel) obtained from $500$ simulation runs.}
\end{figure}

Figure~\ref{fig:schema:case1}--\ref{fig:schema:case6} present the additional results of schema evolution.
The conclusions we draw in Section~\ref{sec:sim:schema} still applies for these results, even though the patterns for high-dimensional (Cases 1--4) and low-dimensional (Cases 5 and 6) results are clearly different.

\subsection{Variable Importance}
\label{supp:sim:varimp}

\begin{figure}[!tb]
\centering
\includegraphics[width=\textwidth]{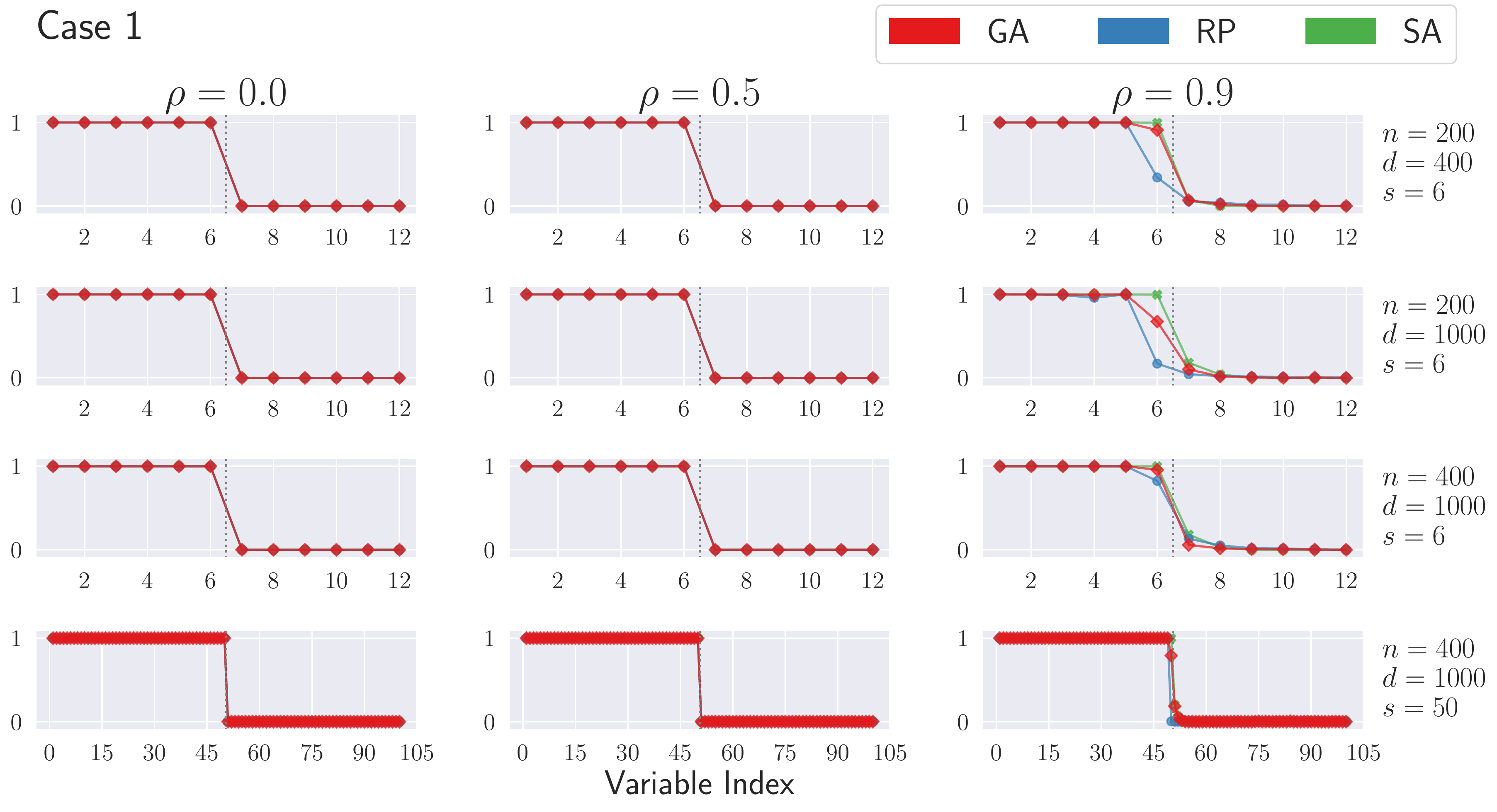}
\caption{\label{fig:soil:case1}
(Case 1)
Averaged SOIL measures.}
\end{figure}

\begin{figure}[!tb]
\centering
\includegraphics[width=\textwidth]{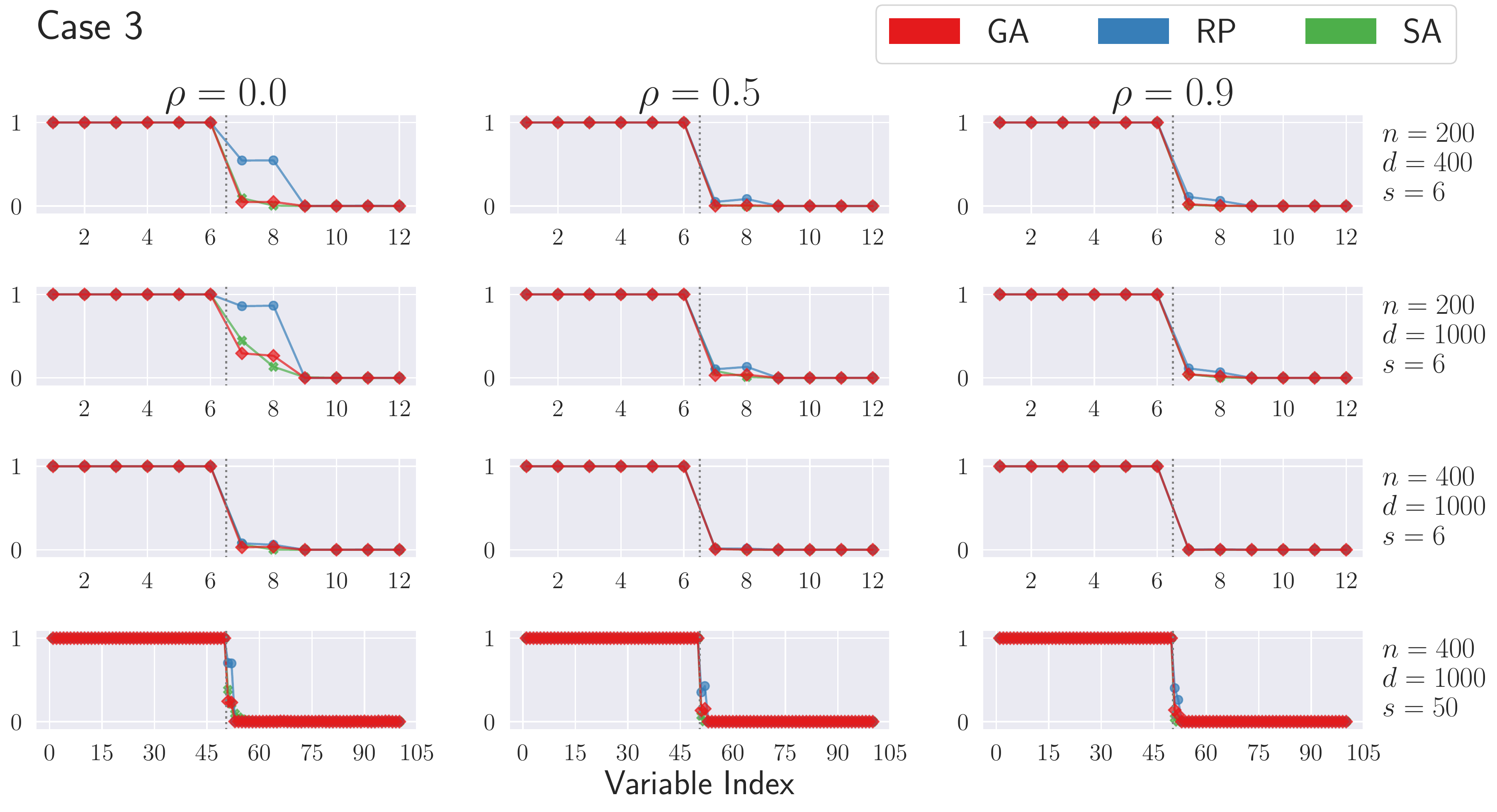}
\caption{\label{fig:soil:case3}
(Case 3)
Averaged SOIL measures.}
\end{figure}

\begin{figure}[!tb]
\centering
\includegraphics[width=\textwidth]{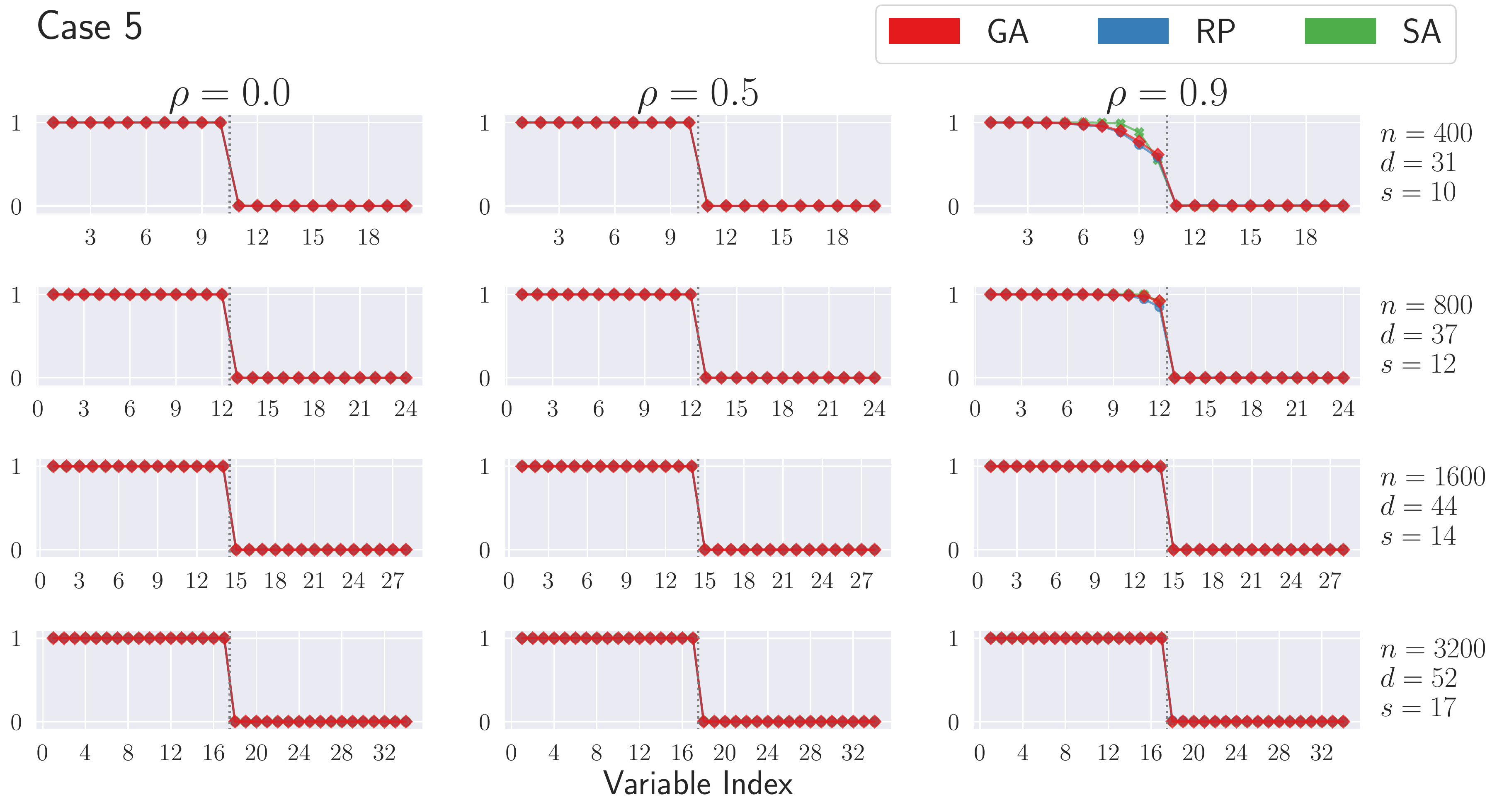}
\caption{\label{fig:soil:case5}
(Case 5)
Averaged SOIL measures.}
\end{figure}

\begin{figure}[!tb]
\centering
\includegraphics[width=\textwidth]{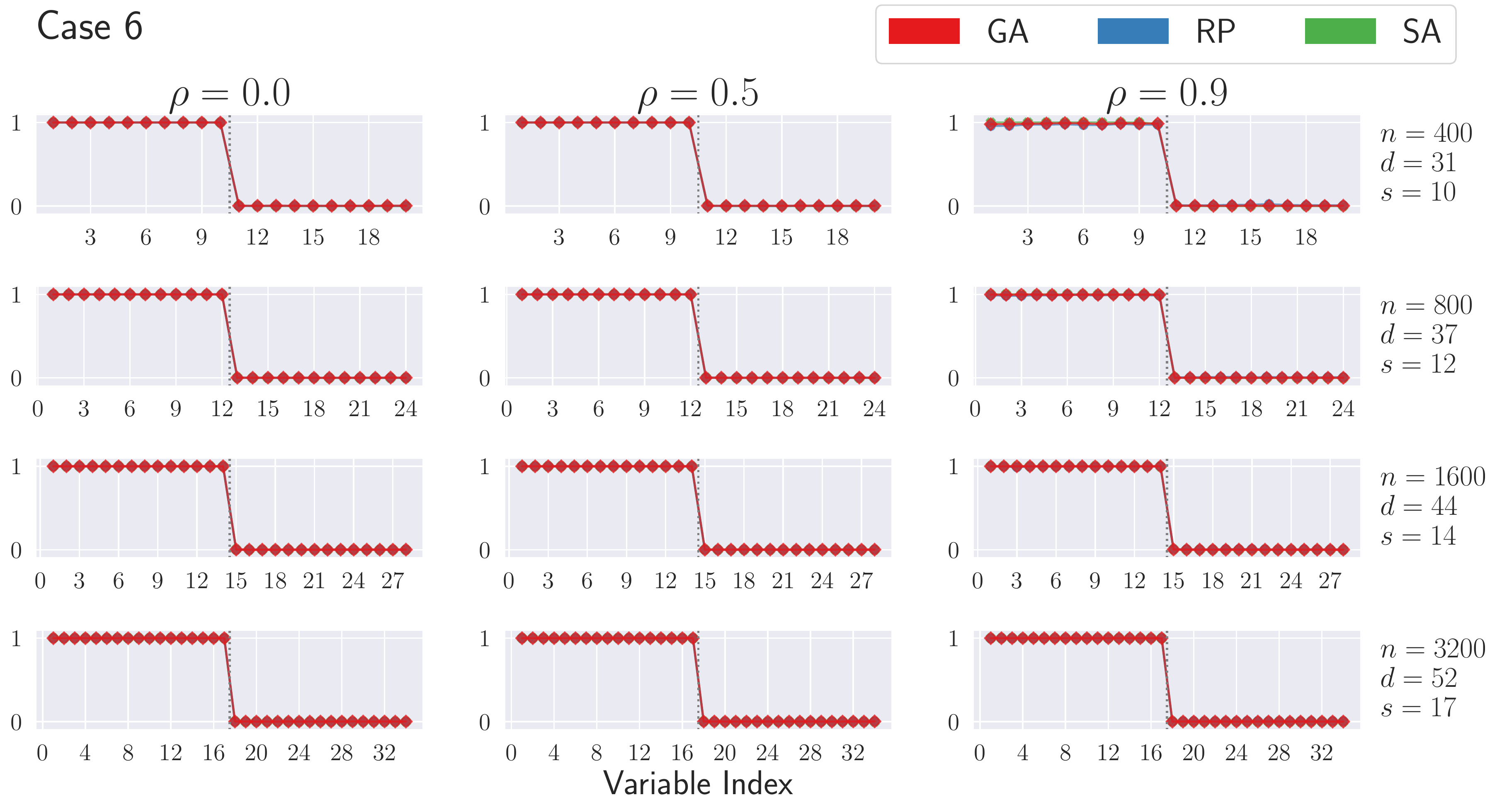}
\caption{\label{fig:soil:case6}
(Case 6)
Averaged SOIL measures.}
\end{figure}

Figure~\ref{fig:soil:case1}--\ref{fig:soil:case6} present additional simulation results of variable importance under Cases 1, 3, 5 and 6.
In Cases 1 and 3 (Figure~\ref{fig:soil:case1} and Figure~\ref{fig:soil:case3}), we see the results of GA and the SA are comparable, and slightly better than the RP in separating the true signals from the rest in some cases (e.g., $X_{s}$ under Case 1 with $\rho=0.9$ and $X_{s+1}$ and $X_{s+2}$ under Case 3 with $\rho=0$).
However, under Cases 5 and 6 (Figure~\ref{fig:soil:case5} and Figure~\ref{fig:soil:case6}) the three methods are just comparable.

\section{Variable Coding for the Residential Building Dataset}
\label{supp:building}

The variable coding with descriptions and units for the residential building dataset used in Section~\ref{sec:realdata:building} is listed in Table~\ref{tab:building:coding}.
Detailed explanations are omitted and can be found in Table~1 of \citet{RA18}.

\begin{table}[!tb]
\centering\small
\caption{\label{tab:building:coding}
Variable coding for the residential building dataset.}
\begin{tabular}{lp{0.7\textwidth}c}
\hline\hline
Variable ID & Descriptions & Unit \\
\hline
\noalign{\vskip 0.5ex}
\multicolumn{3}{l}{\underline{\em Project Physical and Financial (PF) Variables}} \\
\noalign{\vskip 0.5ex}
PF-1 & Project locality defined in terms of zip codes & N/A \\
PF-2 & Total floor area of the building & $m^{2}$ \\
PF-3 & Lot area & $m^{2}$ \\
PF-4 & Total preliminary estimated construction cost based on the prices at the beginning of the project & $10^{7}~\text{IRR}^{m}$ \\
PF-5 & Preliminary estimated construction cost based on the prices at the beginning of the project & $10^{5}~\text{IRR}^{m}$ \\
PF-6 & Equivalent preliminary estimated construction cost based on the prices at the beginning of the project in a selected base year & $10^{5}~\text{IRR}^{m}$ \\
PF-7 & Duration of construction & Time resolution \\
PF-8 & Price of the unit at the beginning of the project per $m^{2}$ & $10^{5}~\text{IRR}^{m}$ \\
\hline
\noalign{\vskip 0.5ex}
\multicolumn{3}{l}{\underline{\em Economic Variables and Indexes (EVI)}} \\
\noalign{\vskip 0.5ex}
EVI-01 & The number of building permits issued & N/A \\
EVI-02 & Building services index (BSI) for a preselected base year & N/A \\
EVI-03 & Wholesale price index (WPI) of building materials for the base year & N/A \\
EVI-04 & Total floor areas of building permits issued by the city/municipality & $m^{2}$ \\
EVI-05 & Cumulative liquidity & $10^{7}~\text{IRR}^{m}$ \\
EVI-06 & Private sector investment in new buildings & $10^{7}~\text{IRR}^{m}$ \\
EVI-07 & Land price index for the base year & $10^{7}~\text{IRR}^{m}$ \\
EVI-08 & The number of loans extended by banks in a time resolution & N/A \\
EVI-09 & The amount of loans extended by banks in a time resolution & $10^{7}~\text{IRR}^{m}$ \\
EVI-10 & The interest rate for loan in a time resolution & ${\%}$ \\
EVI-11 & The average construction cost of buildings by private sector at the time of completion of construction & $10^{5}~\text{IRR}^{m}/m^{2}$ \\
EVI-12 & The average of construction cost of buildings by private sector at the beginning of the construction & $10^{5}~\text{IRR}^{m}/m^{2}$ \\
EVI-13 & Official exchange rate with respect to dollars & $\text{IRR}^{m}$ \\
EVI-14 & Nonofficial (street market) exchange rate with respect to dollars & $\text{IRR}^{m}$ \\
EVI-15 & Consumer price index (CPI) in the base year & N/A \\
EVI-16 & CPI of housing, water, fuel and power in the base year & N/A \\
EVI-17 & Stock market index & N/A \\
EVI-18 & Population of the city & N/A \\
EVI-19 & Gold price per ounce & $\text{IRR}^{m}$ \\
\hline\hline
\end{tabular}
\end{table}

\section{Auxiliary Lemmas}
\label{supp:lemmas}

In this section we provide technical lemmas, with a bit abuse of notations.

\begin{lemma}[Theorem~2 of \citet{R94}]\label{lemma:R94:thm2}
Let $\bP$ be a $n\times n$ reducible stochastic matrix that can be decomposed into
\begin{align*}
\bP=\begin{bmatrix}\bC&\bO\\\bR&\bT\end{bmatrix},
\end{align*}
where $\bC$ is an $m\times m$ primitive stochastic matrix with $m\leq n$ and $\bR$ and $\bT$ are two non-zero matrices with suitable dimensions.
Then there exists an $(n-m)\times n$ positive matrix $\bR_{\infty}$ such that
\begin{align*}
\bP^{\infty}=\lim_{k\to\infty}\bP^{k}=\lim_{k\to\infty}\begin{bmatrix}\bC^{k}&\bO\\\sum_{i=0}^{k-1}\bT^{i}\bR\bC^{k-i}&\bT^{k}\end{bmatrix}=\begin{bmatrix}\bC^{\infty}&\bO\\\bR_{\infty}&\bO\end{bmatrix}
\end{align*}
is a stable stochastic matrix with $\bP^{\infty}=\bm{1}\bpi^{\top}$, where $\bm{1}=(1,\ldots,1)^{\top}$ is the vector of $1$'s with suitable length,  $\bpi=(\pi_{1},\ldots,\pi_{n})^{\top}=\bm{\pi}_{0}^{\top}\bP^{\infty}$ is unique regardless of the initial distribution $\bm{\pi}_{0}$, and $\bpi$ satisfies $\pi_{i}>0$ for $i=1,\ldots,m$ and $\pi_{i}=0$ for $i=m+1,\ldots,n$.
\end{lemma}

\end{document}